\documentclass[11pt]{article}
\usepackage{graphicx,epstopdf,subfigure,mathtools,mathrsfs, amsmath, amssymb} 
\usepackage[font=small,labelfont=bf]{caption}
\usepackage{float, cancel, amssymb}
\usepackage[normalem]{ulem}
\usepackage[title]{appendix}
\PassOptionsToPackage{usenames,dvipsnames}{xcolor}
\usepackage[usenames,dvipsnames]{xcolor}
 
\usepackage[margin=1 in]{geometry}

\usepackage{amsfonts, soul, authblk, enumitem}

\newcommand{\new}[1]{\color{black}#1\normalcolor}

\newcommand{\widebar}[1]{%
   \hbox{%
     \vbox{%
       \hrule height 0.5pt 
       \kern0.5ex
       \hbox{%
         \kern-0.1em
         \ensuremath{#1}%
         \kern-0.1em
       }%
     }%
   }%
} 
\def\onedot{$\mathsurround0pt\ldotp$}
\def\ddot{
  \mathbin{\vcenter{\baselineskip.67ex
    \hbox{\onedot}\hbox{\onedot}}%
  }}%

\graphicspath{{../../BrownianTex/}}
\newfloat{algorithm}{!t}{loa}
\providecommand{\algorithmname}{Algorithm}
\floatname{algorithm}{\protect\algorithmname}

\newcommand{\V}[1]{\boldsymbol{#1}}                 
\newcommand{\M}[1]{\boldsymbol{#1}}

\newcommand{\Lop}[1]{\boldsymbol{\mathcal{#1}}}

\newcommand{\ind}[2]{{#1}^{(#2)}}

\newcommand{\tdisc}[2]{#1^{\left(#2\right)}}
\renewcommand{\tt}[1]{#1}
\newcommand{\Exp}[1]{\mathbb{E}\left[{#1}\right]}                 

\global\long\def\D#1{\Delta#1} 
 
\global\long\def\dt#1{\partial_t#1} 
\global\long\def\ds#1{\partial_s#1}

\global\long\def\norm#1{\left\Vert #1\right\Vert }

\global\long\def\Xhat{\widehat{\V{X}}}

\global\long\def\div{\partial}

\global\long\def\lp{\ell_p}
\global\long\def\Bind#1{{\{#1\}}}
\global\long\def\eqdi{\stackrel{d}{=}}
\global\long\def\Wt{\widetilde{\M{W}}}

\global\long\def\kon{k_\text{on}}

\global\long\def\Slet#1{\M{S}\left(#1\right)}
\global\long\def\Dlet#1{\M{D}\left(#1\right)}

\global\long\def\Xs{\V{\tau}}
\global\long\def\XPoly{\mathbb{X}}
\global\long\def\TauPoly{\mathbb{T}}
\global\long\def\Xmp{\V{X}_\text{MP}}
\global\long\def\Ump{\V{U}_\text{MP}}
\global\long\def\X{\M{\Lop{X}}}
\global\long\def\EPMI{\frac{1}{8\pi\mu}}

\global\long\def\FCL{\V{F}^{\text{(CL)}}}

\global\long\def\kon{k_\text{on}}
\global\long\def\konb{k_\text{on,s}}

\global\long\def\eps{{\hat{a}}}
\global\long\def\epsRS{\hat{\epsilon}}
\global\long\def\rc{a}
\global\long\def\epsc{\epsilon}

\global\long\def\Mfor{\widetilde{\M{M}}}
\global\long\def\MforRef{\widetilde{\M{M}}_\text{ref}}
\global\long\def\Msym{\Mfor}

\global\long\def\Xsbar{\bar{\Xs}}
\global\long\def\Cbar{\bar{\M C}}
\global\long\def\Wproc{\V{\mathcal{W}}}
\global\long\def\gauss{\V{\eta}}

\usepackage{hyperref}
\hypersetup{
    colorlinks=false,
    pdfborder={0 0 0},
}

\title{Bending fluctuations in semiflexible, inextensible, slender filaments in Stokes flow: towards a spectral discretization}
\author[1]{Ondrej Maxian}
\author[2]{Brennan Sprinkle}
\author[1]{Aleksandar Donev}
\affil[1]{Department of Mathematics, Courant Institute, New York University, New York, NY 10012}
\affil[2]{Department of Applied Mathematics and Statistics, Colorado School of Mines, Golden, CO 80401}

\begin{document}
\maketitle

\begin{abstract}
Semiflexible slender filaments are ubiquitous in nature and cell biology, including in the cytoskeleton, where reorganization of actin filaments allows the cell to move and divide. Most methods for simulating semiflexible inextensible fibers/polymers are based on discrete (bead-link or blob-link) models, which become prohibitively expensive in the slender limit when hydrodynamics is accounted for. In this paper, we develop a novel coarse-grained approach for simulating fluctuating slender filaments with hydrodynamic interactions. Our approach is tailored to relatively stiff fibers whose persistence length is comparable to or larger than their length, and is based on three major contributions. First, we discretize the filament centerline using a coarse non-uniform Chebyshev grid, on which we formulate a \emph{discrete} constrained Gibbs-Boltzmann (GB) equilibrium distribution and overdamped Langevin equation for the evolution of the unit-length tangent vectors. Second, we define the hydrodynamic mobility at each point on the filament as an integral of the Rotne-Prager-Yamakawa kernel along the centerline, and apply a spectrally-accurate ``slender-body'' quadrature to accurately resolve the hydrodynamics. Third, we propose a novel midpoint temporal integrator which can correctly capture the Ito drift terms that arise in the overdamped Langevin equation. For two separate examples, we verify that the equilibrium distribution for the Chebyshev grid is a good approximation of the blob-link one, and that our temporal integrator for overdamped Langevin dynamics samples the equilibrium GB distribution for sufficiently small time step sizes. We also study the dynamics of relaxation of an initially straight filament, and find that as few as 12 Chebyshev nodes provides a good approximation to the dynamics while allowing a time step size two orders of magnitude larger than a resolved blob-link simulation. We conclude by applying our approach to a suspension of cross-linked semiflexible fibers (neglecting hydrodynamic interactions between fibers), where we study how semiflexible fluctuations affect bundling dynamics. We find that semiflexible filaments bundle faster than rigid filaments even when the persistence length is large, but show that semiflexible \emph{bending} fluctuations only further accelerate agglomeration when the persistence length and fiber length are of the same order.
\end{abstract}

\section{Introduction}
The closer we look at biological systems, the more we find slender filaments performing important work. These filaments are responsible for cell motility and division in prokaryotes \cite{berg2004coli, lauga2009hydrodynamics} and eukaryotes \cite{mogilner1996cell, alberts, vavylonis2008assembly, murrell2015forcing}, which makes them indispensable for processes like wound healing, stem cell differentiation, and organism growth. In cells, three kinds of filaments can be distinguished: microtubules, which are sufficiently stiff as to behave deterministically \cite{brangwynne2007force}, intermediate filaments, which have small persistence length and are consequently found in entangled networks in vivo \cite{noding2012intermediate, pawelzyk2014attractive, charrier2016mechanical}, and actin filaments, which have visible bending fluctuations driven by thermal forces, but have sufficient stiffness to maintain a relatively smooth appearance of the filament centerline \cite{gittes1993flexural}.

We are motivated here by actin filaments, which have been shown to take on a range of morphologies when combined with cross-linking proteins in vivo \cite{welch1997actin}, in vitro \cite{falzone2012assembly}, and in silico \cite{freedman2019mechanical, ActinCLsRheology}. Our interest in particular is how thermal fluctuations, hydrodynamic interactions, and cross linking compete or cooperate with each other to determine the steady state morphology and stress-strain behavior of actin networks. While there have been a number of coarse-grained theories examining this question \cite{chen2021nonlinear, chen2022motor, mulla2019origin}, there is still a need for detailed simulation of each of the microscopic components in the system, so that assumptions made in deriving coarse-grained theories can be validated and tested more rigorously. This was the idea behind our previous work on \emph{deterministic} actin filaments \cite{ActinCLsRheology}, which looked at how each of the system components contributes to the behavior of the cross-linked actin network on short and long timescales. Our ultimate goal is to extend the study of \cite{ActinCLsRheology} to \emph{fluctuating} (Brownian) actin filaments, so that we can determine how the thermal fluctuations affect the viscoelastic gel behavior. 

With this goal in mind, we turn here to the simulation of semiflexible filaments interacting through a viscous medium. There is already a large body of literature on this topic, which can be analyzed by asking the following two questions: is the chain extensible (with a penalty for stretching) or constrained to be exactly inextensible? And, is the chain being simulated discrete or continuous? That is, is the discrete chain considered the truth itself, or is it simply a discretization of a continuum equation which converges in the limit as the grid spacing goes to zero? This second question is quite difficult, so we begin our review of the literature with discrete chains. 

Discrete extensible chains (known in the literature as bead-spring models) are the most commonly-used models for fluctuating filaments because of their simplicity \cite{barbier2009numerical, kim2009computational, groot2013mesoscale, vargas2018fiber, pincus2022dilute}. In this case, the chain is represented by a series of beads or blobs connected by springs which penalize but do not prohibit extensibility and bending. There are therefore no constraints, and it is relatively straightforward to simulate the overdamped Langevin dynamics using standard algorithms \cite[Sec.~3]{cruz2012review}. Actin filaments, which are nearly inextensible and therefore require a high stretching modulus to accurately simulate, push the boundaries of bead-spring models, as the high stretching and bending moduli restrict the maximum possible time step size, making long simulations of actin filament systems prohibitive \new{(see \cite{krishna2022petascale} for the limit on modern GPUs)}.  

It is therefore attractive to constrain the distance between each bead, which leads to a second class of methods known as bead-link or blob-link models \cite[Sec.~4]{cruz2012review} (we will use the lesser-used ``blob-link'' terminology because of the coupling with ``blob-based'' hydrodynamics \cite{RigidMultiblobs}). This approach, which in theory enables larger time step sizes and longer simulations, is difficult in practice because it requires the formulation of an overdamped Langevin equation with the constraint that the chain is inextensible. In previous work, this has been done by defining a constrained Langevin equation for the \emph{positions} of the blobs which has complicated stochastic drift terms and requires highly specialized algorithms to simulate \cite{morse2003theory, ottinger2012stochastic}. An alternative view, which is one we adopt here, is to view the degrees of freedom as the \emph{link} orientations, as well as the position of the fiber center. For an inextensible chain, this reduces the problem to a series of connected rigid rods, for which we can make use of previous work by some of us on Langevin equations for rigid bodies \cite{BrownianMultiBlobs}.

The main issue with the blob-link model in our context arises when we account for the hydrodynamics of the chain. Denoting the position of bead $i$ by $\V{X}_\Bind{i}$, hydrodynamic interactions give us a mobility matrix $\Mfor \left(\V{X}\right)$ which describes the relationship between the forces on the beads $\V{F}$ and their velocities $\V{U}$ via $\V{U}=\Mfor \V{F}$. For blob-link or bead-spring models in unbounded domains, the obvious choice for the mobility matrix is to use a pairwise hydrodynamic kernel between the blobs/beads. One such kernel is the Rotne-Prager-Yamakawa (RPY) tensor \cite{rpyOG, wajnryb2013generalization}, which is obtained by approximately solving the mobility problem for a pair of spheres of radius $\eps$ an arbitrary distance apart in the fluid. That is, the $3 \times 3$ mobility block associated with beads $i$ and $j$ in a blob-link model is
\begin{equation}
\label{eq:BigRPYMat}
\Mfor_{\Bind{i}\Bind{j}} = \Mfor_\text{RPY}\left(\V{X}_\Bind{i}, \V{X}_\Bind{j}; \eps \right),
\end{equation}
where the $3 \times 3$ matrix $\Mfor_{RPY}\left(\V{x},\V{y}; \eps\right)$ is defined in\ \eqref{eq:rpyknel}. While this kernel and others like it neglect interactions involving more than two beads, it is commonly used in polymer physics to describe the hydrodynamics of a chain \cite{butler2005brownian, beck2006ergodicity, wang2012flipping, keavRPY} because the RPY kernel is symmetric positive definite for any locations of the two particles. This fact can also be seen by using an immersed boundary formulation to justify\ \eqref{eq:BigRPYMat} \cite[Sec.~2.1]{TwistBend}. One contribution of this paper will be to develop a temporal integrator for a fluctuating blob-link chain with mobility given by\ \eqref{eq:BigRPYMat}.

As we will show, however, the blob-link model breaks down in the limit of very slender fibers, such as actin filaments, which have a radius of $\rc=4$ nm \cite{grazi1997diameter} and can have lengths $L$ of a few microns in vivo and tens of microns in vitro  \cite{janmey1994mechanical, mcgrath2000regulation}. This set of parameters gives aspect ratios on the order of $\epsc=\rc/L=10^{-3}$, making the hydrodynamics multiscale in nature. By using a blob-link model with blobs of radius $\sim \rc$, we are implicitly marrying ourselves to the resolution of this smallest lengthscale. \new{If we wanted to accurately simulate very flexible filaments or filaments that are nearly touching, then there would be no way around this. But for smoother filament shapes and non-dense suspensions where we interested in looking at the bulk behavior, a heavy price is paid:} since the hydrodynamics we are really interested in occur on lengthscales $L=\rc/\epsc$, the number of required blobs to resolve the hydrodynamics of a slender filaments is approximately $L/\rc \sim 1/\epsc$ \cite{ttbring08, kallemov2016immersed}. In addition to being prohibitively expensive in the spatial discretization, this model also presents a problem for temporal integration, since having blobs spaced a lengthscale $\rc$ apart often requires us to temporally resolve the fluctuations on that lengthscale; i.e., the time step size is constrained by the \emph{smallest} lengthscale in the system. Thus, for slender filaments a blob-link discretization is prohibitively expensive in both the spatial and temporal variables using existing methods.

\begin{figure}
\centering
\includegraphics[width=\textwidth]{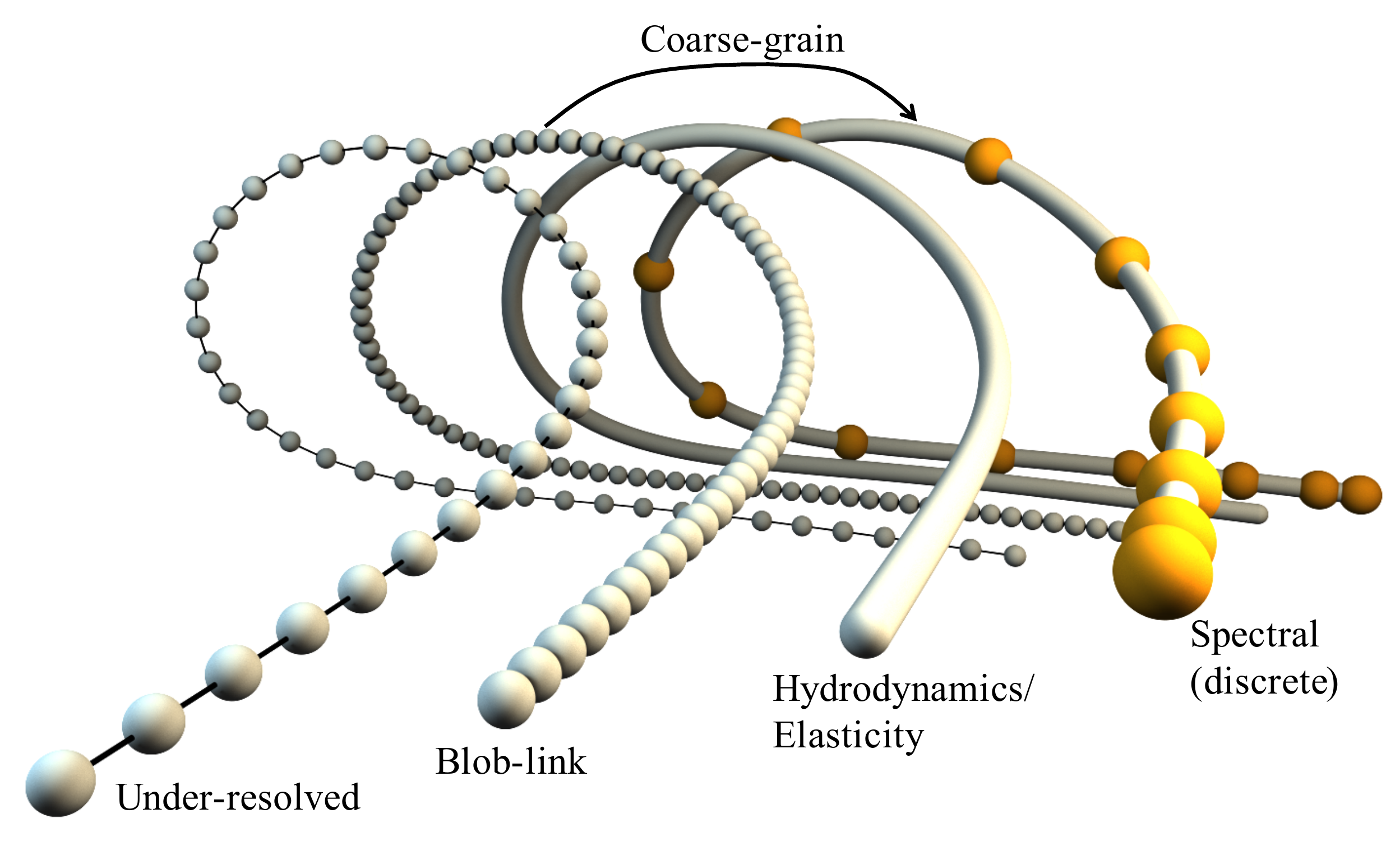}
\caption{\label{fig:Fibs}Conceptual picture of our filament model. We begin with a blob-link chain at left, for which we need a large number of beads (second from left) to resolve the hydrodynamic interactions. A sensible model then becomes a continuum filament (second from right), shown as a solid curve, which we then replace by a fully discrete model based on Chebyshev collocation points (the gold nodes at right show 16 such points). In our approach, the model on the right is a coarse-grained approximation of the fully-resolved blob-link (second-from-left) model with the continuum model (second-from-right) used only in elasticity and hydrodynamic calculations. Figure\ \ref{fig:SchDisc} has a more detailed look at this discretization. }
\end{figure}

Since the hydrodynamics forces such a large number of degrees of freedom per fiber, a sensible resolution in the slender limit is to take a continuum limit of the discrete model, so that the positions $\V{X}$ of the nodes become a function $\XPoly(s)$ of a continuous variable $s$ which describes the arclength of a curve (\cite{shelley1996nonlocal, shelley2000stokesian}, see Fig.\ \ref{fig:Fibs}). Likewise, the links become a continuous function $\TauPoly(s)=\ds{\XPoly}(s)$, and the force on each link becomes a force density $\V{f}(s)$ that is defined everywhere along the filament. The hydrodynamic mobility $\Lop{M}\left[\XPoly(\cdot)\right]$ that maps force density to velocity can then be defined as the integral operator
\begin{equation}
\label{eq:McRPYdef}
\V{U}(s)= \int_0^L \M{M}_\text{RPY}\left(\XPoly(s), \XPoly(s'); \eps \right) \V{f}(s') \, ds':=\left(\Lop{M}\V{f}\right)(s).
\end{equation}
In \cite[Appendix~A]{TwistBend}, we showed that the RPY integral\ \eqref{eq:McRPYdef} is asymptotically equivalent to slender body theory, which is the more commonly-used mobility for continuum curves in Stokes flow \cite{krub,johnson,ts04,ehssan17}, if the relationship between the true filament radius $\rc$ and the regularization radius $\eps$ is given by
\begin{equation}
\eps = \frac{e^{3/2}}{4}\rc \approx 1.1204 \rc.
\end{equation}
This statement is only true for the translation-translation component of the RPY kernel; see  \cite{TwistSBT} for how this choice of radius performs in the context of rotation, as well as \cite{cortez2012slender} for the corresponding analysis for regularized Stokeslets \cite{cortez2018regularized}, and \cite{ttbring08} for a numerical analysis for the immersed boundary method \cite{peskin2002acta}.

To simulate this continuum model, we still ultimately need to choose collocation points at which to represent the fiber positions and tangent vectors. Because we are interested in actin filaments, for which the fiber shapes are relatively smooth, an attractive discretization is a spectral one, where the collocation points are chosen from a Chebyshev grid in the arclength $s$. The motivation for using Chebyshev points comes from classical numerical analysis; since Chebyshev polynomials give spectrally-accurate polynomial interpolants, the error in approximating smooth fiber shapes decreases exponentially in the number of nodes, similar to using a Fourier basis for a closed loop fiber \cite[c.~10]{ascher2011first}. In addition, the well-conditioned property of the Chebyshev interpolation matrix allows us to obtain globally-accurate interpolants to the fiber shape, which are simpler than breaking the fiber into panels. In previous work on deterministic filaments \cite{FibersWeakInextensibility, TwistBend}, we developed numerical methods for such global Chebyshev discretizations \cite{FibersWeakInextensibility} and accompanying quadrature schemes \cite[Appendix~G]{TwistBend} for the integral\ \eqref{eq:McRPYdef} on the spectral grid. Using these schemes, we showed that the number of collocation points required to resolve the integral\ \eqref{eq:McRPYdef} is (roughly) independent of the fiber aspect ratio $\epsc$ \cite[Sec.~4.4]{TwistBend}, which makes this mobility definition, and the spectral grid that accompanies it, an appealing choice for \emph{smooth} slender filaments. In this paper, we will conduct the same analysis for \emph{fluctuating} (Brownian) semiflexible filaments, showing that there are indeed significant advantages to using this model of the mobility over the discrete model, at least for very slender filaments.

The problem with taking the continuum limit in the presence of thermal fluctuations is that it leads to ill-posed constrained Langevin partial-differential equations when dynamics are taken into account. From equilibrium statistical mechanics, it is known that for a free inextensible worm-like fiber the tangent vector $\Xs(s)$ performs Brownian motion in $s$ on the unit sphere with diffusion coefficient $\ell_p^{-1}$, where the persistence length $\lp=\kappa/k_B T$ is defined as the ratio of the bending stiffness $\kappa$ to the thermal energy $k_B T$. On the other hand, we show in Appendix\ \ref{sec:Drift} that the overdamped Ito Langevin equation for a blob-link model with dynamics has stochastic drift terms which are required to obtain the correct equilibrium statistics (see Fig.\ \ref{fig:PoorMan}, where we simulate the Langevin equation with and without the drift terms). These drift terms, which are a consequence of \emph{both} the inextensibility constraint and the hydrodynamics, are only well-defined with respect to the \emph{discrete} orientations of the tangent vectors, and do not converge in the continuum limit. This issue has either not been been treated in previous continuum methods for filament Brownian dynamics \cite{saint1, saint2}, or been partially sidestepped by making the filament extensible through a penalty force \cite{liu2019efficient,ryan2022finite}. The latter approach, which has characterized a number of finite element methods for biopolymers \cite{cyron2009finite, lin2014combined,cyron2012numerical, cyron2013micromechanical} takes us back to an extensible chain (and its aforementioned temporal stiffness), and leaves unclear how to make sense of the continuum limit for constrained fluctuating filaments. And if there is no continuum limit, how is a spectral discretization possible?

While the continuum limit is an interesting object from the standpoint of applied stochastic analysis, for the purposes of simulation it has little relevance, since it requires an infinite number of degrees of freedom. The object that we simulate, like any other statistical mechanics model, must ultimately be discrete, and so we will propose a \emph{fully discrete} or coarse-grained model of a fluctuating fiber based on a spectral representation of the chain (see Fig.\ \ref{fig:Fibs}). The important conceptual leap is then to think of the ``continuum limit'' only as a way to efficiently approximate the ``true'' hydrodynamics of $1/\epsc$ blobs without actually needing to track that many degrees of freedom. As such, for a given value of $\eps$ and $L$ we define our reference result to be a discrete blob-link chain with blobs spaced $\eps$ apart ($L/\eps=1/\epsRS$ blobs). Our hope is then to approximate, to reasonable ($\approx 10\%$) accuracy, the equilibrium statistical mechanics  and dynamics of the blob-link chain using a coarse-grained Chebyshev discretization with as few nodes as possible. The goal of this, which we show is achieved at least in part, is to extract the advantages of the blob-link and continuum discretizations into a single numerical method: by proposing a fully discrete ``spectral'' chain, we can write a well-defined overdamped Langevin equation which samples from a well-defined \emph{discrete} equilibrium Gibbs-Boltzmann distribution. But by using the mobility\ \eqref{eq:McRPYdef} to define the velocity at each point on the fiber, we avoid  having to scale the number of collocation points with the fiber aspect ratio.

Seeing that this is, to our knowledge, the first paper to consider a spectral discretization of constrained Brownian hydrodynamics, our program is very much experimental, in part because our premise is counterintuitive. The conventional wisdom, as taught in a first semester of numerical analysis, is that spectral methods perform well for \emph{smooth} problems (with exponentially-decaying spectra), while for nonsmooth ones (power law spectra) a sparse finite difference method is a better choice. While it is certainly true that Brownian motion induces a power law spectrum in the fiber positions, our hope is that restricting to \emph{semiflexible} fibers (with $\ell_p/L \gtrsim 1$) will make a spectral discretization tractable by reducing the number of modes required for accurate simulation. In other words, as long as the Chebyshev grid resolves the fiber \emph{persistence length} (as opposed to radius), many quantities should be approximated well since scales smaller than $\ell_p$ will not contribute. We demonstrate that this is indeed the case through examples, but it should be noted that this is just a first step, with the conclusion of this paper laying out much of the work that remains.

\section{Langevin equation for semiflexible filaments \label{sec:LangTint}}
In this section, we formulate the overdamped Langevin equation for the evolution of a discrete inextensible semiflexible filament. As discussed in the introduction, we propose a fully discrete spectral discretization composed of tangent vectors $\Xs$, which are defined on a Chebyshev grid of size $N$, and node positions $\V{X}$, which are defined on a Chebyshev grid of size $N+1$. Once we formulate the discrete filament model, we can write the Gibbs-Boltzmann probability distribution for the fiber energy. The overdamped Langevin equation then follows once we define the dynamics of how the fibers evolve deterministically.

\subsection{Discrete fiber model \label{sec:KinInex}}
Let us first formulate a spectral discretization of an inextensible filament which is well-suited for fluctuating hydrodynamics. As discussed in the introduction and shown in Fig.\ \ref{fig:SchDisc}, we choose to track a discrete collection of tangent vectors $\Xs$, which evolve as rigid rods by rotating on the unit sphere. To completely define the fiber, we need to also track the fiber midpoint $\Xmp$. These two quantities give a set of node positions $\V{X}$, which in turn define a smooth polynomial interpolant $\XPoly(s)$ for the fiber centerline. We use this interpolant to compute elastic energy and, importantly, hydrodynamic interactions via\ \eqref{eq:McRPYdef}. Once we define a discretization, we can postulate the constrained Gibbs-Boltzmann distribution in terms of the tangent vectors $\Xs$. This distribution depends on the fiber elastic energy $\mathcal{E}_\text{bend}$, which we discretize in the final part of this section.

\subsubsection{Discretization}
\begin{figure}
\centering

\includegraphics[width=0.9\textwidth]{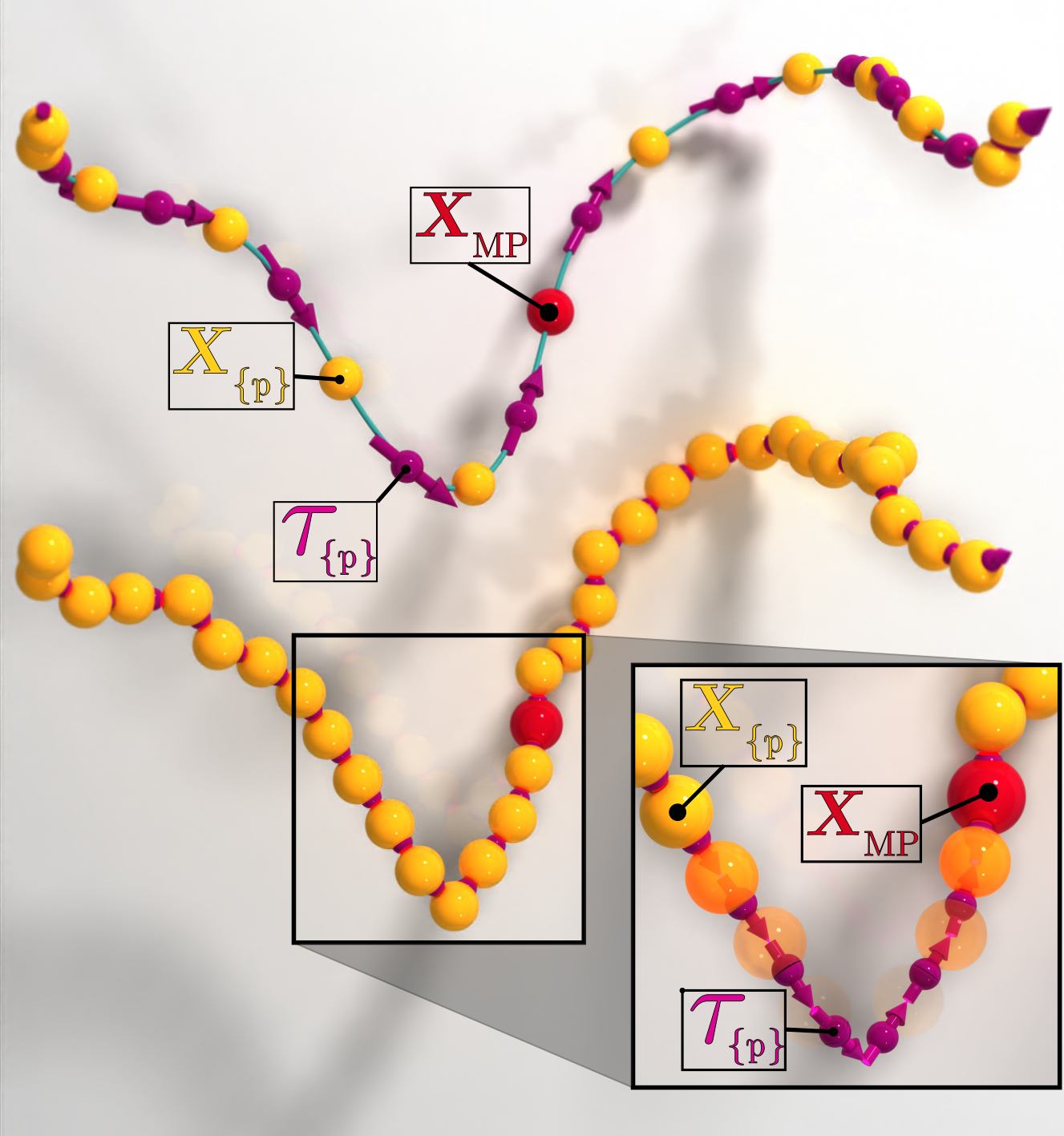}
\caption{\label{fig:SchDisc}Transferring the fiber discretization for the blob-link model (bottom) to a spectral grid (top). We show the tangent vector locations as magenta vectors visible through the ``x-ray'' inset -- for the blob-link discretization these are the points exactly between the beads, while for the spectral discretization they are defined on a type 1 Chebyshev grid with $N$ points ($N=12$ here). The positions, which can be obtained from $\Xs$ and $\Xmp$ (red sphere), are shown as yellow spheres -- for the blob-link discretization these are evenly spaced points, while for the spectral discretization they are defined on a type 2 Chebyshev grid with $N_x=N+1$ points. Note the clustering of both $\Xs$ and $\V{X}$ near the fiber endpoints in the spectral model, which is a characteristic of the Chebyshev grid. This means that the tangent vectors no longer connect individual beads, as they do for the blob-link discretization.}

\end{figure}

We define a spectral discretization of an inextensible filament by transferring the blob-link discretization concept onto a non-uniform Chebyshev grid (see Fig.\ \ref{fig:SchDisc}). We begin with a collection of $N$ unit-length tangent vectors $\Xs=\{ \Xs_{\Bind{p}}\}_{p=1}^N$, with $\Xs_{\Bind{p}}\cdot \Xs_{\Bind{p}}=1$ for $p=1, \dots N$, on a type 1 Chebyshev grid (i.e., a Chebyshev grid that does not include the endpoints).\footnote{Here and throughout this paper, $\Xs$ refers to a $3N \times 1$ column-stacked vector of each $3 \times 1$ tangent vector $\Xs_{\Bind{p}}$. The (scalar) $p$th entry of this vector is denoted $\Xs_p$.} To obtain the fiber position, we oversample $\Xs$ to a (type 2, endpoints-included) Chebyshev grid of size $N+1$ using the evaluation matrix $\M{E}_{N \rightarrow N+1}$, then integrate the result exactly using the Chebyshev integration matrix $\M{D}^\dagger_{N+1}$ (pseudo-inverse of the differentiation matrix $\M{D}_{N+1}$). This can be written as
\begin{gather}
\label{eq:XNp1}
\V{X} = \begin{pmatrix}  \M{D}^\dagger_{N+1} \M{E}_{N \rightarrow N+1} & \M{B} \end{pmatrix} \begin{pmatrix} \Xs \\ \Xmp \end{pmatrix}:=\X \Xsbar.
\end{gather}

The matrix $\M{B}$ is such that the midpoint of $\V{X}$ on the grid of size $N+1$ is $\Xmp$ (when $N$ is even, $\Xmp$ is an actual point on the $N+1$ size Chebyshev grid, but $N$ could also be odd, for which the midpoint is obtained via interpolation). Thus, the nodal points that track the fiber position $\V{X}$ are defined on a grid of size $N+1$, whereas $\Xs$ is defined on a grid of size $N$, just as in the blob-link discretization. To go in the reverse direction, we apply the inverse of $\X$, which is done by differentiating $\V{X}$ on the $N+1$ point grid and then downsampling to the grid of size $N$ via the matrix $\M{E}_{N+1 \rightarrow N}$. The midpoint is determined from the $N+1$ nodes via sampling the $N+1$-term interpolating polynomial at the midpoint using the matrix $\M{E}_{N+1\rightarrow \text{MP}}$. Together this gives, 
\begin{gather}
\label{eq:Xinv}
 \X^{-1} =\begin{pmatrix} \M{E}_{N+1 \rightarrow N}  \M{D}_{N+1} \\ \M{E}_{N+1\rightarrow \text{MP}} \end{pmatrix},
\end{gather}
which is the actual inverse of $\X$ because we carefully handle the integration of $N$-term Chebyshev polynomials by using a grid of size $N+1$. We contrast this with handling everything on a grid of size $N$ as we have done previously \cite{FibersWeakInextensibility, TwistBend}, in which case information is lost when converting the unit-length tangent vectors to the positions $\V{X}$. For the deterministic examples we studied in \cite{FibersWeakInextensibility, TwistBend}, the fibers were relatively smooth, and so the lost (high-frequency) information had a negligible effect. This is no longer the case for Brownian filaments, and so the method here is required to correctly track the high-frequency modes. It is important to note, however, that this formulation only applies to open, two-ended, filaments, and \emph{not} looped filaments, which require further conditions on $\Xs$. 

While we motivated\ \eqref{eq:XNp1} using a spectral discretization, these equations also hold for the standard blob-link discretization \cite{keavRPY}; in that case the map $\X^{-1}$ is defined by taking finite differences of nodal points to give values at the links, and the product $\X$ is defined by summing the values of the links to obtain values at the nodes. Thus, the equations we write in this section are general once the maps $\X$ and $\X^{-1}$ have been specified.

\subsubsection{Gibbs-Boltzmann distribution}
Now that we have introduced our discrete fiber model, we can write the Gibbs-Boltzmann equilibrium distribution as a function of the chain degrees of freedom $\Xsbar$. Letting $\mathcal{E}_\text{bend}[\Xsbar]$ denote the discrete bending energy of the fiber, which we discretize in Section\ \ref{sec:elforce}, we take the Gibbs-Boltzmann distribution to be
\begin{align}
\label{eq:GBDist}
dP_\text{eq}\left(\Xsbar\right) &= Z^{-1} \exp{\left(-\mathcal{E}_\text{bend}(\Xsbar)/ k_B T\right)} d\mu_0\left(\Xsbar\right)\\ \nonumber
P_\text{eq}\left(\Xsbar\right)&= Z^{-1} \exp{\left(-\mathcal{E}_\text{bend}(\Xsbar)/ k_B T\right)} \prod_{p=1}^N \delta \left(\Xs_\Bind{p}^T \Xs_\Bind{p} - 1\right).
\end{align}
What we mean by the product of $\delta$ functions is really that the base measure $d\mu_0\left(\Xsbar\right)$ in the Gibbs-Boltzmann distribution\ \eqref{eq:GBDist} corresponds to the tangent vectors $\Xs_\Bind{p}$ being independently uniformly distributed on the unit sphere for $p=1, \dots N$. For blob-link chains, the tangent vectors $\Xs_\Bind{p}$ are uniformly spaced, as shown in Fig.\ \ref{fig:SchDisc}, and this distribution follows naturally from their inextensibility and statistical independence in a freely-jointed chain, in which case\ \eqref{eq:GBDist} defines the standard worm-like chain model. 

In our spectral discretization, the tangent vectors $\Xs_\Bind{p}$ are defined on a Chebyshev grid, and this model lacks the physical motivation that characterizes the blob-link model (see Fig.\ \ref{fig:SchDisc}). In fact, the base measure $d\mu_0$ on a uniform (blob-link) grid has no continuum limit, since a freely-jointed continuum chain would require the tangent vector $\Xs(s)$ to perform infinitely fast Brownian motion in $s$ on the unit sphere; thus only\ \eqref{eq:GBDist} with the Gibbs-Boltzmann weight makes sense in continuum. Furthermore, a Chebyshev \emph{polynomial} cannot represent an inextensible curve everywhere (since the number of zeros in $\Xs \cdot \Xs-1$ is limited by the polynomial degree), and therefore it is not obvious how to uniquely define a discrete equilibrium distribution for a spectral grid that is a good approximation to the continuum worm-like chain (when the Chebyshev grid resolves the persistence length)\footnote{Of course, we could also take $d\mu_0$ as the base measure and put additional ``entropic'' or ``metric'' factors in $\mathcal{E}_\text{bend}$. Here, however, we take $\mathcal{E}_\text{bend}$ to be the standard elastic bending energy of a semiflexible chain.}. In the first step toward a more mathematically-justified approach, we posulate\ \eqref{eq:GBDist} as a reasonable guess, but do not claim any precise sense of convergence of\ \eqref{eq:GBDist} for a spectral grid to the equilibrium distribution of a continuum worm-like chain. We do show (in Section\ \ref{sec:EEMCMC} and Appendix\ \ref{sec:CurvX0}) that for a sufficient number of Chebyshev nodes over the persistence length, the equilibrium distribution for a spectral chain approximates well the equilibrium distribution for the blob-link chain, in the sense that samples from the two distributions give the same large-scale statistical properties of the chain, e.g., its end-to-end distance. Forthcoming work will justify\ \eqref{eq:GBDist} for spectral chains through the theory of coarse graining \cite{pepCG}.

\subsubsection{Bending energy and forces \label{sec:elforce}}
We now turn to the evaluation of the bending energy $\mathcal{E}_\text{bend}$, and the resulting force and force density that it generates on the fiber centerline for free fibers. In continuum, the fiber bending or curvature energy is given by the squared $L^2$ norm of the fiber curvature vector, 
\begin{equation}
\label{eq:Econt}
\mathcal{E}_\text{bend}\left[\XPoly\left(\cdot\right)\right] = \frac{\kappa}{2}\int_0^L \ds^2 \XPoly(s) \cdot \ds^2 \XPoly(s) \, ds.
\end{equation}
Now, there are two competing ideas for how to generate a force (or force density) from this energy. In the continuum perspective we have used in the past \cite{FibersWeakInextensibility, TwistBend}, we take a functional derivative of $\mathcal{E}_\text{bend}$ in continuum to yield a force density (inside the integral) $\V{f}(s) = - \kappa \ds^4 \XPoly$, subject to the boundary conditions $\ds^2 \XPoly(0,L)=\ds^3\XPoly(0,L)=\V 0$. In a spectral method, we discretize this force density using rectangular spectral collocation (RSC) \cite{tref17, dhale15}, which does \emph{not} guarantee that the total force and torque are zero on a nonsmooth fiber, since the force obtained does not come from differentiating a discrete energy functional \cite{peskin2002acta}. As a result, we cannot write down the equilibrium distribution we would expect the Brownian dynamics to obey, even for extensible fibers.

Because of this, it is useful to discretize the energy functional\ \eqref{eq:Econt} directly, then use the matrix that results to compute forces, implementing the boundary conditions naturally. On a spectral grid, we can put the energy functional\ \eqref{eq:Econt} in the form of an energy \emph{function} $\mathcal{E}_\text{bend}\left(\V{X}\right)=\frac{1}{2}\V{X}^T \M{L}\V{X}$, where
\begin{equation}
\label{eq:Lmat}
\M{L}=\kappa \left(\M{E}_{N_x \rightarrow 2N_x}\M{D}^2\right)^T \M{W}_{2N_x} \left(\M{E}_{N_x \rightarrow 2N_x} \M{D}^2\right):=\kappa \left(\M{D}^2\right)^T \Wt\M{D}^2.
\end{equation}
is a matrix that takes two derivatives of $\V{X}$ on a grid of size $N_x=N+1$, then does an inner product of those derivatives on a $2N_x$ grid using Clenshaw-Curtis weights matrix $\M{W}_{2N_x}$. \new{The inner product on the $2N_x$ grid is computed with an upsampled representation $\M{E}_{N_x \rightarrow 2N_x}\V{X}$, where the upsampling matrix $\M{E}_{N_x \rightarrow 2N_x}$ is applied by using the points $\V{X}$ to form the Chebyshev interpolant $\XPoly(s)$, and then evaluating the interpolant on the upsampled grid.} The purpose of this is to compute the inner product of the Chebyshev polynomials representing $\M{D}^2 \V{X}$ \emph{exactly} (this could be done on a slightly smaller grid because the two derivatives make the polynomial representation $\M{D}^2 \V{X}$ have degree at most $N_x-3$). The $L^2$ inner product on the upsampled grid corresponds to an inner product on $\mathbb{R}^{N_x}$ with the weights matrix 
\begin{equation}
\label{eq:Wt}
\Wt=\M{E}_{N_x \rightarrow 2N_x}^T\M{W}_{2N_x}\M{E}_{N_x \rightarrow 2N_x}.
\end{equation}
Note that an equivalent way to formulate the energy is to use $\Xsbar$ as the degrees of freedom, so that
\begin{equation}
\label{eq:Udef}
\mathcal{E}_\text{bend}\left(\Xsbar\right)=\frac{1}{2}\Xsbar^T \X^T \M{L} \X \Xsbar:=\frac{1}{2}\Xsbar^T \M{L}_{\Xs} \Xsbar.
\end{equation}
The \emph{force} on the fiber nodes is obtained by differentiating the energy function,
\begin{equation}
\label{eq:force}
\V{F} = -\frac{\partial \mathcal{E}_\text{bend}}{\partial \V{X}} = -\M{L}\V{X}=-\M{L}\X \Xsbar.
\end{equation}

When we compute the fiber hydrodynamics, we will need to input the force \emph{density} on the centerline. We can obtain this at the nodes from the force in\ \eqref{eq:force} by 
\begin{equation}
\label{eq:FDenFromF}
\V{f} = \Wt^{-1} \V{F},
\end{equation}
and use the Chebyshev interpolant $\V{f}(s)$ as a continuum force density. Note that the matrix $\Wt$, rather than the diagonal matrix $\M{W}$ of weights on the $N_x$ point Chebyshev grid, must be used for the force density to converge as the spatial discretization is refined. The reason for this is that the weights matrix\ \eqref{eq:Wt} enters in the discrete $L^2$ inner product, and the force density is the representation of the derivative of energy with respect to that inner product \cite[Sec.~6]{li2017towards},\footnote{A more concrete justification for this is as follows: suppose that the filament satisfied the free fiber boundary conditions $\ds^2 \XPoly(0,L)=\ds^3 \XPoly(0,L)=\V 0 $. Then in continuum we could equivalently write the energy as $U=\kappa/2 \int_0^L \XPoly(s) \cdot \ds^4\XPoly(s) \, ds$, which we would discretize as $U=\kappa/2 \V{X}^T \Wt \M{D}^4 \V{X}$. Taking the derivative, we have the force $\V{F}=-\kappa \Wt \M{D}^4 \V{X}$. But we know that force density should be $\V{f}=-\kappa \M{D}^4 \V{X}$, which can only be obtained from $\V{F}$ using\ \eqref{eq:FDenFromF}.} i.e., the function that satisfies $\langle \V{X},\V{f} \rangle = \V{X}^T \V{F}$. 

\subsection{Dynamics}
If we substitute the energy\ \eqref{eq:Udef} into the Gibbs-Boltzmann distribution\ \eqref{eq:GBDist}, we see that our goal is to write an overdamped Langevin equation that is in detailed balance with respect to the distribution
\begin{align}
\label{eq:GBDist2}
P_\text{eq}\left(\Xsbar\right) = Z^{-1} \exp{\left(-\Xsbar^T \M{L}_{\Xs} \Xsbar/\left(2 k_B T\right)\right)} \prod_{p=1}^N \delta \left(\Xs_\Bind{p}^T \Xs_\Bind{p} - 1\right),
\end{align}
and includes hydrodynamic interactions between points on the filament. To do this, we first need to discuss how inextensible filaments evolve \emph{deterministically}. This begins with a description of the kinematics of discrete inextensible filaments, which are analogous to that of a series of elastically-interacting rigid rods. We then describe the evaluation of hydrodynamic interactions and formulate the equations of motion in the deterministic setting. The overdamped Langevin equation can be obtained from these deterministic considerations and\ \eqref{eq:GBDist2} by following a standard formulation.

\subsubsection{Kinematics \label{sec:InexKin}}
Let us first consider the evolution of the fiber tangent vectors. Since, for any $p$, $\Xs_\Bind{p} \cdot \Xs_\Bind{p}=1$ for all time, it follows that the evolution of the tangent vectors can be described by 
\begin{equation}
\label{eq:DtTau}
\dt{\Xs} = \V{\Omega} \times \Xs:=-\M{C}\left[\Xs\right]\V{\Omega},
\end{equation}
where the matrix $\M{C}\left[\Xs\right]$ is such that $\M{C}\left[\Xs\right] \V{\Omega} = \Xs \times \V{\Omega}$; we will drop the explicit $\Xs$ dependence when clear from context. The matrix $\M{C}$ satisfies the following properties, which will be useful when we formulate the overdamped Langevin equation for the fiber evolution:
\begin{gather}
\label{eq:CXs}
\M{C} \Xs = -\M{C}^T \Xs = \V 0, \qquad
\div_{\Xs} \cdot \M{C}^T =\V 0.
\end{gather}
The first equation\ \eqref{eq:CXs} follows from the definition of $\M{C}$ as a cross product with $\Xs$, while the second, which is a divergence with $\Xs$, has $j$th entry
\begin{equation}
\label{eq:DefDiv}
\left(\div_{\Xs} \cdot \M{C}^T\right)_j:=\frac{\partial}{\partial \Xs_k} \M{C}^T_{jk} = \frac{\partial}{\partial \Xs_k} \M{C}_{kj}=0
\end{equation}
because the $k$th row of $\M{C}$ has no entries that depend on $\Xs_k$; here and throughout this paper, repeated indices are summed over using Einstein's convention. Based on these two properties, we can conclude that the time evolution of the tangent vector\ \eqref{eq:DtTau} is analogous to the time evolution of the unit quaternion describing a rigid body's orientation \cite[Eqs.~(5--8)]{BrownianMultiBlobs}, which can be represented as a unit vector on the unit 4 sphere.

The evolution of the tangent vectors according to\ \eqref{eq:DtTau} automatically implies that the evolution of the positions $\V{X}$ is given by 
\begin{equation}
\label{eq:KNp1}
\dt{\V{X}}=\X \begin{pmatrix}-\M{C} & \M{0} \\ \M{0} & \M{I} \end{pmatrix} \begin{pmatrix} \V{\Omega}\\ \Ump \end{pmatrix}:=\X \Cbar\V{\alpha}:=\M{K} \V{\alpha}.
\end{equation}
Note the analogy with\ \eqref{eq:XNp1}, but here we transform \emph{velocity} instead of position using a square $3(N+1) \times 3(N+1)$ non-invertible matrix $\M{K}$. We define a pseudo-inverse of $\M{K}$ as
\begin{equation}
\label{eq:Kinv}
\M{K}^{-1}=\begin{pmatrix} \M{C} & \M{0} \\ \M{0} & \M{I} \end{pmatrix}  \X^{-1}=\Cbar^T\X^{-1}.
\end{equation}
Note that the matrices $\M{K}$ and $\M{K}^{-1}$ have rank $2N+3$, since the $N$ tangent vectors $\Xs$ live in the null space\footnote{This is true even for straight fibers. In continuum, the only nontrivial solutions of $\left(\int_{L/2}^s \V{\Omega}(s') \, ds'\right) \times \Xs+\Ump=\V 0$ are those with $\V{\Omega}=\Omega^\parallel \Xs$. Our discretization preserves this property; since in our method $\int \V{\Omega}(s) \, ds $ (and, for non-straight fibers $\int \V{\Omega}(s) \times \Xs(s) \, ds $) is a polynomial with $N+1$ terms (degree $N$), it can only be exactly equal to $-\Ump$ at $N$ points (and not $N_x=N+1$). This shows why careful handling of high-frequency modes (on a grid of size $N+1$) is necessary.} of $\M{C}$. Thus $\M{K}^{-1}$ is not a true inverse of $\M{K}$, but rather
\begin{equation}
\label{eq:KinvK}
\M{K}^{-1}\M{K} = \begin{pmatrix} -\M{C}^2 & \M{0} \\ \M{0} & \M{I} \end{pmatrix} \rightarrow \Cbar\M{K}^{-1}\M{K} =\Cbar,
\end{equation}
i.e., $\M{K}^{-1}\M{K}$ acts like the identity when applied to $\Cbar$ from the right. \new{The last equality holds because the matrix $\M{C}^2$ is block diagonal with $p$th block $\M{I}_3-\Xs_\Bind{p}\Xs_\Bind{p}$, and when multiplied by $\M{C}$ the term involving $\Xs_\Bind{p}$ becomes zero.} What\ \eqref{eq:KinvK} also means in practice is that 
\begin{equation}
\M{K}^{-1}\M{K} \V{\alpha}=\M{K}^{-1}\M{K} \begin{pmatrix} \V{\Omega} \\ \Ump \end{pmatrix} = \begin{pmatrix} \V{\Omega}^\perp \\ \Ump \end{pmatrix},
\end{equation}
where $\V{\Omega}^\perp_\Bind{p}=\left(\M{I}-\Xs_\Bind{p}\Xs_\Bind{p}\right)\V{\Omega}_\Bind{p}$ \new{(once again, here we see $\M{K}^{-1}\M{K}$ projecting off the parallel parts of $\V{\Omega}$, since the $p$th block of $\M{K}^{-1}\M{K}$ is $\M{I}_3-\Xs_\Bind{p}\Xs_\Bind{p}$)}. Thus, because the fiber evolves according to\ \eqref{eq:DtTau}, it makes no difference whether we use $\V{\alpha}$ or $\M{K}^{-1}\M{K} \V{\alpha}$ for the fiber velocities. It is for this reason that we use this discretization for $\M{K}$, as opposed to those from our prior work \cite{FibersWeakInextensibility, TwistBend}, which viewed $\V{\Omega}$ as a continuous function of $s$, and were consequently plagued by aliasing errors in trying to recover $\V{\Omega}$ from $\M{K}\V{\alpha}$ \cite[Eq.~(110)]{FibersWeakInextensibility}.

While the evolution\ \eqref{eq:DtTau} describes how the tangent vectors evolve in continuous time, in discrete time the update $\tdisc{\Xs}{n+1}=\tdisc{\Xs}{n}-\D t \M{C}\V{\Omega}$ does \emph{not} preserve the unit-length constraint. Thus, in order to keep the dynamics on the constraint in discrete time, we will solve for an angular velocity $\V{\Omega}$, then evolve each tangent vector $\tdisc{\Xs}{n}_\Bind{p}$ by rotating by $\V{\Omega}_\Bind{p}$, 
\begin{equation}
\label{eq:RotUpdate}
\tdisc{\Xs}{n+1}_\Bind{p}=\text{rotate}\left(\tdisc{\Xs}{n}_\Bind{p}, \V{\Omega}_\Bind{p} \D t\right),
\end{equation}
the explicit formula for which is given in \cite[Eq.~(111)]{FibersWeakInextensibility}. We then update the midpoint $\Xmp$ via $\tdisc{\Xmp}{n+1}=\tdisc{\Xmp}{n}+\D t \tdisc{\Ump}{n}$. In deterministic methods \cite{ts04, FibersWeakInextensibility}, rotating the tangent vectors in this way keeps the dynamics on the constraint without needing to introduce penalty parameters \cite{ts04} that lead to additional stiffness in temporal integration.

\subsubsection{Hydrodynamic mobility \label{sec:Mob}}
We now discuss the discretization of the mobility matrix $\Mfor$ that relates the force on the nodal points to their velocities. In this paper, we will not consider hydrodynamic interactions between multiple filaments. Thus the matrix $\Mfor$ can be formed separately as a dense matrix on each fiber, and any operations with $\Mfor$, such as inversions and square roots, done directly by storing the eigenvalue (or Cholesky if $\Mfor \succcurlyeq \M{0}$) decomposition. Extensions to multiple filaments interacting hydrodynamically, in which case $\Mfor$ cannot be formed densely, will be covered in future work.

In \cite{TwistBend}, we developed an efficient quadrature scheme for\ \eqref{eq:McRPYdef}, but this mobility acts on force \emph{densities} by applying a matrix $\M{M}$ to obtain the relationship $\V{U}=\M{M}\V{f}$ (we use the notation $\Mfor$ for the matrix mapping force to velocity and $\M{M}$ for the matrix mapping force density to velocity). Thus, to use our quadrature scheme, we first need to first convert force to force density using\ \eqref{eq:FDenFromF}, which gives $\V{U}=\M{M} \Wt^{-1} \V{F}$. If we were to compute the quadratures \emph{exactly}, for instance by upsampling to a fine grid, $\M{M}\Wt^{-1}$ should be a symmetric positive definite matrix, since the work dissipated in the fluid
\begin{equation}
\int_0^L \V{f}(s) \cdot \V{U}(s) \, ds= \V{f}^T \Wt \V{U}=\V{f}^T \Wt \M{M}\V{f}
\end{equation}
is always positive, which implies that $\Wt \M{M}$ (and $\Wt^{-1} \Wt \M{M}\Wt^{-1}=\M{M}\Wt^{-1})$ is an SPD matrix. Indeed, upsampling to a fine grid gives us a reference SPD mobility 
\begin{equation}
\label{eq:Mref}
\MforRef = \Wt^{-1} \M{E}_u^T \M{W}_u \Mfor_{\text{RPY}, u} \M{W}_u \M{E}_u \Wt^{-1}.
\end{equation}
where the subscript $u$ denotes a matrix on the fine grid. Walking through the steps of this calculation, we first convert force to force density by applying $\Wt^{-1}$. We then extend the force density to the upsampled grid and integrate it against the RPY kernel there (the matrix $\Mfor_{\text{RPY}, u}$ describes the pairwise RPY kernel on the upsampled grid, and $\M{W}_u$ gives the integration weights). The first three matrices, which follow from the symmetry of $\MforRef$, downsample the velocity $\V{U}_{u}$ on the upsampled grid to the $N_x$ point grid by minimizing the $L^2$ difference between $\M{E}_u \V{U}_{N_x}$ and $\V{U}_{u}$. Note that while we write\ \eqref{eq:Mref} in a collocation perspective, it can be shown that the same fundamental matrix appears in a Galerkin method, where we would expect an automatically SPD matrix.

When we use our efficient quadrature scheme to approximate\ \eqref{eq:Mref}, the matrix that results is $\Mfor=\M{M} \Wt^{-1}$, which is not guaranteed to be symmetric positive definite. Indeed, in \cite[Sec.~4.4.1]{TwistBend}, we showed that the mobility matrix $\M{M}$ obtained from quadrature can have negative eigenvalues due to numerical errors, especially for larger $N$ and $\epsRS$. The negative eigenvalues for large $N$ and $\epsRS$ are not altogether surprising; we know that lengthscales in the forcing on the order $\epsRS$ will be filtered by our RPY regularization, resulting in eigenvalues close to zero. Putting more Chebyshev points, which are clustered together near the boundary, brings these lengthscales into play. Combining this with the imperfect accuracy of our quadrature scheme, which is designed for smooth forces, it is not hard to understand why negative eigenvalues result.  

To work around this problem, we define the mobility as the symmetric matrix
\begin{equation}
\label{eq:Msym}
\Mfor = \frac{1}{2}\left(\M{M}\Wt^{-1}+\Wt^{-1} \M{M}^T\right).
\end{equation}
Then, we compute an eigenvalue decomposition of this matrix and set all eigenvalues less than a threshold $\sigma$ to be equal to $\sigma$. The choice of $\sigma$ for a given discretization comes from the smallest eigenvalue of the reference mobility\ \eqref{eq:Mref}. In our tests, our quadrature-based mobility, which is slightly modified for the case of random filaments as described in Appendix\ \ref{sec:QuadMod}, gives the same dynamics for a relaxing filament as the reference mobility\ \eqref{eq:Mref}; see Appendix\ \ref{sec:QuadOS}. Thus the error made in quadrature, symmetrization, and eigenvalue truncation is small in practice, even though the quadrature is only spectrally accurate for smooth fibers. 

\begin{figure}
\centering
\includegraphics[width=0.6\textwidth]{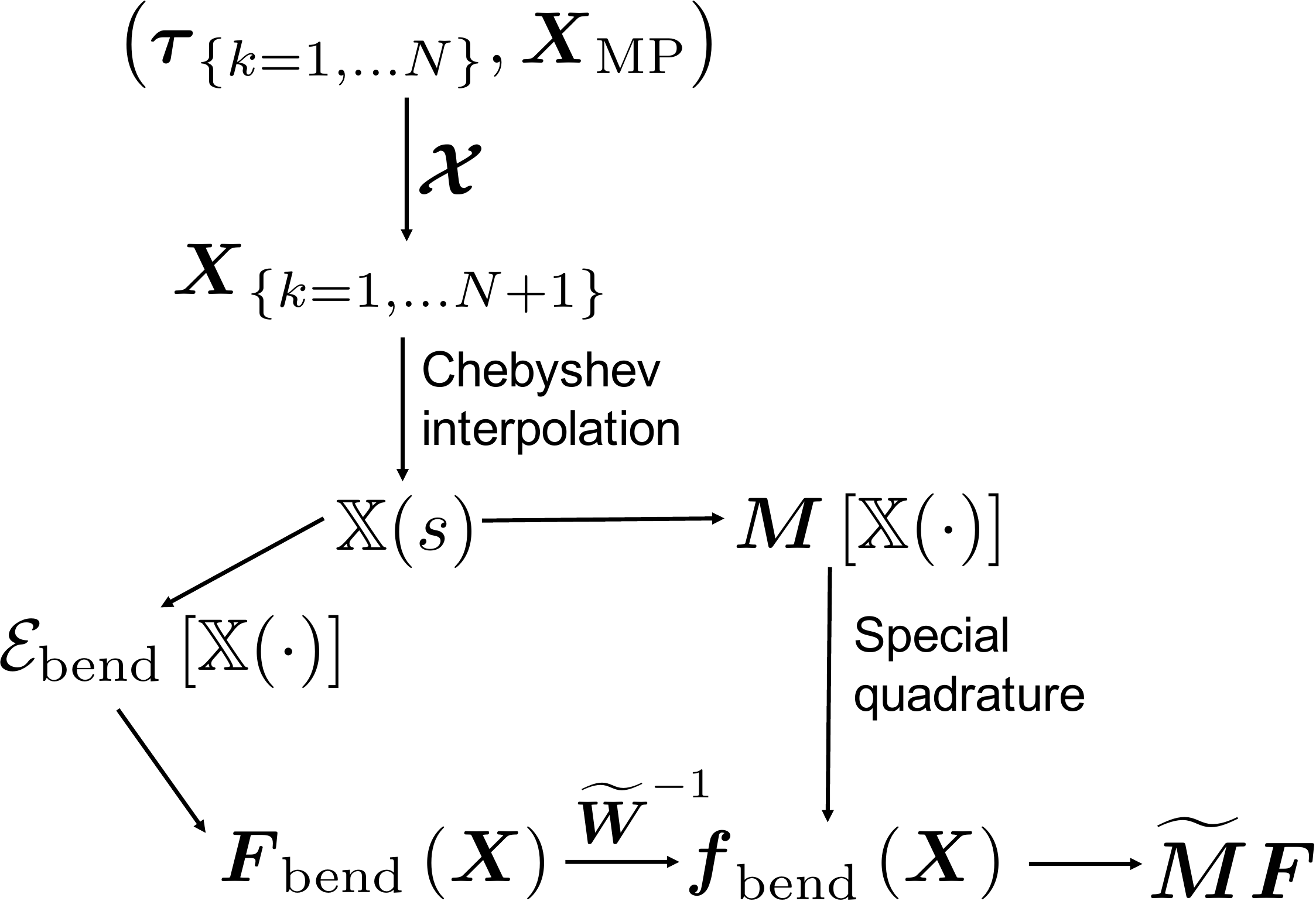}
\caption{\new{\label{fig:Flowchart} Summary of how we use the fiber configuration to compute the forces $\V{F}$ and mobility $\Mfor$. Given $\left(\Xs,\Xmp\right)$, we apply the map $\X$ from\ \eqref{eq:XNp1} to obtain the Chebyshev polynomial $\XPoly(s)$, which we then use to compute the bending energy $\mathcal{E}_\text{bend}$ in\ \eqref{eq:Econt}, force $\V{F}$ in\ \eqref{eq:force}, and force density $\V{f}$ in\ \eqref{eq:FDenFromF}. The velocity $\Mfor \V{F}$ is computed by applying special quadrature to the force density $\V{f}$, see Section\ \ref{sec:Mob} (because the quadrature matrix might have negative eigenvalues, we form $\Mfor$ explicitly to compute and truncate its eigenvalues).}}
\end{figure}

In Appendix\ \ref{sec:AdvSpec}, we also compare the slender body quadrature mobility to two other commonly-used mobilities for slender body hydrodynamics: direct RPY-based quadrature on the Chebyshev grid, and a local slender body theory \cite{saint1, saint2}. Our results show that there are advantages to using the special quadrature, both in terms of resolving spatial scales as $\epsRS \rightarrow 0$, and in terms of the time step required for temporal accuracy. \new{With that in mind, Fig.\ \ref{fig:Flowchart} shows how we use the fiber geometry and special quadrature to apply the mobility $\Mfor \V{F}$.}

\subsubsection{Deterministic evolution \label{sec:DetDyn}}
We are now ready to put the pieces of our discretization together to obtain the deterministic evolution of the fiber centerline. In a deterministic method, the rotation rates and midpoint velocities $\V{\alpha}=\left(\V{\Omega},\Ump\right)^T$ are given by the solution of the saddle point system \cite[Eq.~(54)]{FibersWeakInextensibility}
\begin{gather}
\label{eq:SP}
\M{K}\V{\alpha} = \Mfor \left(-\M{L}\V{X}+\M{\Lambda}\right)\\ \nonumber
\M{K}^T \M{\Lambda} = \V 0,
\end{gather}
where $\M{K}$, $\Mfor$, $\M{L}$ and $\M{K}^T$ are all matrices of size $3N_x=3(N+1)$, and $\V{\Lambda}$ is a constraint force. This is a reformulation of\ \cite[Eq.~(54)]{FibersWeakInextensibility} in terms of \emph{force} rather than force density, as the mobility matrix $\Mfor\left(\V{X}\right)$ relates the velocity of the filament centerline to the \emph{forces} at the nodal points. As shown in Appendix\ \ref{sec:SPLag}, the deterministic dynamics\ \eqref{eq:SP} can also be obtained by minimizing a constrained Lagrangian function, which implies that they are dissipative and have the structure of a gradient descent flow.

Using a Schur complement approach, the Lagrange multipliers $\V{\Lambda}$ can be eliminated from\ \eqref{eq:SP} to yield 
\begin{gather}
\label{eq:dXDet}
\dt{\V{X}} = \M{K}\V{\alpha} = -\left(\M{K}\M{N} \M{K}^T\right) \M{L}\V{X},\\
\label{eq:Ndef}
\text{where} \qquad \M{N}=\left(\M{K}^T \Mfor^{-1} \M{K}\right)^\dagger
\end{gather} 
is the mobility matrix projected onto the space of inextensible motions. We can transform this equation to obtain the evolution of $\Xsbar$ by substituting $\V{X} = \X \Xsbar$ into\ \eqref{eq:dXDet} to obtain
\begin{align}
\label{eq:dTau}
\dt{\Xsbar}&=\X^{-1} \left(\M{K}\M{N}\M{K}^T\right)\M{L}\X \Xsbar = \left(\Cbar \M{N} \Cbar^T\right) \X^T \M{L}\X \Xsbar
= \left(\Cbar \M{N} \Cbar^T \right) \M{L}_{\Xs} \Xsbar.
\end{align}
We note that these equations apply for continuous time. In discrete time, we will solve\ \eqref{eq:SP} for $\V{\alpha}=\left(\V{\Omega},\Ump\right)$, then update the tangent vectors on the unit sphere by $\V{\Omega} \D t$ using the nonlinear rotation\ \eqref{eq:RotUpdate}.

\subsubsection{Overdamped Langevin equation in $\Xsbar$ \label{sec:ODL}}
The deterministic equation\ \eqref{eq:dTau} takes the same form as that for rigid bodies, and in the blob-link picture it can be seen as describing the motion of a series of $N$ connected rigid rods. Following the same process used to derive the overdamped Ito Langevin equation \cite[Eq.~(12)]{BrownianMultiBlobs} in the rigid body case, we have the overdamped Ito Langevin equation 
\begin{gather}
\label{eq:ItoTau}
\dt{\Xsbar} = -\left(\Cbar\M{N}\Cbar^T\right)\M{L}_{\Xs} \Xs+ k_B T \div_{\Xsbar} \cdot \left( \Cbar\M{N}\Cbar^T\right)+ \sqrt{2 k_B T}\Cbar \M{N}^{1/2}\Wproc
\end{gather}
describing the evolution of $\Xsbar$. Here $\Wproc(t)$ is a collection of white noise processes (the formal derivative of Brownian motions), the divergence with respect to $\Xsbar$ is defined in\ \eqref{eq:DefDiv}, and $\M{N}^{1/2}$ satisfies the fluctuation-dissipation relation 
\begin{equation}
\label{eq:FDN}
\M{N}^{1/2}\left(\M{N}^{1/2}\right)^T=\M{N}.
\end{equation}
The arguments in \cite[Sec.~II.B.2]{BrownianMultiBlobs} can be used to show that\ \eqref{eq:ItoTau} is time-reversible with respect to the Gibbs-Bolztmann equilibrium distribution\ \eqref{eq:GBDist2}. As described there, because we formulate\ \eqref{eq:ItoTau} with respect to $\Xs$ and not $\V{X}$, there are no additional drift terms arising from metric/entropic forces.

An important point in practice is that $\M{N}^{1/2}$, which is not unique and only needs to satisfy\ \eqref{eq:FDN}, can be applied by adding Brownian noise with covariance $\Mfor$ to the right-hand side of saddle point system\ \eqref{eq:SP}. In an Euler-Maruyama discretization, this corresponds to solving the saddle point system
\begin{gather}
\label{eq:SPB1}
\M{K}\V{\alpha} = \Mfor \left(-\M{L}\V{X}+\M{\Lambda}\right)+\sqrt{\frac{2k_B T}{\D t}}\Mfor^{1/2}\gauss\\ \nonumber
\M{K}^T \M{\Lambda} = \V 0,
\end{gather}
where $\gauss$ is an i.i.d.\ vector of standard normal random variables. As discussed at length in \cite[Sec.~II(B)]{BrownianMultiblobSuspensions}, solving this saddle point system gives 
\begin{align}
\nonumber
\dt{\V{X}}= \M{K}\V{\alpha} &= -\left(\M{K}\M{N} \M{K}^T\right) \M{L}\V{X}+\sqrt{\frac{2k_B T}{\D t}}\left(\M{K}\M{N}\M{K}^T\right) \Mfor^{-1/2}\gauss\\
\label{eq:Nhalf}
& =  -\left(\M{K}\M{N} \M{K}^T\right) \M{L}\V{X}+\sqrt{\frac{2k_B T}{\D t}}\M{K}\M{N}^{1/2}\gauss,
\end{align}
giving $\M{N}^{1/2}=\M{N}\M{K}^T \Mfor^{-1/2}$. Thus, generating noise of the form $\M{N}^{1/2}\gauss$ reduces to the simpler process of solving a saddle point system with right hand side $\Mfor^{1/2}\gauss$. For hydrodynamics which is localized to each fiber, we do this using the eigenvalue decomposition of $\Mfor$, which already must be computed for the purposes of eigenvalue truncation (see Section\ \ref{sec:Mob}). Given that $\M{N}^{1/2}$ can also be computed via dense linear algebra, the real savings in the saddle point solve come when we need to generate $\M{N}^{1/2}$ with \emph{nonlocal} hydrodynamics (between the many fibers), where dense linear algebra is infeasible, but the action of $\Mfor^{1/2}$ can be computed via the Lanczos algorithm \cite{chow2014preconditioned} or the positively-split Ewald method \cite{PSRPY}. While the case of nonlocal hydrodynamics will not be treated in this paper, the saddle point method is a useful foundation for future work.

\subsubsection{Overdamped Langevin equation in $\V{X}$}
We now use\ \eqref{eq:ItoTau} to derive an overdamped Langevin equation in terms of $\V{X}$. If we multiply\ \eqref{eq:ItoTau} on the left by $\X$ and expand $\M{L}_{\Xs}=\X^T \M{L}\X $, we obtain
\begin{gather}
\label{eq:XEq1}
\dt{ \V X} = -\left(\M{K}\M{N}\M{K}^T\right)\M{L}\V{X}+ k_B T \div_{\Xsbar} \cdot \left( \M{K}\M{N}\Cbar^T\right)+ \sqrt{2 k_B T}\M{K}\M{N}^{1/2}\Wproc.
\end{gather}
Using the chain rule to write differentiation with respect to $\V{X}$ as
\begin{equation*}
\frac{\partial f}{\partial \V{X}_k} = \frac{\partial f}{\partial \Xsbar_p} \frac{\partial \Xsbar_p}{\partial  \V{X}_k}  = \frac{\partial f}{\partial \Xsbar_p} \X^{-1}_{pk}
\end{equation*}
We now rewrite the divergence in\ \eqref{eq:XEq1} as
\begin{equation}
\frac{\partial}{\partial \Xsbar_j}  \left( \M{K}\M{N}\Cbar^T\right)_{ij}=\frac{\partial}{\partial \V{X}_k}  \left( \M{K}\M{N}\Cbar^T\right)_{ij}\frac{\partial \V{X}_k}{\partial \Xsbar_j}=\frac{\partial}{\partial \V{X}_k}  \left( \M{K}\M{N}\Cbar^T\right)_{ij}\X_{kj}=\frac{\partial}{\partial \V{X}_k} \left(\M{K}\M{N}\M{K}^T\right)_{ik}
\end{equation}
so that the Ito equation\ \eqref{eq:ItoTau} could equivalently be formulated in terms of $\V{X}$ as might be expected from the deterministic equation\ \eqref{eq:dXDet},
\begin{align}
\label{eq:ItoX}
\dt{ \V{X}} &= -\left(\M{K}\M{N}\M{K}^T\right)\M{L}\V{X} + k_B T \div_{\V{X}} \cdot \left( \M{K}\M{N}\M{K}^T\right)+ \sqrt{2 k_B T}\M{K}\M{N}^{1/2}\Wproc \\[4 pt]
\label{eq:kinetic}
& \eqdi-\widehat{\M{N}}\M{L}\V{X} +  \sqrt{2 k_B T}\widehat{\M{N}} \circ \widehat{\M{N}}^{-1/2}\Wproc,
\end{align}
where $\widehat{\M{N}}=\M{K}\M{N}\M{K}^T$ and the second equality \new{denotes that paths of\ \eqref{eq:ItoX} and\ \eqref{eq:kinetic} have the same probability distribution.} Equation\ \eqref{eq:kinetic} is the Langevin equation written in a split Stratonovich-Ito \cite{BrownianMultiBlobs} or kinetic \cite{KineticStochasticIntegral_Ottinger} form, where the terms before the $\circ$ are evaluated at the midpoint of a given time step, while the terms after are evaluated at the beginning of the time step (c.f. \cite[Eq.~(26)]{BrownianMultiBlobs}). When we develop our numerical methods, we will do so with the equation for $\V{X}$ in mind, since ultimately we will evolve and track the fiber positions. That said, it will be simpler when analyzing the Langevin equation to work with the equation\ \eqref{eq:ItoTau} for $\Xsbar$, knowing that the one in $\V{X}$ can be obtained by this simple transformation. Note that\ \eqref{eq:ItoX} is much simpler than the overdamped equations derived previously for bead-link models \cite{morse2003theory, ottinger2012stochastic}.

\section{Temporal integration \label{sec:Tint}}
In this section, we discuss our temporal integrator for the overdamped Langevin equation\ \eqref{eq:ItoTau}. The scheme, which is in the spirit of the Fixman method \cite{fixman1978simulation} and similar to that of Westwood et.\ al for rigid bodies \cite{westwood2021generalised}, is able to integrate the overdamped Langevin equation using one saddle-point solve per time step. The key idea is to first move to the midpoint to compute the mobility, then solve a saddle point system using the midpoint values, which generates the required drift term in expectation.

Before introducing our numerical scheme, it is helpful to simplify the drift term in\ \eqref{eq:ItoTau}. We first separate it into three terms, 
\begin{equation}
\label{eq:DriftExpand}
\frac{\partial}{\partial \Xsbar_j} \left( \Cbar\M{N}\Cbar^T\right)_{ij} = \left(\partial_j \Cbar_{ik}\right) \M{N}_{kp} \Cbar^T_{pj}+\Cbar_{ik} \left(\partial_j \M{N}_{kp} \right) \Cbar^T_{pj} + \Cbar_{ik} \M{N}_{kp}\left( \partial_j \Cbar^T_{pj}\right),
\end{equation}
where $\partial_j$ is shorthand for $\partial/\partial \Xsbar_j$. As shown in \cite[Sec.~III(A)]{BrownianMultiBlobs}, rotating the tangent vectors at every time step $n$ using the Euler-Maruyama method 
\begin{equation}
\tdisc{\Xs}{n+1}=\text{rotate}\left(\tdisc{\Xs}{n}, \sqrt{2 k_B T \D t}\left(\tdisc{\M{N}}{n}\right)^{1/2} \tdisc{\gauss}{n}+\mathcal{O}(\D t)\right)
\end{equation}
where $\tdisc{\gauss}{n}$ is a vector of i.i.d.\ standard normal random variables, is sufficient to capture the first drift term. The third term in\ \eqref{eq:DriftExpand} is zero by\ \eqref{eq:DefDiv}. Therefore, our schemes simply need to generate the additional drift term
\begin{equation}
\label{eq:DriftGenC}
\left(k_B T\right) \Cbar_{ik} \left(\partial_j \M{N}_{kp} \right) \Cbar^T_{pj}\D t \quad \text{i.e.,} \quad \left(k_B T\right) \Cbar \left(\nabla_{\Xsbar} \M{N} \ddot \Cbar^T \right)\D t.
\end{equation}
\new{Because $\text{rotate}\left(\tdisc{\Xs}{n},\D t \V{\Omega}\right) = \tdisc{\Xs}{n}+\D t \M{C}\V{\Omega}-\left(\V{\Omega} \cdot \V{\Omega}\right)\D t^2\tdisc{\Xs}{n}/8+\mathcal{O}\left(\D t^3\right)$}, the term\ \eqref{eq:DriftGenC} corresponds to $\mathcal{O}\left(\D t^2\right)$ to a rotate procedure by an angle $\D t \left(k_B T\right)  \nabla_{\Xsbar} \M{N} \ddot \Cbar^T $. Thus our task will be to design numerical methods to produce the stochastic drift term in $\V{\Omega}$
\begin{equation}
\label{eq:DriftGen}
\text{Drift} = \left(k_B T\right) \nabla_{\Xsbar} \M{N} \ddot \Cbar^T  \quad \text{i.e.,} \quad \text{Drift}_k = \left(k_B T\right) \left(\partial_j \M{N}_{kp} \right) \Cbar_{jp}
\end{equation}
in expectation, which will give\ \eqref{eq:DriftGenC} after rotation over time step size $\D t$.

\subsection{Implicit methods \label{sec:ImpMeth}}
We first motivate our temporal discretization of\ \eqref{eq:ItoTau} by considering the discretization of the unconstrained \emph{linearized} SDE
\begin{equation}
\label{eq:UnconsSDT}
\dt{\V{X}} = -\Mfor\left[\V{X}_0\right] \M{L}\D \V{X}+ \V{U}_B,
\end{equation}
where $\V{X}_0$ is an equilibrium position, $\D \V{X}=\V{X}-\V{X}_0$, the matrix $\M{L}$ discretizes the energy $\mathcal{E}_\text{bend}=\D \V{X}^T \M{L}\D \V{X}$, and $\V{U}_B$ is the Brownian velocity given by fluctuation-dissipation balance as $\sqrt{2 k_B T}\left(\Mfor\left[\V{X}_0\right]\right)^{1/2}\Wproc(t)$. Because this SDE is unconstrained, the equilibrium covariance of $\D \V{X}$ is known from statistical mechanics, 
\begin{equation}
\label{eq:CovLin}
\Exp{\D \V{X} \D \V{X}^T}=k_B T \M{L}^{-1}
\end{equation}
In our case, the bending force resulting from the matrix $\M{L}$ is very stiff (fourth derivative), and so we need to discretize it implicitly. Our goal is to design numerical methods that preserve the covariance\ \eqref{eq:CovLin} for arbitrary $\D t$ when applied to\ \eqref{eq:UnconsSDT}. To do this, we follow the analysis of \cite[Sec.~III(B)]{DFDB} to derive the steady state covariance for a given temporal integrator.

We consider an implicit-explicit method for\ \eqref{eq:UnconsSDT} of the form
\begin{equation*}
\tdisc{\V{X}}{n+1} = \tdisc{\V{X}}{n}+\D t \left(-\Mfor \M{L}\left(c\tdisc{\V{X}}{n+1}+(1-c)\tdisc{\V{X}}{n}\right) + \tdisc{\V{U}_B }{n} \right),
\end{equation*}
which can be rearranged to yield
\begin{equation*}
 \left(\M{I}+c \D t \Mfor \M{L}\right)\tdisc{\V{X}}{n+1} =\left[ \left(\M{I}-(1-c) \D t  \Mfor \M{L}\right)\tdisc{\V{X}}{n} + \D t  \tdisc{\V{U}_B}{n}\right]. 
\end{equation*}
We now take an outer product of the two sides of the equation, then substitute the desired covariance from\ \eqref{eq:CovLin}, which at steady state is independent of the time step $n$, to obtain the matrix equation
\begin{gather}
\nonumber
2 k_B T \D t \Mfor+k_B T \left(2c-1\right) \D t^2\Mfor \M{L}\Mfor=\D t^2 \Exp{\V{U}_B\V{U}_B^T}.
\end{gather} 
To obtain the exact covariance for $c=1$ (backward Euler), we therefore set 
\begin{equation}
\label{eq:UB}
\V{U}_B=\sqrt{\frac{2 k_B T}{\D t}}\left(\Mfor^{1/2}\gauss+\sqrt{\frac{\D t}{2}} \Mfor \M{L}^{1/2} \widetilde{\gauss}\right),
\end{equation}
where $\widetilde{\gauss}$ is another standard normal random vector. Another option is to use Crank-Nicolson ($c=1/2$), which gives the exact covariance for arbitrary $\D t$ with the usual Brownian velocity $\sqrt{2k_B T/\D t} \Mfor^{1/2}\gauss$. While this choice has been preferred for other applications \cite{DFDB}, we find it to be less accurate than our ``modified'' backward Euler scheme for the higher order modes, which take too long to equilibrate using $c=1/2$. As such, we will use $c=1$ and the Brownian velocity\ \eqref{eq:UB} throughout this paper, where we can precompute $\M{L}^{1/2}$ via Cholesky or eigenvalue (dense matrix) decomposition (since $\M{L}$ is a block diagonal matrix for a suspension of fibers). When we have constraints and nonlinear updates, the velocity does not generate the exact covariance for arbitrary $\D t$, but it gives a covariance which converges more rapidly to the correct answer. 

\subsection{Midpoint scheme \label{sec:OneSolveMP}}
We can now present our ``midpoint'' method which can produce the drift term\ \eqref{eq:DriftGen} in expectation with only one saddle-point solve per time step. Given that an inextensible chain can be viewed as a collection of interacting rigid rods (tangent vectors), our method is similar in spirit to that of Westwood et al.\ \cite{westwood2021generalised} for rigid body suspensions, but differs in the ways we detail in Appendix\ \ref{sec:Brennan}. The method as presented here is optimized for dense linear algebra, in the sense that it requires only two mobility evaluations per time step if those mobilities can be stored as dense matrices. See Appendix\ \ref{sec:MobDrAlt} for modifications when the mobility cannot be stored as a dense matrix.

At each time step $n$, we perform the following steps
\begin{enumerate}
\item Compute a rotation rate for the tangent vectors based on the Brownian velocity $\Mfor^{1/2} \tdisc{\gauss}{n}$
\begin{equation}
\label{eq:OmegaTilde}
\tdisc{\V{\alpha}}{n,*} =\sqrt{ \frac{2k_B T}{\D t} }\left(\tdisc{\M{K}}{n}\right)^{-1} \left(\tdisc{\Msym}{n}\right)^{1/2} \tdisc{\gauss}{n}=
\begin{pmatrix} \tdisc{\V{\Omega}}{n,*}\\ \tdisc{\Ump}{n,*} \end{pmatrix}
\end{equation}
where $\M{K}^{-1}$ is defined in\ \eqref{eq:Kinv}.
\item Rotate the tangent vectors by $\tdisc{\V{\Omega}}{n,*} \D t/2$ to generate a new configuration
\begin{align}
\tdisc{\Xsbar}{n+1/2,*} &= \begin{pmatrix}\text{rotate}\left(\tdisc{\Xs}{n}, \left(\D t / 2\right) \tdisc{\V{\Omega}}{n,*}\right) \\ \tdisc{\Xmp}{n}+\left(\D t/2\right) \tdisc{\V{U}}{n,*}_\text{MP}\end{pmatrix}\\[4 pt]
\label{eq:TauTilde}
&= \sqrt{\frac{ k_B T \D t}{2}} \tdisc{\Cbar}{n}\left(\tdisc{\M{K}}{n}\right)^{-1}\left(\tdisc{\Msym}{n}\right)^{1/2}\tdisc{\gauss}{n}+\mathcal{O}(\D t^{3/2})
\end{align}
\item Evaluate the mobility $\tdisc{\Mfor}{n+1/2,*}$ and use it to compute the additional drift velocity using the random finite difference (RFD) \cite{BrownianBlobs} with $\delta \sim \sqrt{\D t}$, 
\begin{equation}
\label{eq:UMD}
\tdisc{\V{U}}{n}_\text{MD}=\sqrt{ \frac{2k_B T}{\D t} }\left(\tdisc{\Mfor}{n+1/2,*}-\tdisc{\Mfor}{n}\right) \left(\left(\tdisc{\Mfor}{n}\right)^{-1/2}\right)^T \tdisc{\gauss}{n}.
\end{equation}
This term might be impractical for large systems because it is based on solving a resistance problem to obtain $\left(\Mfor^{-1/2}\right)^T$.  Appendix\ \ref{sec:MobDrAlt} has an alternative approach which generates the same drift term via an RFD in which $\Mfor$ only has to be applied rather than inverted. 
\item To obtain the tangent vector rotation rates, solve the saddle point system 
\begin{gather}
\label{eq:SPmidp}
\tdisc{\begin{pmatrix} -\Msym & \left(\M{I}+c \D t \Msym \M{L}\right)\M{K}\\ \M{K}^T & \V 0 \end{pmatrix}}{n+1/2,*}
\tdisc{\begin{pmatrix} \V{\Lambda} \\ \V{\alpha} \end{pmatrix}}{n+1/2} \\ \nonumber =
\begin{pmatrix} -\tdisc{\Msym}{n+1/2,*} \M{L}\tdisc{\V{X}}{n}+\tdisc{\V{U}}{n}_B +\tdisc{\V{U}}{n}_\text{MD}\\  \V{0}\end{pmatrix}
\end{gather}
for $\tdisc{\V{\Lambda}}{n+1/2}$ and $\tdisc{\V{\alpha}}{n+1/2}=\left(\tdisc{\V{\Omega}}{n+1/2},\tdisc{\Ump}{n+1/2}\right)$. The Brownian velocity $\V{U}_B$ is defined in\ \eqref{eq:UB}, and the first part of it is used in\ \eqref{eq:OmegaTilde} to generate the midpoint configuration, i.e., the same $\gauss$ is used in steps 1, 3, and 4.
\item Update the fiber via\ \eqref{eq:RotUpdate},
\begin{gather}
\label{eq:OneSolveUpdate}
\tdisc{\Xsbar}{n+1}=\begin{pmatrix}
\text{rotate}\left(\tdisc{\Xs}{n}, \D t \tdisc{\V{\Omega}}{n+1/2}\right)\\ \tdisc{\Xmp}{n}+\D t \tdisc{\Ump}{n+1/2} \end{pmatrix}.
\end{gather}
\end{enumerate}
Solving\ \eqref{eq:SPmidp} yields (to leading order in $\D t$)
\begin{equation}
\label{eq:OmMidP}
\tdisc{\V{\alpha}}{n+1/2} = \tdisc{\left(\M{N}\M{K}^T \Msym^{-1}\right)}{n+1/2,*} \left(- \tdisc{\Msym}{n+1/2,*}\M{L}\V{X}+\tdisc{\V{U}}{n}_B +\tdisc{\V{U}}{n}_\text{MD}\right)
\end{equation}
In Appendix\ \ref{sec:DriftOS}, we show that using this value of $\V{\alpha}$ in the nonlinear update\ \eqref{eq:OmMidP} generates the drift term \eqref{eq:DriftGen} in expectation. Thus, after the rotation, we obtain dynamics consistent with\ \eqref{eq:ItoTau}.

In this paper, we discretize filaments with at most 30 Chebyshev nodes and consider hydrodynamics on one filament at a time. We therefore use direct solvers for the saddle point system\ \eqref{eq:SPmidp} in our implementation, which is available (along with python files for all dynamic examples in this paper) at \url{https://github.com/stochasticHydroTools/SlenderBody}. To implement the midpoint method efficiently for a blob-link discretization, where there are many blobs even on a single filament, we need to use an iterative solver for\ \eqref{eq:SPmidp}. In fact, if such solvers are accessible, then the midpoint temporal integrator can be applied as is using the positively split Ewald method \cite{PSRPY} to generate random displacements with covariance $\Mfor$, as has been done for suspensions of rigid bodies \cite{BrownianMultiblobSuspensions}. Further details will be provided elsewhere; here we will only use such an implementation of the blob-link discretization to compare to our spectral results.

\section{Equilibrium statistical mechanics \label{sec:EqStat}}
This section is devoted to the equilibrium statistical mechanics of fluctuating inextensible filaments. We focus first on comparing the Gibbs-Boltzmann distribution\ \eqref{eq:GBDist} for blob-link and spectral chains through \new{Markov Chain Monte Carlo (MCMC)} calculations, then transition to showing that our midpoint temporal integrator can also sample from the Gibbs-Boltzmann distribution if the time step size is sufficiently small. We focus on two examples: small fluctuations of a filament held near a curved base state, and free fibers. The former example is advantageous because it allows us to linearize both the inextensibility constraint (for static calculations) and the SDE\ \eqref{eq:ItoX} (for dynamics) and then break the dynamics into a set of modes. We can then compare the variance of each mode between the blob-link and spectral discretizations. Free fibers, however, are far more common and relevant in practice, and so they are our focus in this section. We relegate the details on the curved filament to Appendix\ \ref{sec:CurvX0}, and here give only a few summary statements to point to the parallels between the two examples. 

For free fibers, we will use end-to-end distance
\begin{equation}
r(t) = \frac{1}{L}\norm{\XPoly(s=0,t)-\XPoly(s=L,t)}.
\end{equation}
 as a metric to compare statistics across different discretizations. When configurations are sampled from the equilibrium distribution of inextensible filaments, the distribution of $r$ is approximately \cite{wilhelm1996radial}
\begin{equation}
\label{eq:EETh}
G(r) = \frac{1}{Z} \sum_{\ell = 1}^{\infty} \frac{1}{\left(\ell_p^* \left(1-r\right)\right)^{3/2}} \exp{\left(-\frac{\left(\ell-1/2\right)^2}{\ell_p^*(1-r)}\right)}H_2 \left(\frac{\ell-1/2}{\sqrt{\ell_p^* \left(1-r\right)}}\right)r^2,
\end{equation}
where $\ell_p^*=\ell_p/L$ is the dimensionless persistence length, $H_2(x)=4x^2-2$ is the second Hermite polynomial, and we have multiplied by the Jacobian factor $r^2$ to effectively compare\ \eqref{eq:EETh} to a one-dimensional histogram of distances that we will generate from our data.  We will also look at other distance metrics, including the distance from the fiber end to its middle, the end to the quarter point, and the distance between the two interior quarter points (the middle half of the fiber). 

\subsection{Quantifying the Gibbs-Boltzmann distribution with MCMC sampling \label{sec:EEMCMC}}

\begin{figure}
\centering
\includegraphics[width=\textwidth]{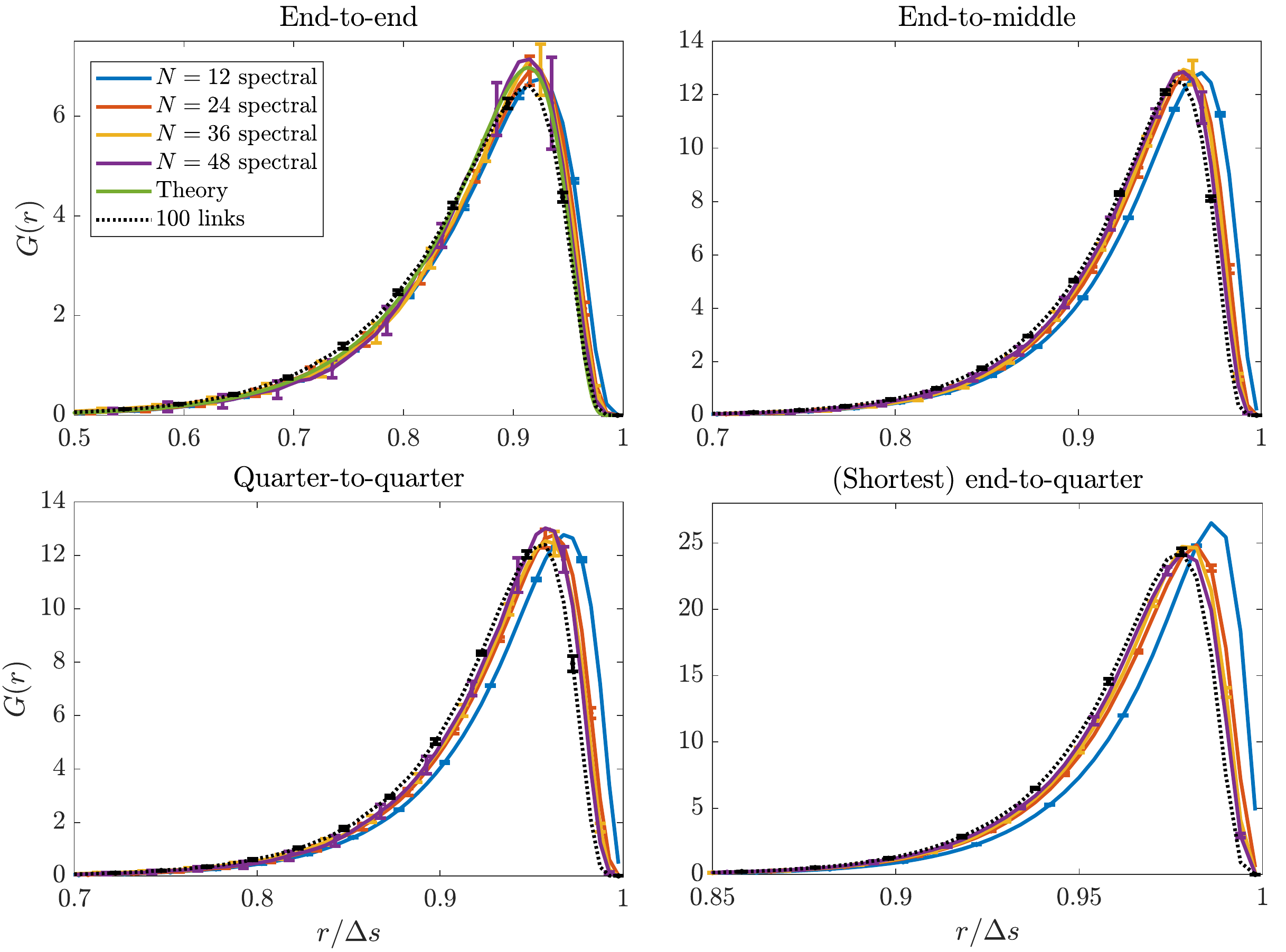}
\caption{\label{fig:MCMCFreeLp1} MCMC results for a freely fluctuating spectral fiber with $\ell_p/L=1$. We compute distances along the fiber including end-to-end ($\D s=L$, top left), end-to-middle ($\D s=L/2$, top right), quarter-to-quarter ($\D s=L/2$, bottom left) and the nearest end-to-quarter ($\D s=L/4$, bottom right). Note that there are twice as many observations in the right column as the left, since we measure from both the left and right endpoint. We compare the blob-link data for 100 links (which is the same with 200 links; the mean end-to-end distance is $\left(0.851 \pm 0.001\right)L$ in both cases) with several spectral discretizations and the theory\ \eqref{eq:EETh}. The spectral discretization results are approaching the 100 link results as $N$ increases (to within error bars), with a more rapid approach for larger $\D s$.}
\end{figure}

We first examine if the equilibrium distribution\ \eqref{eq:GBDist} for a freely fluctuating \emph{spectral} filament approximates well that of a blob-link chain. To do this, we use MCMC to sample the Gibbs-Boltzmann measure
\begin{gather}
\label{eq:EqdX}
d\pi\left( \V{X}\right) \propto \exp{\left(-\frac{\mathcal{E}\left(\V{X}\right)}{k_B T}\right)} d\mu_0 \left(\V{X}\right):=\widebar{\pi}\left(\V{X}\right)d\mu_0 \left(\V{X}\right),
\end{gather}
where $\mathcal{E}\left(\V{X}\right)= \frac{1}{2} \V{X}^T \M{L} \V{X}$ and the base measure $d\mu_0 \left( \V{X}\right)$ is defined in\ \eqref{eq:GBDist}. We generate a proposed configuration $\widetilde{\V{X}}$ by randomly rotating the tangent vectors (see Appendix\ \ref{sec:MCMCSpec} for details), which preserves the base measure. In this case, the probability of acceptance is simply the Metropolis factor
\begin{equation}
\label{eq:pacc}
p_\text{acc} =  \max{\left(\frac{\widebar{\pi}\left(\D \widetilde{\V{X}}\right)}{\widebar{\pi}\left(\D {\V{X}}\right)},1\right)}.
\end{equation}
We use this MCMC procedure to generate $10^6$--$10^7$ sample chains, removing the first 20\% of the chains as a burn-in period. Repeating this ten times to generate error bars, we report the distribution of end-to-end, end-to-middle, quarter-to-quarter, and end-to-quarter distances for both the spectral and blob-link chains in Fig.\ \ref{fig:MCMCFreeLp1}. We show only $\ell_p/L=1$, as the relative errors for $\ell_p/L=10$ are the same as $\ell_p/L=1$. The spectral discretization has a relatively small error even when $N=12$, and the distributions it generates move towards the blob-link ones as $N$ increases. This occurs at a faster rate for larger scales (end-to-end) than for smaller scales (end-to-quarter), as expected.

\begin{figure}
\centering
\includegraphics[width=\textwidth]{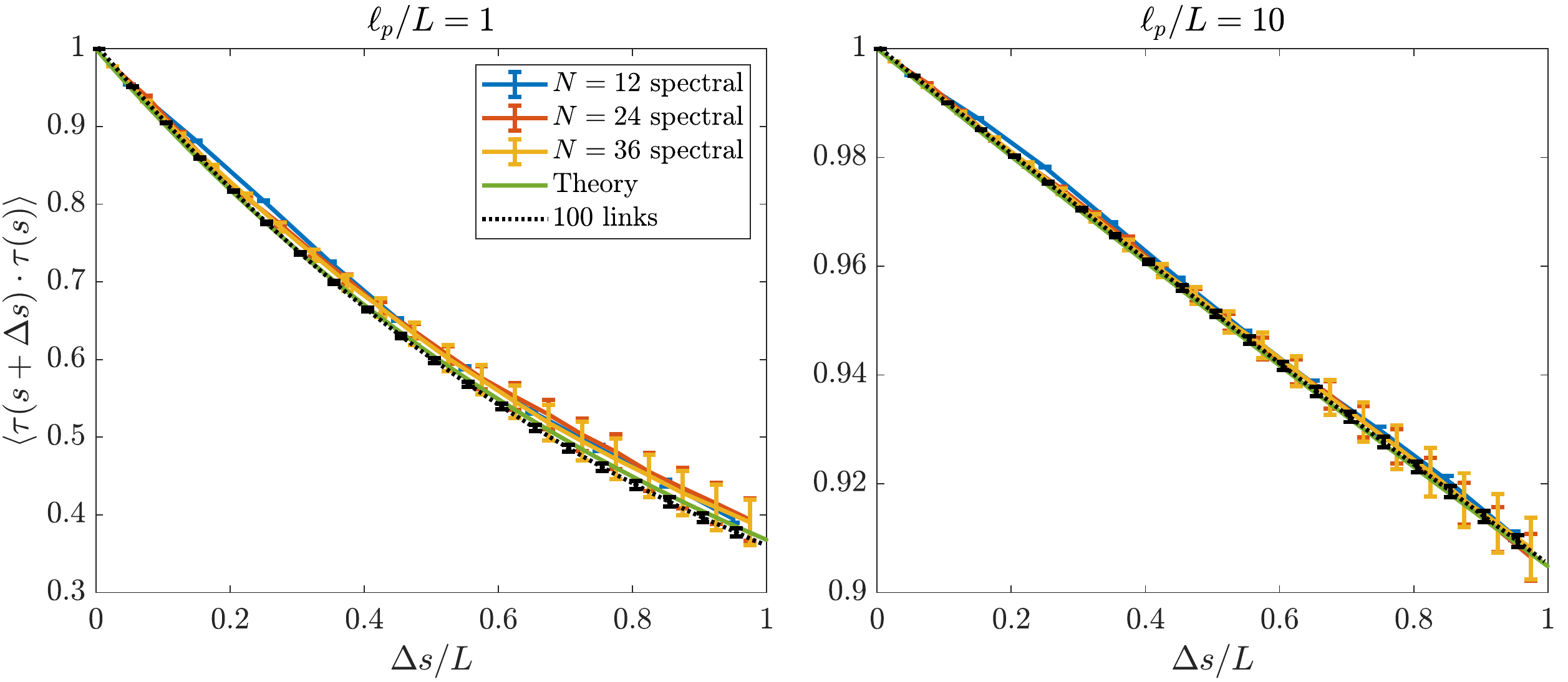}
\caption{\label{fig:TauCorrsSpec} Tangent vector correlation function $\langle \Xs\left(s + \D s\right) \cdot \Xs(s)\rangle$, computed from MCMC sampling. We show the data for a blob-link discretization with 100 links in black, and $N=12$, 24, and 36 spectral nodes in blue, red, and yellow. The theory $e^{-\Delta s/\ell_p}$ is shown in green. }
\end{figure}

An additional metric we can use to study the equilibrium distribution of a freely fluctuating chain is the correlation in the tangent vectors $\langle \Xs\left(s + \D s\right) \cdot \Xs(s)\rangle$ for $\D s \in [0,L]$. According to the definition of persistence length, this correlation should decay exponentially as $e^{-\D s/\ell_p}$. To measure the correlation function in the spectral discretization, we compute the correlation for all $\D s$ on the type 1 Chebyshev grid on which $\Xs$ is defined (see Fig.\ \ref{fig:SchDisc}). We then assign these measurements into bins corresponding to 10 (for $N=12$) or 20 (for $N\geq 24$) uniformly-spaced values of $\D s$ on $[0,L]$. Figure\ \ref{fig:TauCorrsSpec} shows how our spectral results compare to the theory (and 100 link discretization). While the blob-link chain has a correlation function which matches the theory exactly, there is an apparent small bias in our spectral chain at both small and large distances, with the correlation being larger than expected for small $N$. This bias is especially noticeable for $\ell_p/L=1$, but larger error bars for larger $N$ make it harder to make a definitive statement. 

In Appendix\ \ref{sec:MCMC2nd}, we perform a similar analysis for a filament undergoing small thermal fluctuations around a curved base state. By breaking the dynamics into a set of modes of the linearized covariance matrix, we show that the spectral method with $N=12$ nodes can successfully give the correct variance of the first ten modes (with a larger error for $\ell_p/L=1$ than for $\ell_p/L \gtrsim 10$), while $N=24$ and $N=36$ are sufficient to give the correct variance of the first 25 modes (see Fig.\ \ref{fig:MCMCPen}). Thus whether we consider small or large fluctuations, the Gibbs-Bolztmann distribution\ \eqref{eq:GBDist} for \emph{spectral} filaments is a good approximation of the more physical one for blob-link chains.

\subsection{Sampling with the midpoint temporal integrator}
We now discuss how the midpoint integrator of Section\ \ref{sec:OneSolveMP} can also give samples from the Gibbs-Bolztmann distribution\ \eqref{eq:GBDist}. To make our analysis in this regard universal, we need to understand how a certain time step size generalizes to a set of arbitrary parameters. We are once again aided by our example of a filament with small fluctuations, where we can linearize the SDE\ \eqref{eq:ItoX} around a certain state and compute a set of eigenmodes and associated timescales for the dynamics. This analysis, which we carry out in Appendix\ \ref{sec:LinTScales}, is a discrete version of that carried out by Kantsler and Goldstein \cite{kantsler2012fluctuations} in continuum, the difference being that the mobility in the latter case was approximated by local drag, so that the calculations could be done semi-analytically. For free filaments, the largest timescale in the problem is associated with the first ``fundamental'' bending mode \cite{kantsler2012fluctuations}, which we show in the inset of Fig.\ \ref{fig:Timescales}. The timescale associated with this mode is roughly
\begin{equation}
\label{eq:TauSm}
\tau_\text{fund} = 0.003 \frac{4 \pi \mu L^4}{\kappa \ln{\left(\epsRS^{-1}\right)}},
\end{equation}
and so we will report time in units of $\tau_\text{fund}$. There is a slight complication, however, as Fig.\ \ref{fig:Timescales} shows that the linearized timescales for two different $\epsRS$ do \emph{not} collapse onto the same curve when rescaled by the estimate of\ \eqref{eq:TauSm}. In fact, the expected log scaling, which comes from slender body theory \cite{johnson, krub}, approximately holds only for the smoothest modes ($k \lesssim 10$), with the timescales of the high-frequency modes scaling at a much sharper rate. Indeed, at the shortest scales we expect to see $\epsRS^{-1}$ scaling, corresponding to the timescales on which individual blobs relax (Stokes drag law).

\subsubsection{Required time step for midpoint integrator}
To examine the accuracy of the midpoint integrator relative to our MCMC calculations, we run Langevin dynamics from $t=0$ to $t=10\tau_\text{fund}$ on an initially straight filament using the RPY mobility\ \eqref{eq:McRPYdef} with slenderness $\epsRS=10^{-3}$. We record a histogram of the end-to-end distance after the first $\tau_\text{fund}$ (burn-in), ignoring the other distance metrics which we have already seen behave similarly.

\begin{figure}
\centering
\includegraphics[width=\textwidth]{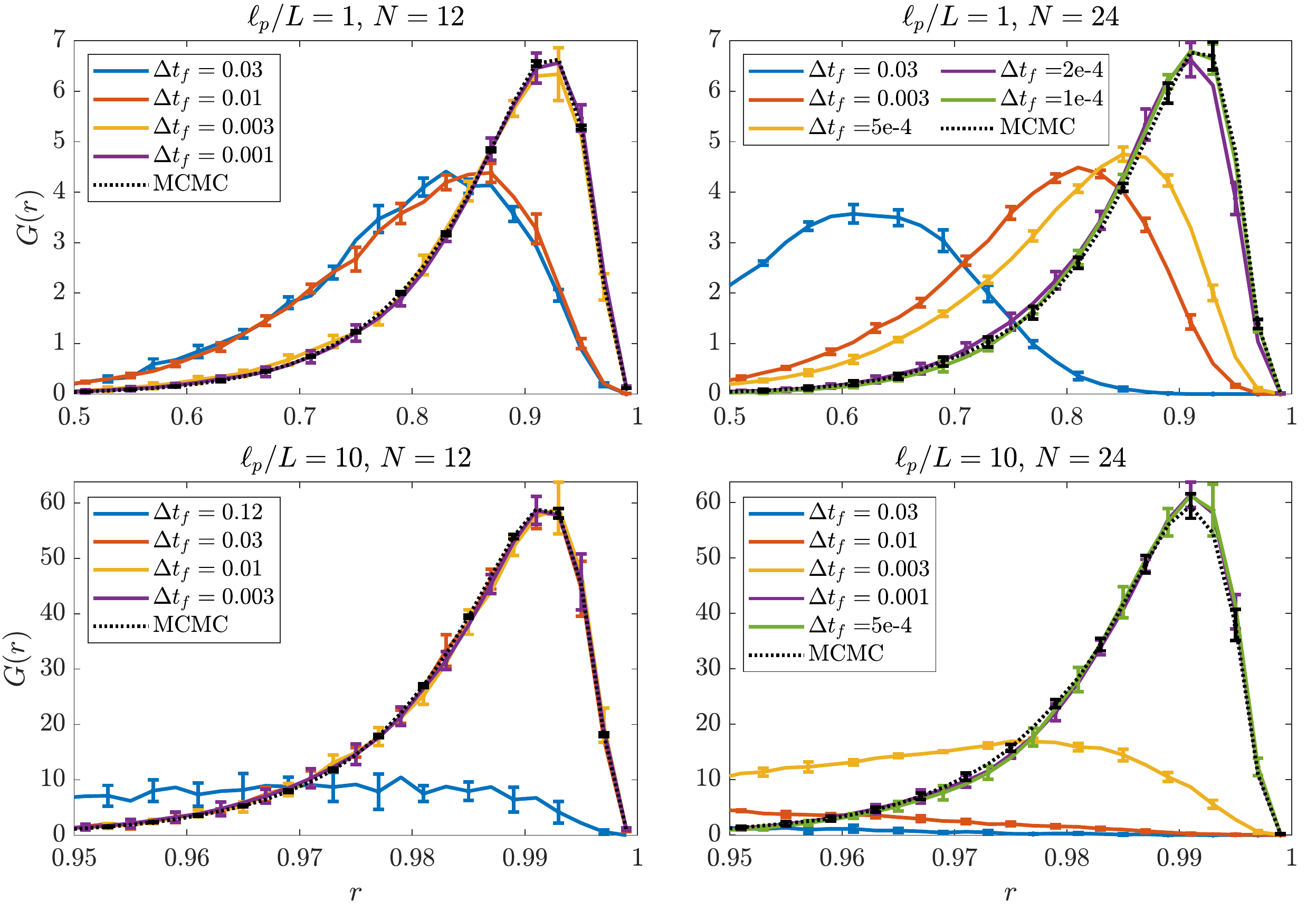}
\caption{\label{fig:OneSolveEERPY} Distribution of end-to-end distance for the midpoint temporal integrator discussed in Section\ \ref{sec:OneSolveMP} with the RPY mobility and $\epsRS=10^{-3}$. Time step sizes are reported as $\D t_f=\Delta t/\tau_\text{fund}$, where $\tau_\text{fund}$ is the slowest system relaxation time, defined in\ \eqref{eq:TauSm}. For each $\ell_p/L$ and $N$, we show the MCMC results in black, together with the distributions obtained from the midpoint integrator with various time step sizes. }
\end{figure}

Figure\ \ref{fig:OneSolveEERPY} shows that for sufficiently small $\D t$ the end-to-end distributions from Langevin dynamics converge to those of MCMC, validating our temporal integrator. It also gives an indication of how the required time step sizes change with $N$ and $\kappa$ (reported in terms of $\ell_p/L$). Focusing on $N$ first, we see that for a fixed $\ell_p/L$ the required time step size decreases by a factor of about 15 as $N$ doubles from 12 to 24. In Appendix\ \ref{sec:LinTInt}, we re-interpret this in terms of modes, finding that we need to resolve roughly twice as many modes when we double $N$. The scaling of the timescale of each mode $\tau_k \sim k^{-4}$ then implies a decrease in the time step size of 16, which points to the limitations of our temporal integrator for larger $N$. 

Switching our focus to $\kappa$, Fig.\ \ref{fig:OneSolveEERPY} shows that the \emph{relative} time step size required as we increase from $\ell_p/L=1$ to $\ell_p/L=10$ increases by a factor of roughly 10. In terms of modal analysis, the number of modes we need to resolve decreases by about 3 for every factor of 10 increase in $\kappa$ (see Appendix\ \ref{sec:LinTInt}), hence the increase in relative time step size. However, since the timescale $\tau_\text{fund}$ (and the timescale of each of the modes) scales like $1/\kappa$, the net effect of this behavior is no change in the \emph{absolute} time step size required for accurate equilibrium statistics. 

Appendix\ \ref{sec:LinTInt} also shows how the required time step size changes with $\epsRS$, although the analysis is less straightforward since there is no simple rescaling of time in this case. Our results show that the number of modes we need to resolve increases weakly as $\epsRS$ decreases, so that our time step size drops by a factor of roughly 5 when we drop from $\epsRS=10^{-2}$ to $\epsRS=10^{-3}$. 

\section{Dynamics of relaxation to equilibrium \label{sec:Relax}}
So far, we have only examined equilibrium statistical mechanics, finding that samples from the spectral and blob-link Gibbs-Boltzmann distributions generate similar statistics for a given set of parameters. But what about dynamics, and in particular, resolving the hydrodynamic interactions in slender filaments? 

The temporal integrator we developed here performs similarly regardless of the spatial discretization, in the sense that the number of modes we need to resolve scales with $N$. Thus, if we want to simulate slender filaments without having to take unreasonably small time steps, our only hope is to resolve the hydrodynamic interactions with a small number of collocation points, or a number of collocation points that is independent of $\epsRS$. For a direct blob-link discretization, it has already been established that $1/\epsRS$ beads are required to resolve hydrodynamics \cite{ttbring08, kallemov2016immersed}, although one could use an asymptotic theory like slender body theory (SBT) to model the hydrodynamics approximately \cite{saint1, saint2}. But in our spectral discretization, we can resolve \emph{deterministic} hydrodynamics with $\mathcal{O}(1)$ points (i.e., $N$ independent of $\epsRS$) \cite[Sec.~4.4]{TwistBend}; in this section (and Appendix\ \ref{sec:AdvSpec}), we verify that this is also the case for \emph{Brownian} hydrodynamics. 

To do this, we consider the dynamic problem of an initially straight semiflexible chain relaxing to its equilibrium fluctuations \cite{poelert2012analytical,obermayer2009freely, kantsler2012fluctuations,dimitrakopoulos2004longitudinal}. Our focus here is on the relaxation of the mean end-to-end distance to its mean value, i.e., to the mean of the distributions shown in Fig.\ \ref{fig:OneSolveEERPY}. As discussed in \cite{obermayer2009freely}, the scenario that we simulate is not really physical, since it is not possible for a fluctuating chain to ever reach an exactly straight configuration. As such, the more physically-relevant timescales are those that correspond to long-wavelength modes, where the shorter wavelength modes (which affect the end-to-end distance relatively little) have already reached their equilibrium state. Thus we will accept errors in the end-to-end distance on short timescales and concentrate on long-time behavior \cite{dimitrakopoulos2004longitudinal}. 

\begin{figure}
\centering
\includegraphics[width=\textwidth]{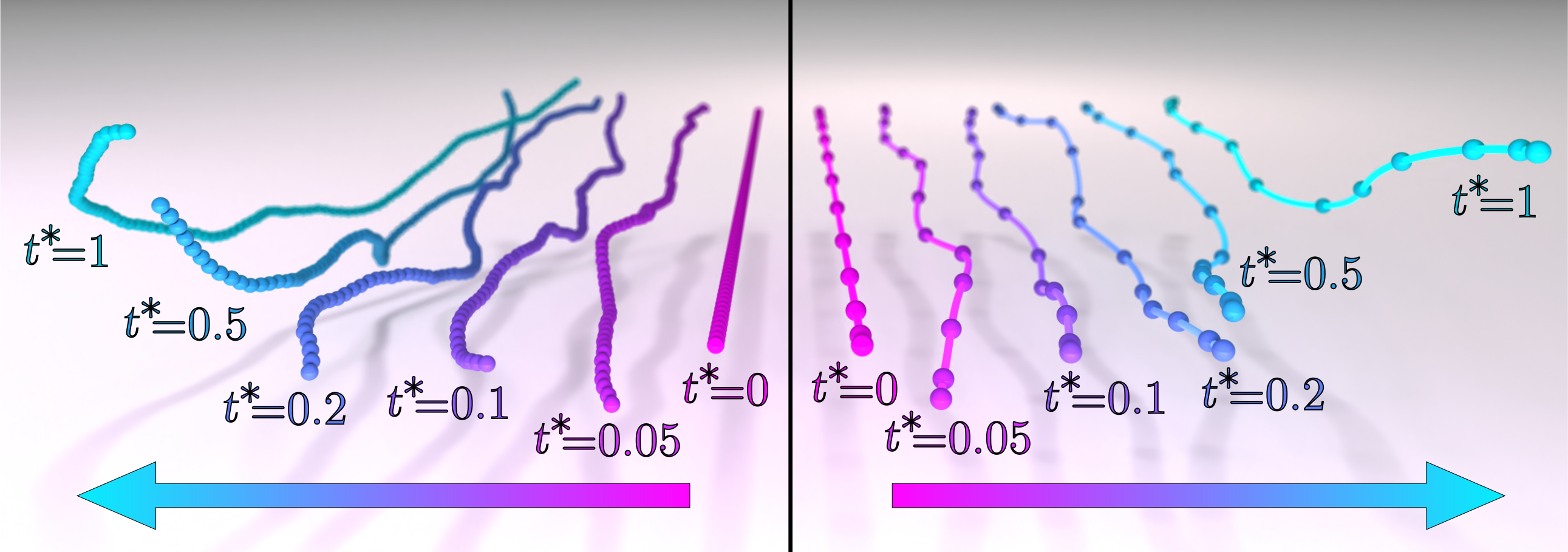}
\caption{\label{fig:RelFibPics} Two samples of free fibers with $\ell_p/L=1$ relaxing to their equilibrium fluctuations. In both cases $\epsRS=10^{-2}$, but at left we use a blob-link discretization, while at right we use a spectral discretization with $N=12$ nodes (we do not try to match the random forcing between the two discretizations, so the comparison is only in the qualitative look of the chain). The fibers are initialized straight (magenta $t^*=0$ lines), and then assume their equilibrium end-to-end length on the timescale $\bar{t}$ (defined in\ \eqref{eq:tbar}). In both cases, realizations of the fibers at different times are artificially staggered along the direction of the arrows for visual clarity.}
\end{figure}

On long timescales, numerical results verify that the data for various $\mu$, $L$, $\kappa$, $\epsRS$, and $k_B T$ can be (roughly) collapsed onto a single master curve with rescaled time and end-to-end variables
\begin{gather}
\label{eq:tbar}
t^* = \frac{t}{\bar{t}}, \qquad \bar{t} = 0.0008 \frac{4 \pi \mu L^4}{\kappa \ln{\left(\epsRS^{-1}\right)}}\approx 0.27 \tau_\text{fund},\qquad
r^* = \frac{r-\bar{r}}{1-\bar{r}},
\end{gather}
where $\bar{r}$ is the mean end-to-end distance computed in Section\ \ref{sec:EEMCMC} and $\bar{t}$ is the long-time decay rate, i.e., $r^*(t^*=1) \approx e^{-1} \approx 0.3$. In terms of our modal analysis, the timescale $\bar{t}$ is between the longest and second-longest timescales in the system (see Fig.\ \ref{fig:Timescales}), meaning that all of the modes except the first should be relaxed by $\bar{t}$. Simulating until $\bar{t}$ thus provides a set of intermediate times at which we can measure non-equilibrium statistics (see Fig.\ \ref{fig:RelFibPics} for pictures of the relaxation process).

To compare the spectral method to a blob-link method, we fix $\ell_p/L=1$ and compare four different discretizations of the chain: the spectral discretization with $N=12$ (which requires a time step size $\D t_f \approx 0.003$ for accurate dynamics; this time step is the same as that needed for equilibrium statistical mechanics), $N=24$ (time step size $\D t_f \approx 2 \times 10^{-4}$), and $N=36$ ($\D t_f \approx 7 \times 10^{-5}$) and the blob-link discretization with 100 blobs (required time step size $\D t^*\approx 10^{-4}$). The blob-link discretization is considerably more expensive to simulate for small $\epsRS$ (even with a GPU-accelerated implementation), which requires us to limit our comparison to $\epsRS=10^{-2}$.

\begin{figure}
\centering
\includegraphics[width=\textwidth]{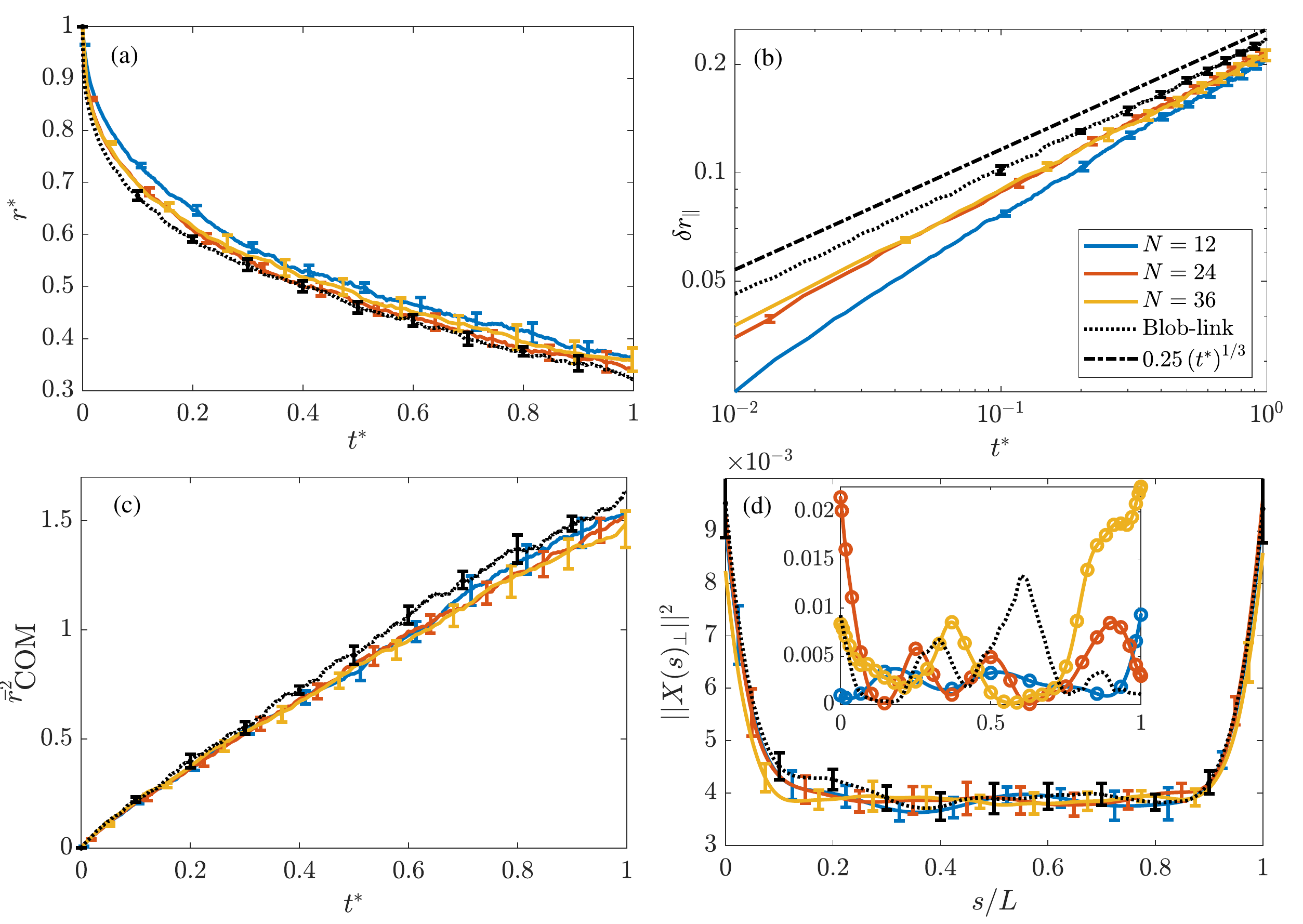}
\caption{\label{fig:SpecSecCompare} Comparing the spectral and blob-link trajectories for relaxation of a straight filament to its equilibrium fluctuations. We set $\ell_p/L=1$ and $\epsRS=10^{-2}$ and simulate until $t^*=1$, as defined in\ \eqref{eq:tbar}, comparing spectral discretizations with $N=12$ (blue), $N=24$ (red), and $N=36$ (yellow) to the blob-link discretization (black). (a and b) The end-to-end distance over time. In (a), we isolate the \emph{dynamic} error by normalizing each curve by its average $\bar{r}$ obtained from MCMC in Section\ \ref{sec:EEMCMC} (for $N=12$ this is $\bar{r}=0.865 \pm 0.001$, for $N=24$ it is $\bar{r}=0.862 \pm 0.005$, for $N=36$ it is $\bar{r}=0.859 \pm 0.006$, and for the blob-link it is $\bar{r}=0.851 \pm 0.001$). In (b), we plot $\delta r_\parallel$, which is the change in end-to-end distance in the parallel direction, for the four different discretizations, observing a $1/3$ power law scaling at long times. (c) The squared-norm of the displacement of the center of mass, normalized by the theoretical value for a rigid fiber (see text for details). (d) The squared perpendicular displacement at $t^*=0.1$. The inset shows a single sample, where we observe some high-frequency behavior in the blob-link results that is smoothed by the spectral method.}
\end{figure}

Figure\ \ref{fig:SpecSecCompare} shows how the two discretizations compare with each other for three different statistics: the average end-to-end distance $r^*=\langle \norm{\XPoly(0,t)-\XPoly(L,t)}\rangle$, the mean-square displacement of the center-of-mass $r^2_\text{COM}=\langle \norm{\XPoly(L/2,t)-\XPoly(L/2,0)}^2 \rangle$, and the average square perpendicular displacement $\langle \norm{\XPoly_\perp(s,t)}^2\rangle$ at $t^*=0.1$. We normalize the center-of-mass MSD by the value at $t=\bar{t}$ for a rigid fiber, which is $\left(2k_B T \bar{t}\right)\, \text{trace}\left({\M{N}_\text{tt}}\right)$, where $\M{N}_\text{tt}$ is the $ 3 \times 3$ matrix relating forces on a rigid fiber to its translational velocity. The normalized displacement is denoted $\bar{r}^2_\text{COM}$.

To separate the error in the \emph{dynamics} from that of equilibrium statistical mechanics, in Fig.\ \ref{fig:SpecSecCompare}(a), we normalize $r^*$ by the average $\bar{r}$ obtained from MCMC in Section\ \ref{sec:EEMCMC}. With this normalization, there is little difference between the spectral method with $N \geq 24$ and the blob-link method at later times, and the difference between the spectral method with $N=12$ and the blob-link method is small. Combining this with Fig.\ \ref{fig:SpecSecCompare}(c), which shows that the diffusion of the center of mass is the same (within error bars) across the different discretizations, we can conclude that a small number of spectral nodes can indeed resolve the dynamics at later times, as desired. This is an important statement because the spectral method with $N=12$ (resp.\ $N=24$) uses a time step that is two (resp.\ one) orders of magnitude larger than that of the blob-link chains, as well as a number collocation points that is an order of magnitude fewer. 

To compare our results to theory, in Fig.\ \ref{fig:SpecSecCompare}(b), we plot the shortening of the end-to-end distance projected onto the initial tangential direction, $\delta r_\parallel(t) = \langle L- \left(\XPoly(L,t)-\XPoly(0,t)\right) \cdot \Xs(t=0) \rangle$. Because we no longer normalize by the equilibrium end-to-end distance, we see a larger difference between the spectral and blob-link codes. Focusing on long times, we see that the data approach a $1/3$ power law for $t^* \gtrsim 0.1$, with a faster growth for short times, matching what is observed in \cite[Fig.~5(d)]{obermayer2009freely}. This 1/3 exponent is predicted to be universal independent of the way the initial state is prepared; see the second column in \cite[Table~I]{obermayer2009freely}. 

At early times, we see a more significant difference between the blob-link and spectral code, which can be explained by the fast relaxation of high-frequency modes. In Fig.\ \ref{fig:SpecSecCompare}(d) we examine the perpendicular displacement along the curve at an early time of $t^*=0.1$. In the inset, we show a single sample of the (squared) perpendicular displacement, and observe the small-length fluctuations in the blob link code which appear at short times. These fluctuations, which control the early-time relaxation, are smoothed out by the spectral code, and therefore treated incorrectly. At late times, their contribution is sufficiently small for the spectral and blob-link codes to match. 

\begin{figure}
\centering
\includegraphics[width=\textwidth]{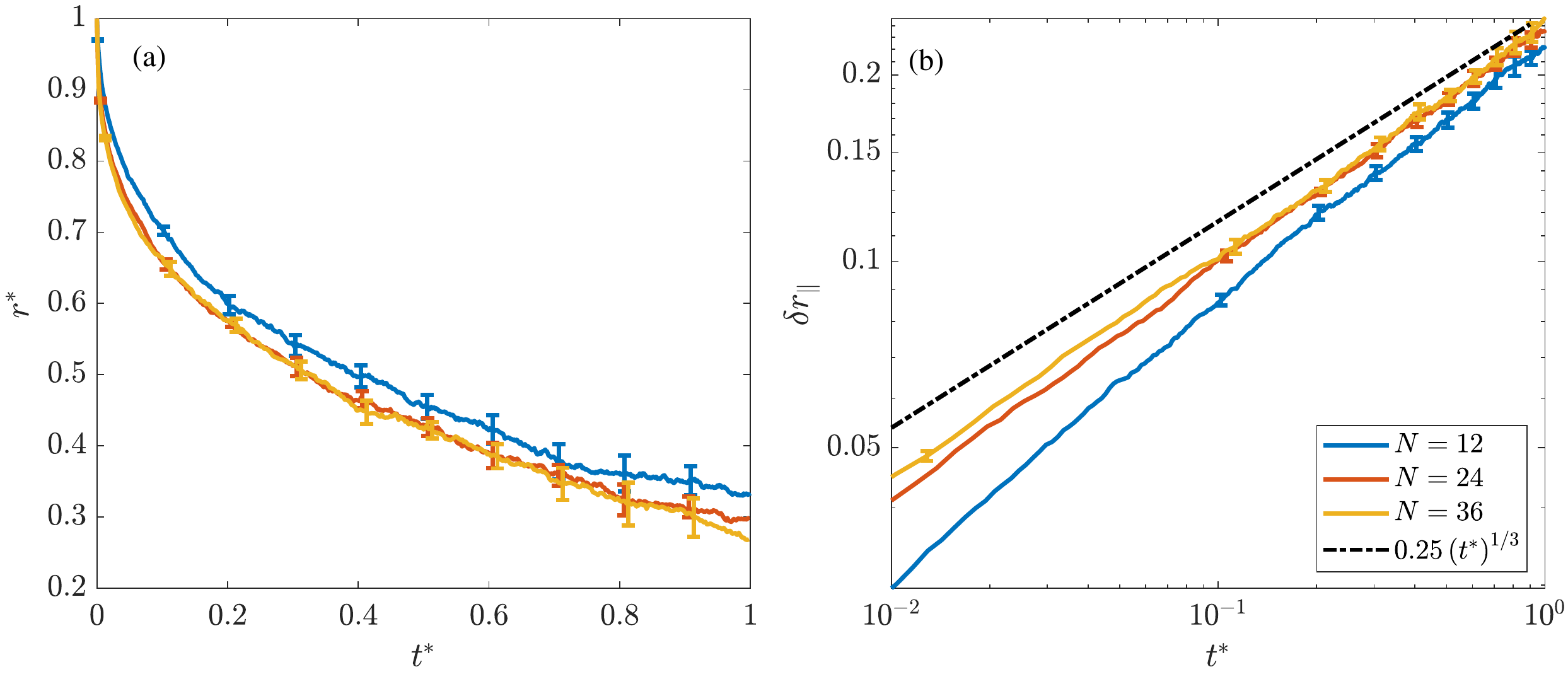}\caption{\label{fig:Eps3Conv} Fiber relaxation results when $\epsRS=10^{-3}$. In this case the blob-link method is too expensive to simulate, so we show the results in the end-to-end distance for the spectral method only with $N=12$ (blue), $N=24$ (red), and $N=36$ (yellow). These plots are analogous to Fig.\ \ref{fig:SpecSecCompare}(a) and (b), and the difference between the two discretizations is similar to the $\epsRS=10^{-2}$ case, as desired.}
\end{figure}

While we cannot obtain statistics for the blob-link code when $\epsRS=10^{-3}$ because of the expense in resolving hydrodynamics, we can still repeat our fiber relaxation test using $\epsRS=10^{-3}$ in the spectral method. Figure\ \ref{fig:Eps3Conv} shows the results of this test in end-to-end distance (compare to Fig.\ \ref{fig:SpecSecCompare}(a) and (b)). We see that the error between the three different values of $N$ is roughly the same across the two aspect ratios, which indicates that the number of points required for a given accuracy is \emph{not} sensitive to the fiber aspect ratio. In addition, we continue to observe the $1/3$ universal power law scaling for $\delta r_\parallel$ at long times \cite{obermayer2009freely}. Without blob-link data, it is difficult to say this for certain how the error in the spectral code scales with $\epsRS$. Still, the data strongly suggest that we can effectively resolve hydrodynamics of slender fibers with the same number of points, independent of the fiber aspect ratio.

\section{Bundling of transiently cross-linked semiflexible filaments \label{sec:Bundle}}
In previous work \cite{ActinBundlingDynamics}, some of us examined the dynamics of bundling in transiently-cross-linked actin networks. These networks form an important part of the cytoskeleton, and their arrangement into bundled structures is a topic of interest for cell division and motility \cite{alberts, chaubet2020dynamic}. As such, they have also been reconstituted in vitro \cite{gardel2004elastic, falzone2012assembly, schmoller2009structural, lieleg2009structural}, and the focus of our previous study \cite{ActinBundlingDynamics} was to compare our computational results to those observed experimentally in \cite{falzone2012assembly}. In \cite{ActinBundlingDynamics}, we included fluctuations by modeling the filaments as rigid, and made statements about the role of thermal fluctuations in bundling dynamics for networks of rigid filaments, leaving semiflexible filaments for future work. Having developed a temporal integrator for semiflexible filaments, we are now ready to complete our investigation by studying the role of semiflexible fluctuations in bundling of cross-linked fiber networks. We emphasize that there are \emph{no} hydrodynamic interactions between different fibers in these simulations.

To simulate cross linking, we couple the filament model developed here with a Markov chain describing the transiently-bound cross linkers (CLs). Our model of cross linking is laid out in detail in \cite{ActinCLsRheology, ActinBundlingDynamics}, so here we present only a summary. Each filament is divided into $N_u$ uniformly-spaced binding sites separated by distance $\Delta s_u = L/(N_u-1)$. Assuming that the CLs diffuse rapidly relative to the filaments, the binding of one end of a CL to one of these sites can be approximated by a single rate $\kon$ with units 1/(length$\times$time). Once the first end is bound, the second end can bind to a nearby filament with rate
\begin{equation}
\label{eq:koneqn}
\konb\left(\ell_k'\right)= \konb^0 \exp{\left(-\frac{K_c}{2}\frac{\left(\ell_k'-\ell_c\right)^2}{k_B T}\right)},
\end{equation}
where $\ell_k'$ is the deformed length of the CL (distance between the pair of binding sites), $\ell_c$ is the rest length of the CL, and $K_c$ is its stiffness. The relationship\ \eqref{eq:koneqn} ensures that the CL dynamics are in detailed balance; that is, the links are passive and do not consume energy. To efficiently search for nearby pairs of filaments, we limit $\ell_k'$ to two standard deviations of the Gaussian\ \eqref{eq:koneqn}; that is, we only search for pairs of binding sites $2\sqrt{k_B T/K_c}$ apart. Each of the binding reactions has an associated unbinding (reverse) reaction with a rate on the order 1/s \cite{kuhlman1994kinetics}, so that there are a total of four possible reactions which are simulated using a version of the standard Stochastic simulation / Gillespie algorithm \cite{gillespie2007stochastic,  anderson2007modified,ActinCLsRheology}.

We use a time splitting algorithm to update the filaments and cross linkers in sequence. At each time step, we take a step $\D t$ of the stochastic simulation algorithm with the filament positions fixed. This gives pairs of binding sites that are bound together, and consequently a force $\FCL$ exerted on the corresponding filament pairs. In previous work \cite[Sec.~6.1]{FibersWeakInextensibility}, we spread this force as a smoothed delta function around the cross linker binding location, ensuring smoothness of the cross linking force density and the subsequent fiber shapes. Since the smoothness assumption doesn't apply to fluctuating filaments, it is more physical to use instead the spring cross linking energy
\begin{equation}
\label{eq:ECL}
\mathcal{E}_\text{CL} = \frac{K_c}{2}\left(\norm{\ind{\XPoly}{i}\left(s_i^*\right)-\ind{\XPoly}{j}\left(s_j^*\right)}-\ell_c \right)^2
\end{equation}
between points $\ind{\XPoly}{i}\left(s_i^*\right)$ (on fiber $i$) and $\ind{\XPoly}{j}\left(s_j^*\right)$ (on fiber $j$). If we introduce the matrix $\M{R}_u$ which resamples the Chebyshev interpolant $\XPoly$ at uniformly-spaced binding sites, this energy can be rewritten in terms of the Chebyshev collocation points $\ind{\V{X}}{i}$ and $\ind{\V{X}}{j}$, and the resulting force computed by differentiating the energy with respect to $\V{X}$. The final expression for the force $\FCL$ at point $p$ for a CL attached to binding site $k$ then becomes the standard force for a spring (equal and opposite at the two fibers $i$ and $j$) multiplied by the $(k,p)$ entry of $\M{R}_u$. These forces at each time step become additional forces in the Langevin equation\ \eqref{eq:ItoX}, so that at each time step we solve (for each fiber independently)
\begin{gather}
\label{eq:ItoCL}
\dt{ \V{X}} = \M{K}\M{N}\M{K}^T\left(-\M{L}\V{X} +\FCL\right)+ k_B T \div_{\V{X}} \cdot \left( \M{K}\M{N}\M{K}^T\right)+ \sqrt{2 k_B T}\M{K}\M{N}^{1/2}\Wproc.
\end{gather}
by replacing $-\M{L}\tdisc{\V{X}}{n}$ with $-\M{L}\tdisc{\V{X}}{n}+\FCL\left(\tdisc{\V{X}}{n}\right)$ on the right hand side of the saddle point solve\ \eqref{eq:SPmidp} (that is, we treat the spring forces explicitly in time). Quantitative comparison of simulations with the smoothed forcing from \cite[Sec.~6]{FibersWeakInextensibility} and the new energy-based forcing from\ \eqref{eq:ECL} shows little difference between the two models.

Throughout this section, we will use the parameters given in \cite[Table~1]{ActinBundlingDynamics}. Just as in that study, we consider filaments with initial mesh size $0.2$ $\mu$m, which corresponds to 200 filaments of length $L=1$ $\mu$m in a periodic domain of edge length $L_d=2$ $\mu$m and 675 filaments of length $L=1$ $\mu$m in a periodic domain of edge length $L_d=3$ $\mu$m (as discussed in \cite{ActinBundlingDynamics}, the results are repeatable as we increase the domain size until the structure begins to collapse into 1 or 2 bundles). The question we will examine is how the behavior changes as we increase $\ell_p/L$, so we will leave all parameters constant except the bending stiffness $\kappa$. This includes $k_B T = 4.1 \times 10^{-3}$ pN$\cdot \mu$m, which corresponds to the thermal energy at room temperature. For our spatial and temporal discretization, we use $N=12$ and $\D t = 10^{-4}$ s over all simulations, having verified that doubling the number of points and halving the time step size does not change the results within statistical error. Our explicit treatment of the forces from the CLs, whose base stiffness of $K_c=10$ pN/$\mu$m is an order of magnitude estimate for the effective stiffness of $\alpha$-actinin \cite{le2017mechanotransmission, falzone2012assembly}, limits the time step size. In particular, resolving the spring dynamics automatically resolves the equilibrium statistical mechanics.\footnote{If we convert the time step sizes from Fig.\ \ref{fig:TintLinConvLp} to the units we use here, we need a time step size of $3 \times 10^{-3}$ s, $1.4 \times 10^{-3}$ s, and $1.2 \times 10^{-3}$ s to resolve the equilibrium fluctuations for $\ell_p/L=1$, 10, and 100, respectively. All of these are at least an order of magnitude larger than that required to resolve the explicit treatment of CLs.}

\subsection{Visualizing the bundling process}

In previous work \cite{ActinBundlingDynamics}, we showed that bundling of filaments occurs via a thermal zippering mechanism, where cross linkers stretch to bind nearby pairs of filaments, then contract to their rest length. This contraction pulls the filaments closer together, which allows for binding of additional cross linkers. The resulting equilibrium configuration contains filaments which are aligned in parallel and spaced roughly a distance equal to the cross linker rest length \cite[Fig.~1]{ActinBundlingDynamics}. In \cite{ActinBundlingDynamics}, we showed that this small-scale ratcheting mechanism also leads to large scale bundling for networks of rigid diffusing filaments, and that the bundling process roughly occurs in two stages: first, individual filaments come together into small bundles of a few filaments. Then, the bundles begin to coalesce, forming bundles of bundles and eventually one very large bundle. To quantify this process, we map the filament connections to a connected graph, where a connection in the graph exists when at least two CLs spaced $L/4$ apart connect the two filaments. The ``bundles'' are then connected regions in this graph. We track over time the bundle density, defined as the number of bundles per unit volume, and see a peak before the bundles start to coalesce. As discussed in \cite{ActinBundlingDynamics}, our definition of bundle density, while arbitrary, is a good way to compare the dynamics across multiple systems.

\begin{figure}
\includegraphics[width=\textwidth]{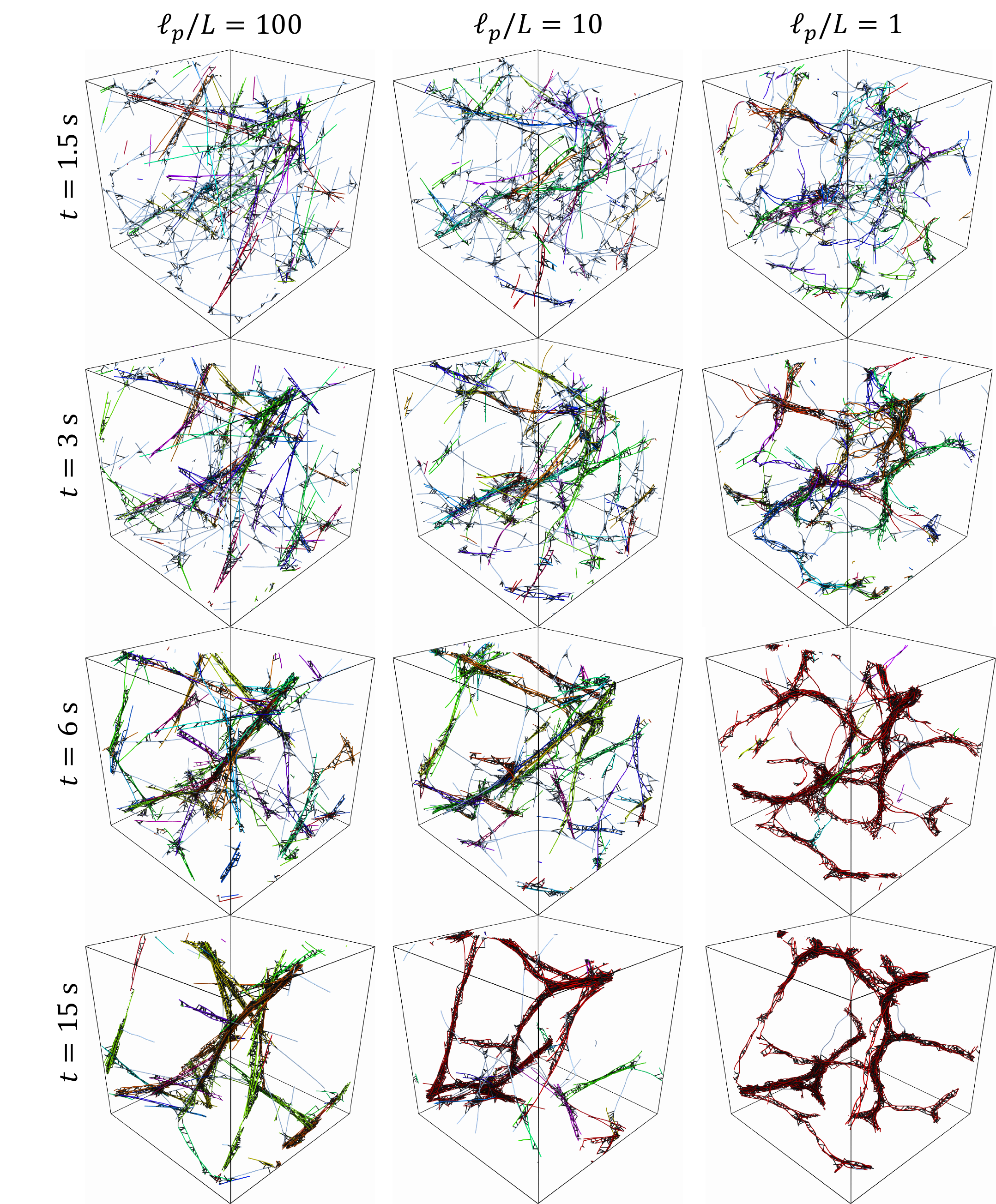}
\caption{\label{fig:BundlePic} Snapshots of the bundling process in networks of varying fiber stiffness. Colored fibers are actin filaments (colored by bundle), while black fibers are the cross linkers. The bundling process is fastest for the most flexible fibers, where the final morphology shows curved bundles.}
\end{figure}

The same fundamental process plays out in networks of semiflexible filaments, as shown in Fig.~\ref{fig:BundlePic}, where we show snapshots of the bundling process at different time points for three different orders of magnitude of $\ell_p/L$ (these plots show 200 filaments with $L_d=2$). The top plots ($t=1.5$~s) correspond to initial stage of bundling, where there are many bundles of a few filaments, while subsequent plots begin to show coalescence of the bundles. There are a few takeaways here: first, we see that the bundle morphology looks qualitatively different as we decrease $\ell_p/L$, with smaller persistence length having more curved fibers and therefore more curved bundles. Furthermore, the smaller persistence length bundles agglomerate faster, and at a given time they appear more clumped (especially $t=6$ s). It is not clear, however, to what extent these differences in bending deformations are driven by CL forces vs.\ the thermal bending fluctuations.\footnote{A similar question was studied in \cite{brangwynne2007force} for microtubule networks, where the authors showed that large \emph{nonthermal} forces combine with polymerization dynamics to generate bent microtubule shapes in cells.} Indeed, we did show in previous work \cite{ActinBundlingDynamics} that agglomeration (the second stage of the bundling process) happens faster for non-fluctuating filaments that are less stiff, since they are able to be bent easier by the cross linkers. Thus the main question here is whether semiflexible bending fluctuations themselves speed up (or slow down) the bundling process, or whether CL forces dominate.

\subsection{Quantifying the role of semiflexible bending fluctuations}
To get at this question, we need to dissect the evolution of the actin filaments in a cross-linked network into the three possible ways they can move: action by CL forces, thermal rotation and translation (keeping the fiber shape fixed), and bending fluctuations. In previous work \cite{ActinBundlingDynamics}, we studied the first two of these, showing in particular that deterministic filaments (those that can only move by CL forces) behave the same as rigid ones when $\ell_p/L \sim 10$. We used this assumption to justify neglect of thermal bending fluctuations for \emph{actin} filaments of length 1 $\mu$m, for which $\ell_p/L \approx 17$ \cite{gittes1993flexural}. In fact, we neglected any bending and treated the filaments as rigid, so that they evolved under both CL forces and thermal rotation and translation. We showed that translational and rotational diffusion accelerates the bundling process significantly, since there is more mixing and filaments are able to find each other faster (assuming there are sufficiently many CLs to link two filaments that are close together). 

\begin{figure}
\centering
\includegraphics[width=\textwidth]{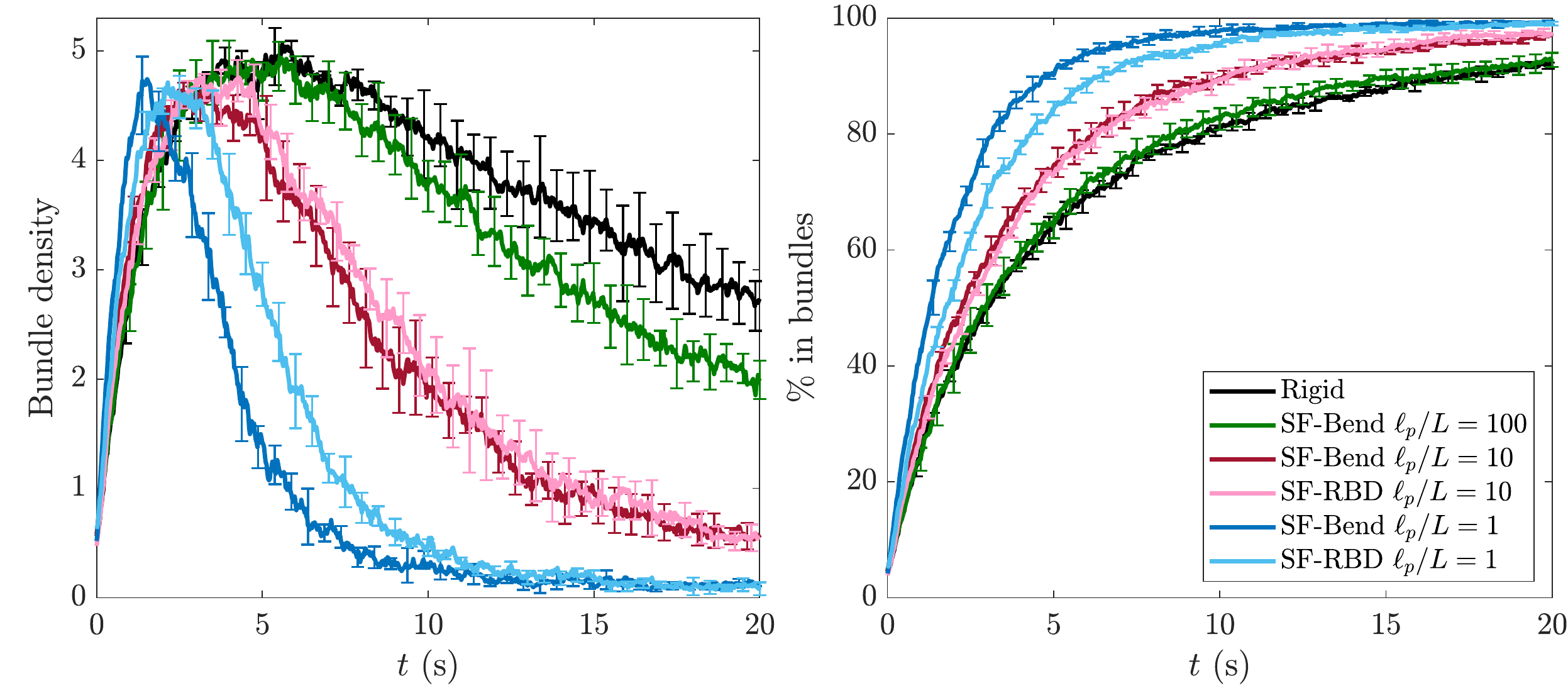}
\caption{\label{fig:BundleStats} Quantifying the dynamics of bundling in semiflexible fiber networks, and determining the role of bending fluctuations. We use the bundle density (left, number of bundles per unit volume) and percent of fibers in bundles (right) as metrics for bundling dynamics (see \cite{ActinBundlingDynamics} for discussion of these). We compare the semiflexible filaments with $\ell_p/L=1$ (blue), $\ell_p/L=10$ (red), and $\ell_p/L=100$ (green) to the rigid filament networks we simulated in \cite{ActinBundlingDynamics} (black), finding that $\ell_p/L=100$ is ``stiff enough'' to be rigid. We then separate the effect of bending fluctuations from CL forces and rigid body diffusion by comparing our semiflexible filament simulations (SF-Bend) to filaments which only fluctuate by rigidly rotating and translating (SF-RBD, lighter colors, see text for details). The CL forces and rigid body fluctuations are sufficient to account for the speed-up in bundling when $\ell_p/L = 10$, but not when $\ell_p/L=1$.}
\end{figure}

In this work, we can finally consider the full system, and how the third possible motion, transverse bending fluctuations, impacts bundling. It will be important, however, to separate these fluctuations from the rotational and translational diffusion that we have considered previously. Since our temporal integrator in Section\ \ref{sec:OneSolveMP} does not distinguish between the two kinds of fluctuations, for comparison we consider an alternative model where the only fluctuations are translational and rotational diffusion, as if the fibers were rigid in their current configuration. That is, instead of\ \eqref{eq:ItoCL}, we consider the dynamics
\begin{gather}
\label{eq:ItoCL2}
\dt{ \V{X}} = \M{K}\M{N}\M{K}^T\left(-\M{L}\V{X} +\FCL\right)+  \sqrt{2 k_B T}\M{K}_r\M{N}_r^{1/2}\Wproc,
\end{gather}
where $\M{K}_r$ is the kinematic matrix for a \emph{rigid} fiber which acts on a 6-vector $\left(\V{\Omega},\Ump\right)$ to give velocity in an analogous way to\ \eqref{eq:KNp1}. The $6 \times 6$ rigid-body mobility is\footnote{The psuedo-inversion in the calculation of $\M{N}_r$ is problematic when the fibers are nearly straight and there are eigenvalues near, but not equal to, zero in the resistance matrix $\M{K}_r^T \Mfor^{-1} \M{K}_r$. We find that the dynamics for SF-RBD filaments with $\ell_p/L =100$ are quite sensitive to the tolerance we use, and so we do not report them here. We find it sufficient to only consider up to $\ell_p/L=10$, for which we obtain dynamics that are not sensitive to the tolerance (see Fig.\ \ref{fig:BundleStats}). For $\ell_p/L=100$, fibers with semiflexible bending fluctuations behave almost identically to rigid fibers (see Fig.\ \ref{fig:BundleStats}), which we can simulate without difficulty \cite{ActinBundlingDynamics}.} $\M{N}_r =\left(\M{K}_r^T \Mfor^{-1} \M{K}_r\right)^\dagger$, exactly as in\ \eqref{eq:Ndef}. Our rationale for\ \eqref{eq:ItoCL2} is that if such semiflexible rigid body motion (``SF-RBD'') simulations give the same results as those with semiflexible bending fluctuations (``SF-Bend''), then bending fluctuations are not important to the bundling process. 

To advance the dynamics of SF-RBD filaments, we first take a step of the stochastic simulation algorithm for the CLs (as for SF-Bend fibers), then solve\ \eqref{eq:ItoCL2} using a splitting scheme. The splitting scheme is to first perform a rigid body rotation and translation to add the random term in\ \eqref{eq:ItoCL2} (see \cite[Eq.~(15)]{ActinBundlingDynamics}), then compute the cross-linking forces and perform a deterministic saddle point solve to capture the deterministic term in\ \eqref{eq:ItoCL2}. 

We begin by looking at the differences in bundling dynamics between \emph{semiflexible} filaments and the rigid filaments we considered in \cite{ActinBundlingDynamics}. This is actually the same test we performed in \cite{ActinBundlingDynamics}, but this time we consider filament fluctuations in addition to CL forces. The bundle density (number of bundles divided by periodic cell volume) and percent of fibers in bundles over time are shown using darker colors in Fig.\ \ref{fig:BundleStats}, where we consider $\ell_p/L=1$, 10, and 100, and compare the results to the case when the fibers are actually rigid, the dynamics of which are governed by the overdamped Ito Langevin equation
\begin{gather}
\label{eq:ItoCL3}
\dt{ \V{X}} = \M{K}_r\M{N}_r\M{K}_r^T\left(-\M{L}\V{X} +\FCL\right)+  \sqrt{2 k_B T}\M{K}_r\M{N}_r^{1/2}\Wproc.
\end{gather}
We observe almost complete overlap between the trajectory for $\ell_p/L=100$ and that for rigid fibers, (except for the bundle density curve at late times, which is when cross linkers exert extremely large forces on the filaments). This is not a surprise given that we see bundles of \emph{straight} filaments when $\ell_p/L=100$ in Fig.\ \ref{fig:BundlePic}, but it does showcase the ability of our temporal integrator to remain accurate in the stiff limit (this is not the case for SF-RBD dynamics, as discussed in footnote 7). Dropping to $\ell_p/L=10$, the pictures in Fig.\ \ref{fig:BundlePic} show curved bundles, and the plots in Fig.\ \ref{fig:BundleStats} for $\ell_p/L=10$ show significant deviations from rigid fibers, especially in the later stages of bundling. The curves for $\ell_p/L=1$ do not even match rigid fibers at early times, which indicates that the semiflexible fluctuations impact the first stage of bundle formation, in contrast to larger persistence lengths where the fluctuations appear to only accelerate later stages of bundle agglomeration. Since filaments are weakly cross-linked at early times, these results suggest that semiflexible bending fluctuations are accelerating the bundling process when $\ell_p/L \sim 1$. For $\ell_p/L \gtrsim 10$, the deviations from rigid fibers come only when the fibers are strongly cross linked, suggesting that CL forces combine with fiber flexibility to accelerate bundling.

To make this statement more precise, we compare simulations with (SF-Bend) and without (SF-RBD) thermal bending forces using lighter colors in Fig.\ \ref{fig:BundleStats}. When $\ell_p/L = 10$, we see identical dynamics between SF-Bend and SF-RBD filaments, which means that bending fluctuations contribute minimally to the bundling process for $\ell_p/L \geq 10$, and demonstrates that the curvature of the bundles we see in Fig.\ \ref{fig:BundlePic} when $\ell_p/L=10$ is indeed driven primarily by cross linking forces. When $\ell_p/L=1$, by contrast, we see faster bundling dynamics with SF-Bend filaments than with SF-RBD filaments, and we also see bundles in Fig.\ \ref{fig:BundlePic} that appear to have wavy spatial shapes. This implies that thermal bending fluctuations can accelerate bundling both in the initial and later stages, but only when the persistence length is comparable to the contour length of the fiber, in which case the transverse fluctuations effectively increase the probability that a CL (which can only stretch a finite amount) can bind two filaments. However, since actin filaments have persistence length on the order 10--20 $\mu$m \cite{gittes1993flexural}, we can conclude that the bundling dynamics of filaments with length 1--2 $\mu$m are \emph{not} significantly impacted by semiflexible bending fluctuations.


\section{Conclusions}
This paper represents a first step in applying spectral methods to fluctuating inextensible filaments in Stokes flow. While the advantages of spectral methods for simulating \emph{smooth} fibers are well known \cite{trefethen2000spectral, ehssan17, FibersWeakInextensibility}, this paper is to our knowledge the first attempt to use Chebyshev polynomials to simulate fibers that are inherently nonsmooth due to Brownian bending fluctuations. Our main motivation for doing this comes from hydrodynamics: our slender body quadrature scheme developed in \cite{TwistBend} allows us to compute the hydrodynamic mobility using $\mathcal{O}(1)$ points per filament (with respect to the aspect ratio $\epsc$), as opposed to the $1/\epsc$ points required in traditional blob-link (bead-link) methods \cite{ttbring08, kallemov2016immersed, keavRPY}. The quadrature scheme works on a spectral grid, since it requires a global interpolating function for the fiber centerline. Thus to include fluctuations we postulated the Gibbs-Boltzmann distribution\ \eqref{eq:GBDist} on the spectral grid of $N$ nodes, then constructed spatially-discrete overdamped Langevin dynamics that is in detailed balance with respect to that distribution. 

We showed through a series of equilibrium and non-equilibrium tests that spectral methods can be advantageous in the regime where the fibers are semiflexible $\ell_p \gtrsim L$ and slender $\epsc =\rc/L \ll 1$, which is therefore the regime where the persistence length (smallest lengthscale on which fluctuations are visible) is much larger than the fiber radius. Since our primary interest is in relatively stiff (and quite slender) fibers like actin, there is some promise for spectral methods to give good approximations of blob-link chains by faithfully modeling the hydrodynamics with fewer degrees of freedom. In this case, the fluctuations on lengthscales $\rc$, which are smoothed out by our spectral method, do not impact the fiber dynamics (other than on very short timescales) or equilibrium statistics, and so we are able to approximate both well using a small number of Chebyshev nodes. Our tests on dynamics showed that the error in the rate of relaxation of a straight chain for $N \geq 24$ is small and dominated by the error from \emph{equilibrium} statistical mechanics, which means that the main driver of the difference between the blob-link and spectral methods comes from the difference between our coarse-grained Gibbs-Boltzmann distribution\ \eqref{eq:GBDist} and the GB distribution of a refined blob-link chain. Nevertheless, for $\ell_p/L=1$ and $\epsc  \approx 10^{-2}$, we showed that even $N=12$ provides a reasonable approximation to the relaxation dynamics of a stretched chain, with a time step size two orders of magnitude larger than the corresponding blob-link discretization with 100 blobs. This sort of spectral coarse-grained (fully discrete rather than continuum) representation of semiflexible fibers can allow for large-scale simulations of networks of fibers over physically relevant timescales, unlike fully or finely resolved blob-link models, \new{which are better suited for simulations where the small $\mathcal{O}(\eps)$ lengthscales need to be resolved.} 

Our approach to the spatial discretization as presented in Section\ \ref{sec:KinInex} seems straightforward, but was actually the product of several months of trial and error. The reason is that the nonsmoothness of thermal fluctuations exposes all of the weaknesses of a particular spectral discretization, especially when that discretization fails to properly track high-frequency modes. It is therefore vital that our discretization be free from aliasing errors. In particular, it is important that we be able to recover the original tangent vector rotation rates from the fiber velocity field by ``inverting'' the kinematic matrix, and thereby develop temporal integrators that capture the right drift terms in expectation. An open question is whether this level of scrupulousness, which is also necessary because Chebyshev polynomials cannot define everywhere inextensible curves, could have been avoided by using a different (rapidly converging) series representation for the fiber centerline. Other, non-Chebyshev, representations for (everywhere) inextensible filaments like curvature-torsion require higher order derivatives of $\Xs=\ds{\V{X}}$ which have non-decaying spectra, since for Brownian filaments $\Xs$ is almost everywhere differentiable but its derivative is white noise. An optimal representation would be in terms of some basis functions whose contributions to the fiber shape decay fast, especially for large persistence lengths. \new{Such a representation would ideally nearly diagonalize the bending elasticity operator and be everywhere inextensible, the latter being vital to remove ambiguity in how to represent the fiber (tangent vectors vs.\ positions).}

We showed that our spectral method with slender body RPY quadrature can reduce the number of nodes required to resolve hydrodynamics, and also substantially increase the time step size required for accuracy relative to both blob-link methods and spectral methods with other hydrodynamic models (such as local slender body theory). This is an important feature because our midpoint temporal integrator, which is based on taking a Fixman-like predictor step to the midpoint to capture the drift term in expectation, can only generate accurate equilibrium statistics and dynamics when the time step is sufficiently small. By breaking the dynamics into a set of modes, we showed that to generate the correct equilibrium statistics the temporal integrator must resolve a number of modes that scales linearly with $N$. Seeing as the timescales of the modes scale like $k^{-4}$ for large node index $k$ (see Fig.\ \ref{fig:Timescales}), this means that our time step size scales roughly as $N^{-4}$, which means we must hold the number of modes down if we want to simulate with a reasonable time step size. An alternative approach would be to try to increase the required time step size by using an exponential integrator. For rotations on the unit sphere, this requires Lie integrators like that used by Schoeller et al.\ in \cite{keavRPY}. Stochastic exponential Lie integrators have not yet been developed (to our knowledge), and so they represent an interesting avenue of exploration for future work. 

We ended this paper by applying our algorithm to the bundling dynamics of transiently-cross linked fluctuating actin filaments. When the fibers do not turnover, we showed previously \cite{ActinBundlingDynamics} that there is nothing stopping them from forming large bundles through a thermal ratcheting mechanism where cross linkers (CLs) stretch to bind filaments, then contract to pull them closer together and allow additional CL bindings. Here we showed that the bundling of fluctuating actin filaments, which have persistence length on the order 10 $\mu$m \cite{gittes1993flexural}, significantly differs from that of rigid fibers, which is in contrast to the conclusion we drew previously when filaments did not fluctuate (and only moved by CL forces) \cite{ActinBundlingDynamics}. However, when $\ell_p/L \gtrsim 10$, we showed that these differences can be explained by a reduced model where actin filaments rotate and translate randomly, and bend via CL forces. In this case the combination of translational and rotational rigid-body diffusion with fiber flexibility and CL forces drives the observed differences from rigid fiber dynamics. Semiflexible bending fluctuations only accelerate bundling dynamics when $\ell_p/L \approx 1$, in which case the differences from rigid fibers can be observed even in the regime when the fibers are weakly cross linked.

Our study of bundling dynamics, like our previous work \cite{ActinBundlingDynamics}, did not consider hydrodynamic interactions between filaments. In previous work \cite{ActinCLsRheology}, we showed that such interactions slow down the bundling process for non-fluctuating fibers. An interesting extension of the work here would be to see if such a conclusion also holds for fluctuating filaments. Simulating hydrodynamic interactions with fluctuations is in fact possible in (log) linear time with the midpoint temporal integrator developed here once suitable iterative solvers are developed. For the blob-link discretization, this has already been done by some of us for triply periodic systems and will be shared in separate work. For spectral methods, our previous work \cite{FibersWeakInextensibility}, which was restricted to the deterministic setting, relied on splitting the mobility into intra- and inter- fiber hydrodynamics, then computing the second piece fast using oversampled quadrature on a GPU by the Positively Split Ewald (PSE) method \cite{PSRPY, ActinCLsRheology}. This approach might break down for fluctuating fibers unless we can guarantee that each piece is separately symmetric positive definite (SPD) \cite{PSRPY}. A simple solution would be to use oversampled quadrature for the mobility on \emph{all fibers} (which Appendix\ \ref{sec:QuadOS} shows is what our special quadrature approximates), so that $\M{M}$ and its square root can be applied fast using the PSE method. Other possible approaches are to use a Galerkin approach, which we have shown is consistent with the discretization presented here (see\ \eqref{eq:Mref}), to define a consistent SPD mobility matrix, or use the positive splitting of the RPY kernel in the PSE method to split the action of the mobility into a near field special quadrature and a far-field fast summation. \new{After we formulate the mobility, it will be interesting to see how many tangent vectors $N$ are required to resolve the hydrodynamics, since in the case of bundled suspensions the lengthscale of filament interaction decreases over time. Presumably more than $N=12$ tangent vectors would be required to resolve the individual hydrodynamic interactions, but whether the error made in using a small number of points affects the overall suspension behavior (e.g., the bundling time) significantly is hard to predict.}

From an application standpoint, the rheology of transiently cross-linked actin networks with fiber turnover is an exciting open problem in biological physics. Recent experiments have shown that the stress relaxation in the actin cytoskeleton is slow, with a viscous modulus that scales as $\omega^{1/2}$ on long timescales \cite{yao2011JMB}. Since then, there have been a number of physical theories seeking to explain this behavior \cite{mulla2019origin, muller2014PRL, broedersz2010PRL,  yao2013PRL}, all of which are based to some degree on a coarse-grained relationship between the cross linker distribution and the stress relaxation timescales. Thus the ability to simulate many-fiber suspensions for long times is vital to justify more rigorously the assumptions made by continuum theories, and study the role of Brownian motion and bending fluctuations in stress relaxation in the actin cytoskeleton. For the bundled networks in Section\ \ref{sec:Bundle}, we have already shown \cite{ActinCLsRheology} that both intra- and inter-fiber hydrodynamic interactions change the viscoelastic moduli in the deterministic setting, and it is therefore interesting to see how or if this conclusion changes with fluctuations. Doing these calculations will require care, however, since in rheology we are looking for \emph{quantitative} values of stress, which includes a contribution from Brownian motion \cite[Eq~(3.169)]{doi1988theory} that needs to be formulated carefully for overdamped inextensible fibers, and computed efficiently numerically. We have already seen how fluctuating fiber numerics can be quite sensitive to spatial and temporal discretizations. When we combine this with other (non-equilibrium) dynamics like cross linker attachment/detachment and oscillatory shear flow, we will need extra care to ensure that the results are not dominated by discretization artifacts. The stress also has a contribution from steric interactions, which are neglected here but can be added in future work by approximating each fiber with $\mathcal{O}(1)$ spherocylinders \cite{chen2021rheology} which resolve the smallest lengthscale on which the fiber bends (roughly $\ell_p$). Robust quantitative results will give us the ability to connect microscale Brownian filament hydrodynamics to macroscopic cytoskeletal behavior.

\section*{Acknowledgments}
We thank Eric Vanden-Eijnden, Pep Español, Ravi Jagadeeshan, Isaac Pincus, and Hans C.\ Ottinger for useful discussions about the continuum vs.\ discrete inextensible chain. We also thank Raul P.\ Pelaez for help with the implementation of the blob-link algorithm, and Alex Mogilner for reviewing the section on actin bundling dynamics. Ondrej Maxian is supported by the NYU Dean's Dissertation Fellowship and the NSF via GRFP/DGE-1342536. This work was also supported by the National Science Foundation through Division of Mathematical Sciences award DMS-2052515, and through a Research and Training
Group in Modeling and Simulation under award RTG/DMS-1646339.

\subsection*{Data availability}
All of the python/C++ codes and corresponding input files to reproduce our results
can be found at \url{https://github.com/stochasticHydroTools/SlenderBody}.

\subsection*{Author declarations}
The authors have no conflicts to disclose.

\begin{appendices}
\setcounter{equation}{0}
\renewcommand{\theequation}{\thesection.\arabic{equation}}
\section{Saddle point system as gradient descent dynamics \label{sec:SPLag}}
In this appendix, we show that the deterministic dynamics\ \eqref{eq:SP} take the form of a constrained gradient descent. Let us consider a backward Euler integrator for\ \eqref{eq:SP} and the Lagrangian 
\begin{gather}
\label{eq:LagSP}
\mathcal{L}\left[\V{X},\V{\alpha},\V{\lambda}\right] = \frac{1}{2} \V{X}^T \M{L} \V{X} + \frac{1}{2 \D t} \left(\V{X}-\tdisc{\V{X}}{n}\right)^T \left(\tdisc{\Mfor}{n}\right)^{-1}\ \left(\V{X}-\tdisc{\V{X}}{n}\right)\\ \nonumber +\V{\lambda}^T \Wt \left( \tdisc{\M{K}}{n}\V{\alpha}-\frac{1}{\D t} \left(\V{X}-\tdisc{\V{X}}{n}\right)\right).
\end{gather}
Here the matrices $\M{K}$ and $\Mfor$ are evaluated at time $n$, and consequently are constant with respect to $\V{X}$. In\ \eqref{eq:LagSP}, the second term is $\D t/2$ times the rate of dissipation in the fluid (velocity $\times$ force), and $\V{\lambda}$ represents a force density, which is a Lagrange multiplier for the inextensibility constraint. We formulate that constraint as an $L^2$ inner product using the discrete weights matrix $\Wt$, since this gives the work done by the constraint force density on the fluid. 

We now differentiate the Lagrangian with respect to the three inputs to arrive at our equations of motion
\begin{gather}
\label{eq:diffLX}
\frac{\delta \mathcal{L}}{\delta \V{X}} = \V 0 \rightarrow -\M{L}\V{X}+\Wt \V{\lambda} = \frac{1}{\D t} \left(\tdisc{\Mfor}{n}\right)^{-1} \left(\V{X}-\tdisc{\V{X}}{n}\right)\\
\label{eq:diffLlam}
\frac{\delta \mathcal{L}}{\delta \V{\lambda}} = \V 0 \rightarrow \tdisc{\M{K}}{n}\V{\alpha} = \frac{1}{\D t} \left(\V{X}-\tdisc{\V{X}}{n}\right) \\
\frac{\delta \mathcal{L}}{\delta \V{\alpha}} = \V 0 \rightarrow \left(\tdisc{\M{K}}{n}\right)^T \Wt \V{\lambda} = \V 0
\end{gather}
Combining\ \eqref{eq:diffLX} and\ \eqref{eq:diffLlam} and taking $\D t \rightarrow 0$, we arrive at the deterministic saddle point system
\begin{gather}
\label{eq:SPInex}
\M{K}\V{\alpha} = \Msym \left(-\M{L}\V{X}+\V{\Lambda}\right) \\
\label{eq:KTLam}
\M{K}^T \V{\Lambda} = \V 0.
\end{gather}
Here we have set the force $\M{\Lambda} = \Wt \V{\lambda}$, which represents another instance of the conversion between force and force density defined in\ \eqref{eq:FDenFromF}. In this case it is easy to see how this conversion arises from the $L^2$ inner product. Equation\ \eqref{eq:KTLam} is the principle of virtual work invoked in our previous work \cite[Sec.~3.4]{FibersWeakInextensibility}.

\subsection{Adding fluctuations}
With fluctuations, the same variational technique can be used to derive the saddle point system\ \eqref{eq:SPB1} by adding the work done (entropy dissipated) by the random force $\V{F}_\text{stoch}$
$$\left(\frac{\V{X}-\tdisc{\V{X}}{n}}{\D t}\right)^T\V{F}_\text{stoch}\D t = \left(\frac{\V{X}-\tdisc{\V{X}}{n}}{\D t}\right)^T\sqrt{\frac{2 k_B T}{\D t}}\left(\tdisc{\Mfor}{n}\right)^{-1/2} \tdisc{\gauss}{n}\D t,$$
to the Lagrangian\ \eqref{eq:diffLX}. However, solving the system\ \eqref{eq:SPB1} produces only the first and last term of the Langevin equation\ \eqref{eq:ItoX}, and not the additional stochastic drift terms, which are required to ensure time-reversibility (detailed balance). So far, we have not been able to produce those terms via a variational argument.

\section{RPY Mobility \label{sec:RPYAppen}}
\setcounter{figure}{0}    
\renewcommand{\thefigure}{B\arabic{figure}}
\setcounter{equation}{0}
In this appendix, we give a more detailed discussion of the RPY mobility, beginning with the definition of the mobility for two particles. We then give a few details on modifications of our ``slender-body'' quadrature scheme \cite{TwistBend} for nonsmooth filaments. Specifically, since our quadrature scheme was developed with smooth filaments in mind, we assumed that the tangent vectors obtained by differentiating $\XPoly$ on the $N_x$ point grid have unit length. This is not necessarily the case for random filaments, as the tangent vectors are unit length only when confined to the $N$ collocation points tracking $\Xs$. We therefore have to make a few small modifications, as detailed in Appendix\ \ref{sec:QuadMod}. In Appendix\ \ref{sec:MatrixSym}, we consider the error incurred in symmetrizing the mobility\ \eqref{eq:Msym} and in truncating its negative eigenvalues. We do this by looking at the errors with respect to the symmetric positive definite upsampled mobility\ \eqref{eq:Mref}. Finally, in Appendix\ \ref{sec:AdvSpec}, we compare our special quadrature mobility to other mobility options for the relaxing filament in Section\ \ref{sec:Relax}, showing that there are signficant advantages to using special quadrature both from a spatial and temporal perspective.

\subsection{RPY kernel \label{sec:RPYkern}}
To define the RPY kernel between two points $\V{x}$ and $\V y$, we first let $\V r = \V x - \V y$ and $r = \norm{\V r}$ with $\hat{\V r}=\V r/r$. The RPY kernel is based on linear combinations of the Stokeslet and doublet, which are defined as 
\begin{equation}
\label{eq:Slet}
\Slet{\V x, \V y} = \EPMI \left(\frac{\M{I}+\hat{\V{r}}\hat{\V{r}}}{r}\right) \qquad \text{ and } \qquad 
\Dlet{\V x, \V y}  = \EPMI \left(\frac{\M{I}-3\hat{\V{r}}\hat{\V{r}}}{r^3}\right). 
\end{equation}
The RPY kernel for an unbounded fluid is then given by  \cite{rpyOG, wajnryb2013generalization, PSRPY}
\begin{gather}
\label{eq:rpyknel}
 \Mfor_\text{RPY}\left(\V{x},\V{y}; \eps\right) = 
\begin{cases} 
\displaystyle
\Slet{\V x, \V y}+\frac{2\eps^2}{3} \Dlet{\V x, \V y}& r > 2\eps\\
\displaystyle
\EPMI \left[\left(\frac{4}{3\eps}-\frac{3r}{8\eps^2}\right)\M{I}+\frac{1}{8\eps^2 r}\V{r}\V{r}\right]&  r \leq 2 \eps
\end{cases}.
\end{gather}

\subsection{Modifications to slender body quadrature for nonsmooth filaments \label{sec:QuadMod}}
A subtle point in computing the mobility on the grid of size $N_x$ is that the quadrature schemes we use to compute\ \eqref{eq:McRPYdef} are based on having unit-length tangent vectors; i.e., $\left(\M{D}\V{X}\right)_\Bind{p} \cdot \left(\M{D}\V{X}\right)_\Bind{p}=1$ for all $p$ on the grid of size $N_x$. For random filaments, the length of the tangent vectors is 1 on the $N$ point grid, but might differ substantially (by as much as 50\% or so) when these polynomials are resampled to the $N_x=N+1$ point grid. As such, we need to modify the quadrature schemes in \cite[Appendix~G]{TwistBend} to obtain a quadrature scheme that converges to the continuum integral\ \eqref{eq:McRPYdef} as the number of points becomes large. 

This is quite straightforward to do. In \cite[Appendix~G]{TwistBend}, we break the RPY integral into three pieces: the Stokeslet on $|s-s'| > 2\eps$, the doublet on $|s-s'| > 2\eps$ and the piece for $|s-s'| \leq 2\eps$. The last piece is done by direct Gauss-Legendre quadrature, so there are no changes to that in this case. The first two pieces are done using special quadrature schemes, which involve subtracting the leading order singularity and integrating what remains. The leading order singularities in the case of a tangent vector that does not have norm one are
\begin{equation*}
\M{S}_1(s,s') = \EPMI \left(\frac{\M{I}+\widehat{\ds{\XPoly}}(s)\widehat{\ds{\XPoly}}(s)}{\norm{\ds{\XPoly}(s)}|s-s'|}\right) \qquad \M{D}_1(s,s') = \EPMI \left(\frac{\M{I}-3\widehat{\ds{\XPoly}}(s)\widehat{\ds{\XPoly}}(s)}{\norm{\ds{\XPoly}(s)}^3|s-s'|^3}\right),
\end{equation*}
where the hat denotes a normalized vector. Thus the remaining integrals are
\begin{gather}
\nonumber
\tt{\V{U}}^\text{(int, S)}(s)=\int_{D(s)} \left(\Slet{\XPoly(s),\XPoly(s')} \V{f}\left(s'\right) -  \M{S}_1(s,s')\V{f}(s)\right) \, ds' ,\\ \nonumber
\tt{\V{U}}^\text{(int, D)}(s)=\int_{D(s)} \left(\Dlet{\XPoly(s),\XPoly(s')} \V{f}\left(s'\right) -\M{D}_1(s,s') \V{f}(s)\right) \, ds'.
\end{gather}
Here $D(s)$ is the domain on which $|s-s'| \geq 2\eps$ for each $s$, defined in \cite[Eq.~(G.2)]{TwistBend}. 

In \cite[Appendix~G]{TwistBend}, we develop a special quadrature scheme which is based on putting the integrals in the form
\begin{gather} 
\nonumber
 \tt{\V{U}}^\text{(int, S)}(s) =\int_{D(s)} \V{g}_\text{tt}(s,s') \frac{(s'-s)}{|s'-s|} \, ds' \\ \nonumber
 \tt{\V{U}}^\text{(int, D)}(s) =\int_{D(s)} \V{g}_D(s',s)  \frac{(s'-s)}{|s'-s|^3} \,  ds'.
\end{gather}
Our special quadrature scheme can be implemented as described there, but with the following modified definitions for tangent vectors that do not have norm one, 
\begin{gather}
\nonumber
\V{g}_\text{tt}(s,s') = \left(\Slet{\XPoly(s),\XPoly(s')} \V{f}\left(s'\right) -\M{S}_1(s,s')\V{f}(s) \right)\frac{|s'-s|}{s'-s} \\ \nonumber
\V{g}_D(s',s)=\left(\Dlet{\XPoly(s),\XPoly(s')} \V{f}\left(s'\right) -\M{D}_1(s,s')\V{f}(s) \right) \frac{|s'-s|^3}{(s'-s)}.
\end{gather}
These functions are nonsingular at $s=s'$, with the finite limits 
\begin{gather*}
\lim_{s' \to s} \V{g}_\text{tt}(s,s') = \EPMI  \Bigg{[} \frac{1}{2\norm{\ds{\XPoly}(s)}^3} \bigg{(} \ds{\XPoly}(s)\ds^2{\XPoly}(s)+\ds^2{\XPoly}(s)\ds{\XPoly}(s)\\ -\left(\ds{\XPoly}(s) \cdot \ds^2{\XPoly}(s)\right) \left(\M{I}+3\widehat{\ds{\XPoly}}(s)\widehat{\ds{\XPoly}}(s)\right)\bigg{)}\V{f}(s) + \left(\frac{\M{I}+\widehat{\ds{\XPoly}}(s)\widehat{\ds{\XPoly}}(s)}{\norm{\ds{\XPoly}(s)}}\right) \ds{\V{f}}(s) \Bigg{]}\\ 
\lim_{s'\rightarrow s} \V{g}_D(s',s) = \EPMI  \Bigg{[} \frac{1}{2\norm{\ds{\XPoly}(s)}^5} \bigg{(} -3\left(\ds{\XPoly}(s)\ds^2{\XPoly}(s)+\ds^2{\XPoly}(s)\ds{\XPoly}(s)\right)\\ -\left(\ds{\XPoly}(s) \cdot \ds^2{\XPoly}(s)\right) \left(3\M{I}-15\widehat{\ds{\XPoly}}(s)\widehat{\ds{\XPoly}}(s)\right)\bigg{)}\V{f}(s) + \left(\frac{\M{I}-3\widehat{\ds{\XPoly}}(s)\widehat{\ds{\XPoly}}(s)}{\norm{\ds{\XPoly}(s)}^3}\right) \ds{\V{f}}(s) \Bigg{]}.
\end{gather*}

\subsection{Eigenvalues of slender body quadrature mobility matrix \label{sec:MatrixSym}}
In this part of the appendix, we look in more detail at how to determine a systematic threshold for eigenvalue truncation of our slender body quadrature. We do this by comparing the eigenvalues of the symmetrized mobility\ \eqref{eq:Msym} with the reference mobility\ \eqref{eq:Mref}. We will use $N_u=1000$ points for the reference mobility, having verified that using even more reference points gives the same results.

\begin{figure}
\centering
\includegraphics[width=\textwidth]{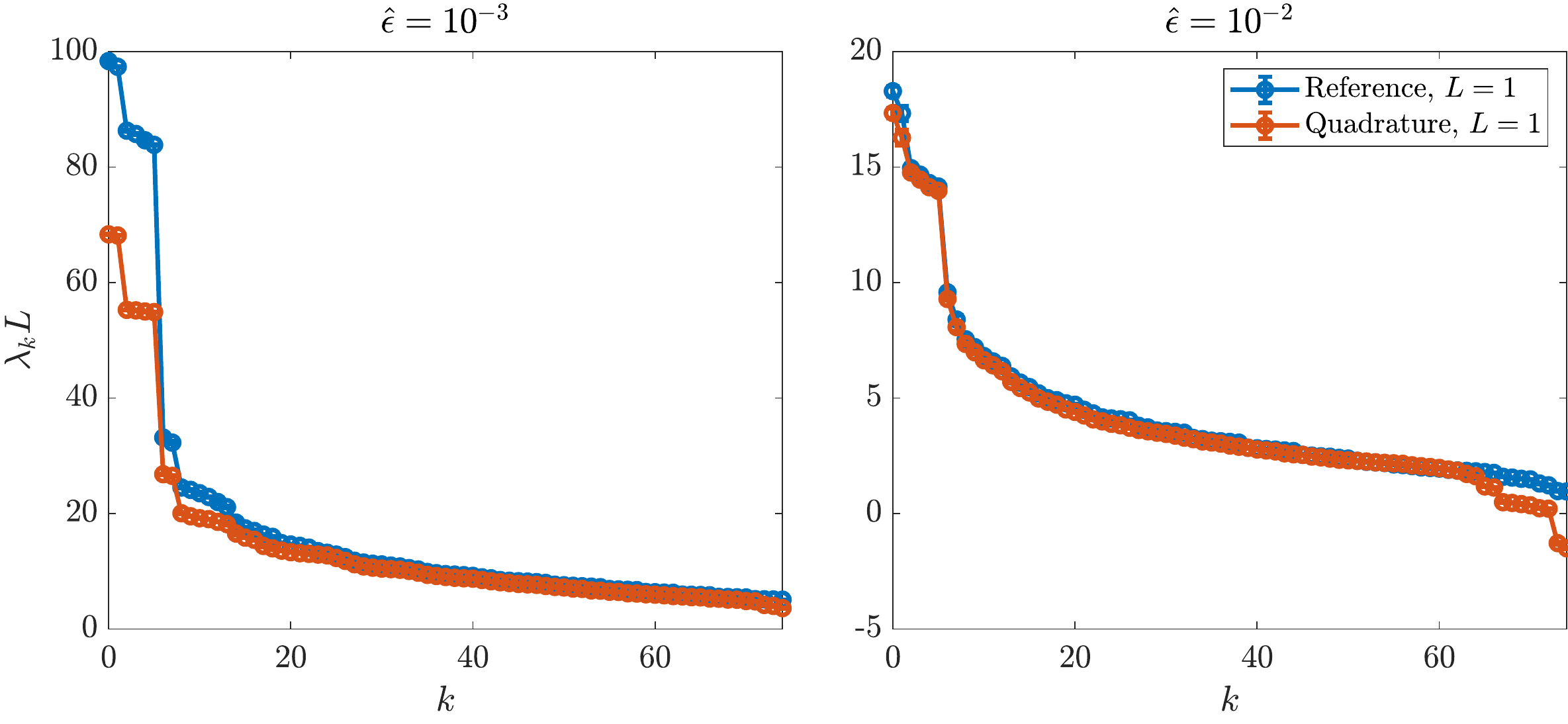}
\caption{\label{fig:EigvalsRQ} Eigenvalues of $\Mfor$ for 10 random filaments, sampled from the equilibrium distribution for $\ell_p/L=1$ and $N=24$ (here $\mu=1$). We show the eigenvalues of the reference mobility\ \eqref{eq:Mref} in blue and those of the symmetrized quadrature mobility\ \eqref{eq:Msym} in red. The largest eigenvalues are incorrect for small $\epsRS$ because they correspond to highly nonsmooth modes which are flat in the fiber interior and large at the endpoints. For, $\epsRS=10^{-2}$, we start to see negative eigenvalues in the symmetrized mobility\ \eqref{eq:Msym}. }
\end{figure}

To motivate the truncation of the negative eigenvalues of\ \eqref{eq:Msym}, we first compare the eigenvalues of\ \eqref{eq:Msym} (without any truncation) to those of the reference mobility\ \eqref{eq:Mref}. To do this, we fix $N=24$ and $\ell_p/L=1$ and consider a set of ten samples from the equilibrium distribution studied in Section\ \ref{sec:EEMCMC}. The mean eigenvalues of $\MforRef$ and our approximate $\Msym$ are shown in Fig.\ \ref{fig:EigvalsRQ}, where we have repeated the test three times to generate error bars. Fortunately, as we've observed previously \cite[Sec.~4.4.1]{TwistBend}, the eigenvalues of the mobility are not particularly sensitive to the fiber shape; hence the error bars in Fig.\ \ref{fig:EigvalsRQ} are smaller than the symbol size. For $\epsRS=10^{-3}$, we observe large errors in the largest eigenvalues, which correspond to modes that are flat in the fiber interior and peak at the fiber endpoint. These modes are rapidly damped, and are associated with fast timescales. 

When we drop to $\epsRS=10^{-2}$, we begin to see negative eigenvalues in the quadrature-based mobility. These negative eigenvalues, which separate themselves from the positive eigenvalues of the reference mobility, are associated with highly-oscillatory modes which have a small Rayleigh quotient $\V{v}^T \MforRef \V{v}$ using the reference mobility. Because modes with negative eigenvalues in the quadrature mobility are high frequency, and therefore give small Rayleigh quotients in the reference mobility, a sensible way to truncate the eigenvalues of $\Msym$ in\ \eqref{eq:Msym} is to set any eigenvalues less than the smallest eigenvalue of $\MforRef$ to that smallest eigenvalue. That is, we post-process the eigenvalues of\ \eqref{eq:Msym} to set any initial $\lambda_k < \sigma$ to $\lambda_k=\sigma$, where $\sigma$ is \emph{the smallest eigenvalue of the reference mobility matrix} $\MforRef$ for that $N$, $L$, and $\epsRS$. \new{In Fig.\ \ref{fig:EigvalsRQ}, we do this by averaging over thirty samples, but, as discussed previously, this makes little difference, and the eigenvalues for a straight filament could also be used. Thus, in practice our python code sets the threshold $\sigma$ by computing the reference mobility\ \eqref{eq:Mref} for a straight filament with $1/\epsc$ upsampling points and finding the smallest eigenvalue of the resulting matrix.}

\subsection{Comparing special quadrature to other mobilities \label{sec:AdvSpec}}
This part of the appendix is concerned with the different options for the matrix $\Mfor$, which maps force to velocity in fluctuating fibers. We consider four possibilities:
\begin{enumerate}[wide, labelwidth=!, labelindent=0pt]
\item RPY special quadrature as discussed in Section\ \ref{sec:Mob}, in which we first compute the matrix $\M{M}$, then symmetrize $\M{M} \Wt^{-1}$ and truncate its negative eigenvalues. 
\item RPY oversampled quadrature on a grid of size $N_u$, which is the SPD matrix\ \eqref{eq:Mref}.
The main part of this computation is the action of the matrix $\Mfor_{\text{RPY}, u}$, which maps forces to velocities on the upsampled grid. 
\item Direct evaluation of the RPY kernel on the spectral grid, i.e., $\Mfor =  \Mfor_{\text{RPY}}$, which is automatically SPD. This makes the action of the mobility the same as that of the blob-link method, but with uneven spacing of the grid points. While this calculation does not formally see the Chebyshev weights if we work on forces directly, it can still be viewed as Clenshaw-Curtis quadrature on force \emph{density,} with the conversion given in\ \eqref{eq:FDenFromF}.
\item The local drag approximation to the matrix $\M{M}$, which is a block-diagonal matrix with $3 \times 3$ blocks
\begin{gather}
\label{eq:UinnerA}
8 \pi \mu \M{M}_\Bind{i,i} = \left( \frac{a_S(s_i)+a_\text{CLI}(s_i)}{\norm{\ds{\XPoly}(s_i)}}+ \frac{2 \eps^2}{3}\frac{a_D(s_i)}{\norm{\ds{\XPoly}(s_i)}^3}\right)\M{I}\\ \nonumber +\left( \frac{a_S(s_i)+a_\text{CLT}(s_i)}{\norm{\ds{\XPoly}(s_i)}} -2 \eps^2\frac{a_D(s_i)}{\norm{\ds{\XPoly}(s_i)}^3}\right)\Xs_\Bind{i}\Xs_\Bind{i},
\end{gather}
where $\Xs_\Bind{i}=\ds{\XPoly}(s_i)/\norm{\ds{\XPoly}(s_i)}$ and 
\begin{gather}
a_S(s) = \begin{cases}
\ln{\left(\dfrac{(L-s)s}{4\eps^2}\right)} & 2\eps < s < L-2\eps\\[4 pt]
\ln{\left(\dfrac{(L-s)}{2\eps}\right)}& s \leq 2\eps \\[4 pt]
\ln{\left(\dfrac{s}{2\eps}\right)} & s \geq L-2\eps
\end{cases} \\[4 pt] \nonumber
a_D(s) = \begin{cases}
\dfrac{1}{4 \eps^2}-\dfrac{1}{2s^2}-\dfrac{1}{2(L-s)^2}& 2\eps < s < L-2\eps\\[4 pt]
\dfrac{1}{8 \eps^2}-\dfrac{1}{2(L-s)^2}& s \leq 2\eps \\[4 pt]
\dfrac{1}{8 \eps^2}-\dfrac{1}{2s^2}& s \geq L-2\eps
\end{cases}\\[4 pt] \nonumber
a_\text{CLI}(s) = \begin{cases}
\dfrac{23}{6}& 2\eps < s < L-2\eps\\[4 pt]
\dfrac{23}{12}+\dfrac{4s}{3\eps} - \dfrac{3s^2}{16\eps^2}& s \leq 2\eps \\[4 pt]
\dfrac{23}{12}+\dfrac{4(L-s)}{3\eps} - \dfrac{3(L-s)^2}{16\eps^2}& s \geq L-2\eps
\end{cases}
a_\text{CLT}(s) = \begin{cases}
\dfrac{1}{2} & 2\eps < s < L-2\eps\\[4 pt]
\dfrac{1}{4}+\dfrac{s^2}{16\eps^2} & s \leq 2\eps \\[4 pt]
\dfrac{1}{4}+\dfrac{(L-s)^2}{16\eps^2}& s \geq L-2\eps
\end{cases}
\end{gather}
The derivation of this formula is given in \cite{TwistSBT}, but here we have made some modifications since the tangent vectors $\ds{\XPoly}(s_i)$ might not be exactly norm one when the fibers are Brownian. If $\norm{\ds{\XPoly}(s_i)}=1$, then this formula reduces to \cite[Eq.~(B7)]{TwistSBT}. We note that this local drag theory is qualitatively the same as the more commonly used $8\pi \mu \M{M}_\Bind{i,i}=\ln{\left(\epsc^{-2}\right)}\left(\M{I}+\Xs_\Bind{i} \Xs_\Bind{i} \right) + \left(\M{I}-3\Xs_\Bind{i} \Xs_\Bind{i}\right)$, which describes filaments with ellipsoidally-tapered radius functions. Since our filaments are not ellipsoidally-tapered, the RPY-based theory from \cite{TwistSBT} is quantiatively more accurate (the trends we describe about spatial and temporal convergence are unchanged when we use ellipsoidal SBT). As in step 1, we symmetrize the matrix $\Mfor = \M{M}\Wt^{-1}$ and truncate its negative eigenvalues (the matrix is not automatically symmetric) to obtain the force-velocity mobility matrix.
\end{enumerate}

Throughout this appendix, we will consider the test of an initially-straight single filament relaxing to its equilibrium fluctuations from Section\ \ref{sec:Relax}, utilizing\ \eqref{eq:tbar} for temporal and spatial rescaling (we will rescale each trajectory by its equilibrium mean from MCMC). We will first consider only trajectories that are converged in time, so that our goal is to isolate the difference between spatial discretizations. We will then zero in on potential temporal convergence problems for the mobilities that appear advantageous at first.

\subsubsection{Oversampling converges to special quadrature \label{sec:QuadOS}}
\begin{figure}
\centering
\includegraphics[width=\textwidth]{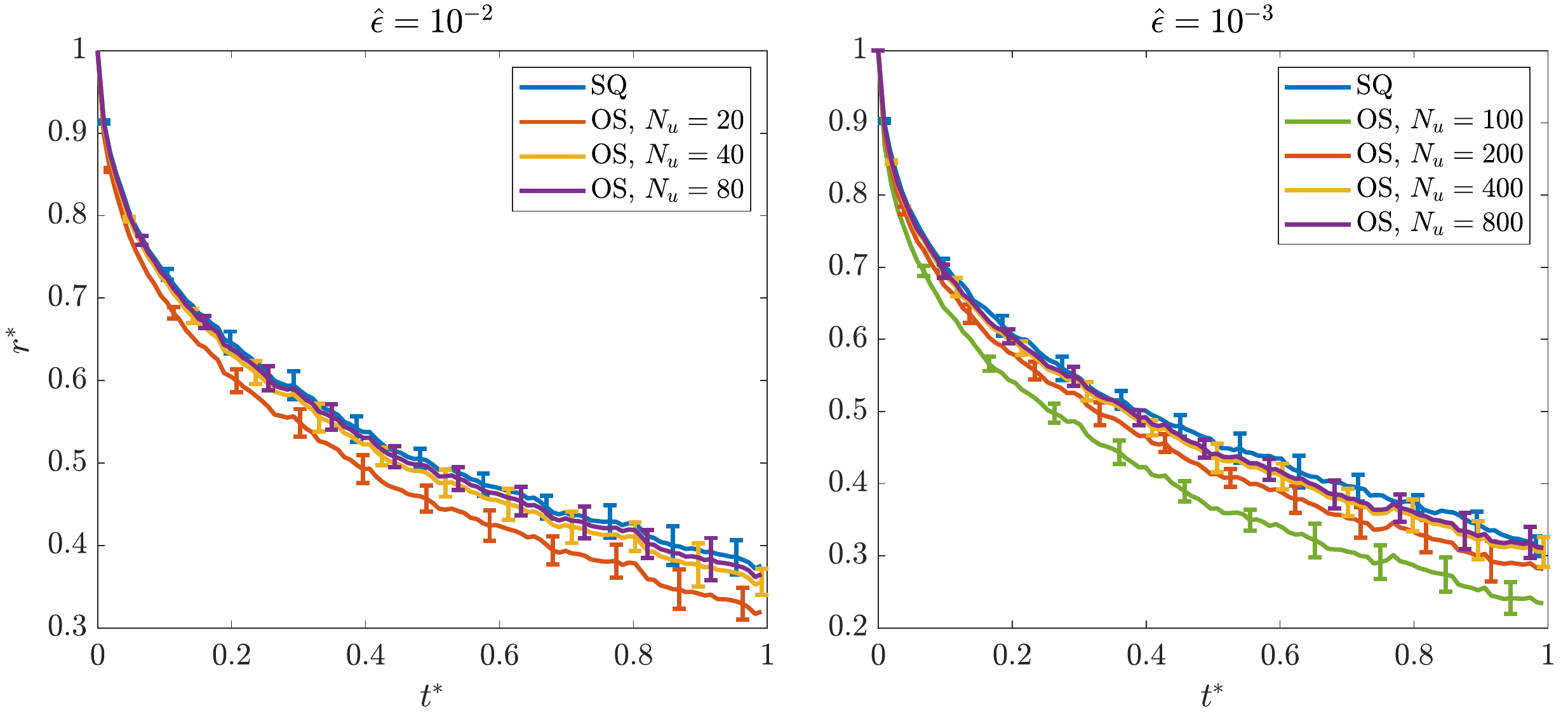}
\caption{\label{fig:OvSp} End-to-end distance trajectory comparing special quadrature (SQ) with oversampled quadrature (OS) with various numbers of oversampling points $N_u$, and a fixed number of $N=12$ collocation points. The left plot shows $\epsRS=10^{-2}$, for which we can use a small number of oversampled points. The more slender $\epsRS=10^{-3}$ at right requires many more oversampled points, but eventually converges to the special quadrature results. } 
\end{figure}

We can eliminate the ambiguity between special and oversampled quadrature by showing that trajectories using oversampled quadrature converge to those of symmetrized special quadrature as $N_u$ increases. Results for this are shown in Fig.\ \ref{fig:OvSp}, where we fix $N=12$ and show that the number of oversampled points $N_u$ required for the special quadrature to match oversampled quadrature depends on $\epsRS$. We require roughly $0.4/\epsRS$ points to obtain equivalence of the two trajectories: thus 40 points are acceptable for $\epsRS=10^{-2}$, which may make the special quadrature worthless for this value of $\epsRS$. Decreasing to $\epsRS=10^{-3}$, we need about 400 points for the trajectory to match special quadrature. Considering that the cost of the quadrature is constant with respect to $\epsRS$, this tell us that, for sufficiently slender fibers, special quadrature will be advantageous, and going forward we will only compare special quadrature to direct quadrature and local drag.

\subsubsection{Special vs.\ direct quadrature}
\begin{figure}
\centering
\includegraphics[width=\textwidth]{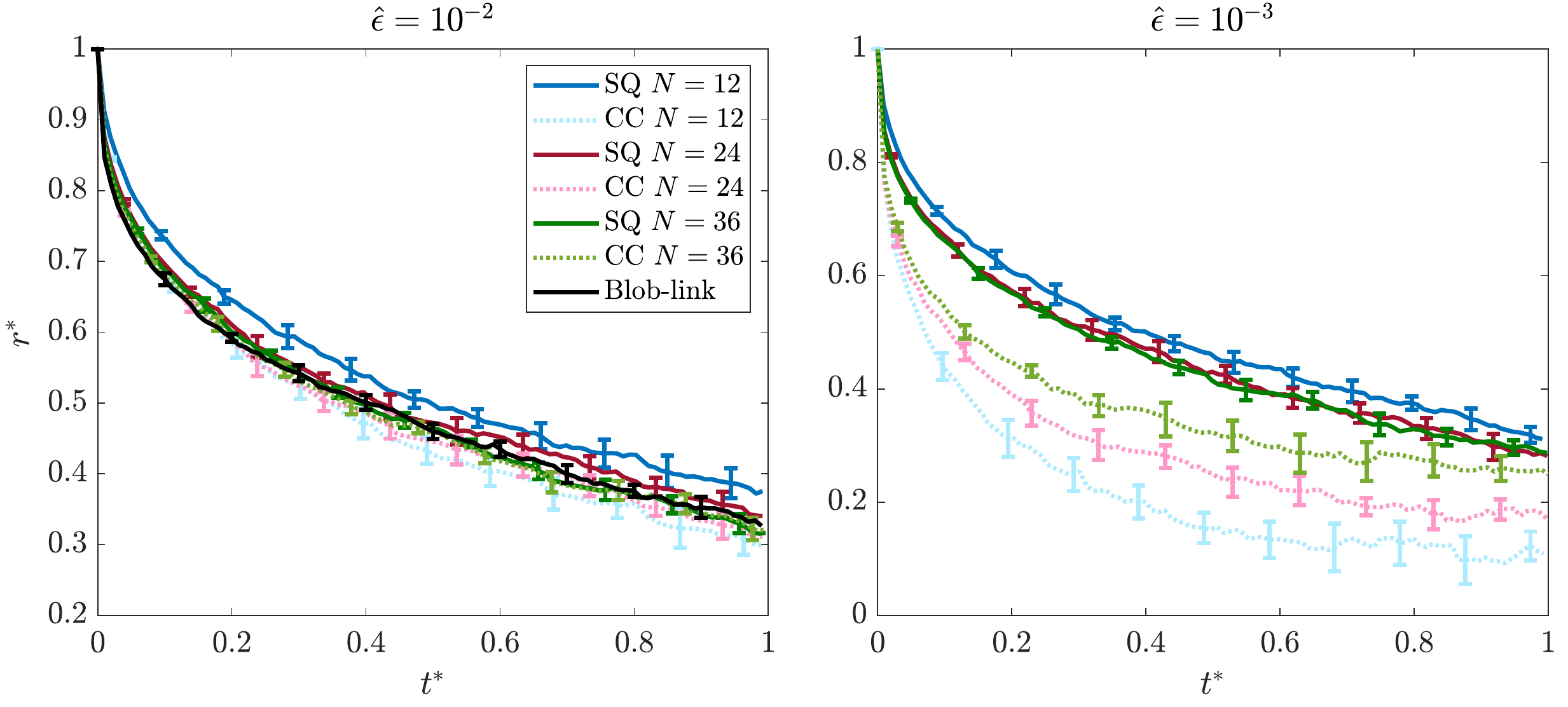}
\caption{\label{fig:DirSp} End-to-end distance trajectory comparing special quadrature (SQ) with direct quadrature (CC) with various numbers of collocation points. The left plot shows $\epsRS=10^{-2}$, for which we can use a small number of collocation points. The more slender $\epsRS=10^{-3}$ at right shows the weaknesses of direct quadrature (similar to an under-resolved bead link model).} 
\end{figure}

Let's now move on to comparing special quadrature to direct quadrature on the spectral grid. With direct quadrature, the spectral mobility becomes analogous to a blob-link mobility, and so we would expect fewer collocation points to be acceptable for large $\epsRS$. Figure\ \ref{fig:DirSp} confirms this result; there we see that when $\epsRS=10^{-2}$, the results for direct quadrature and special quadrature for $N \geq 24$ are all essentially the same. When $N=12$, we see that the direct quadrature is actually more accurate than special quadrature (relative to the blob-link or self-refined solutions), though this may be accidental (for this specific example).

As in the previous section, the special quadrature scheme shines only when we drop to $\epsRS=10^{-3}$. This time (right panel of Fig.\ \ref{fig:DirSp}), we see that special quadrature converges in space (to within statistical error) for $N \geq 24$, while direct quadrature fails to give an accurate trajectory even for $N=36$. This is not surprising; for $\epsRS=10^{-3}$ we would expect at least several hundred points (not computationally feasible due to the small required time step size) would be needed to accurately resolve the hydrodynamics via direct quadrature. Still, it it reassuring that those trajectories move toward the converged special quadrature ones.

\subsubsection{Local drag}
So far, we have established that special quadrature is only effective in the slender limit, which is also the regime where we expect local drag to give a correct answer. So what about the difference between special quadrature and local drag? This is what we study in this section. We will compare both local drag and special quadrature to special quadrature with $N=36$, which matches the blob-link calculation that is available to us, and gives a reference answer for smaller $\epsRS$.

\begin{figure}
\centering
\includegraphics[width=\textwidth]{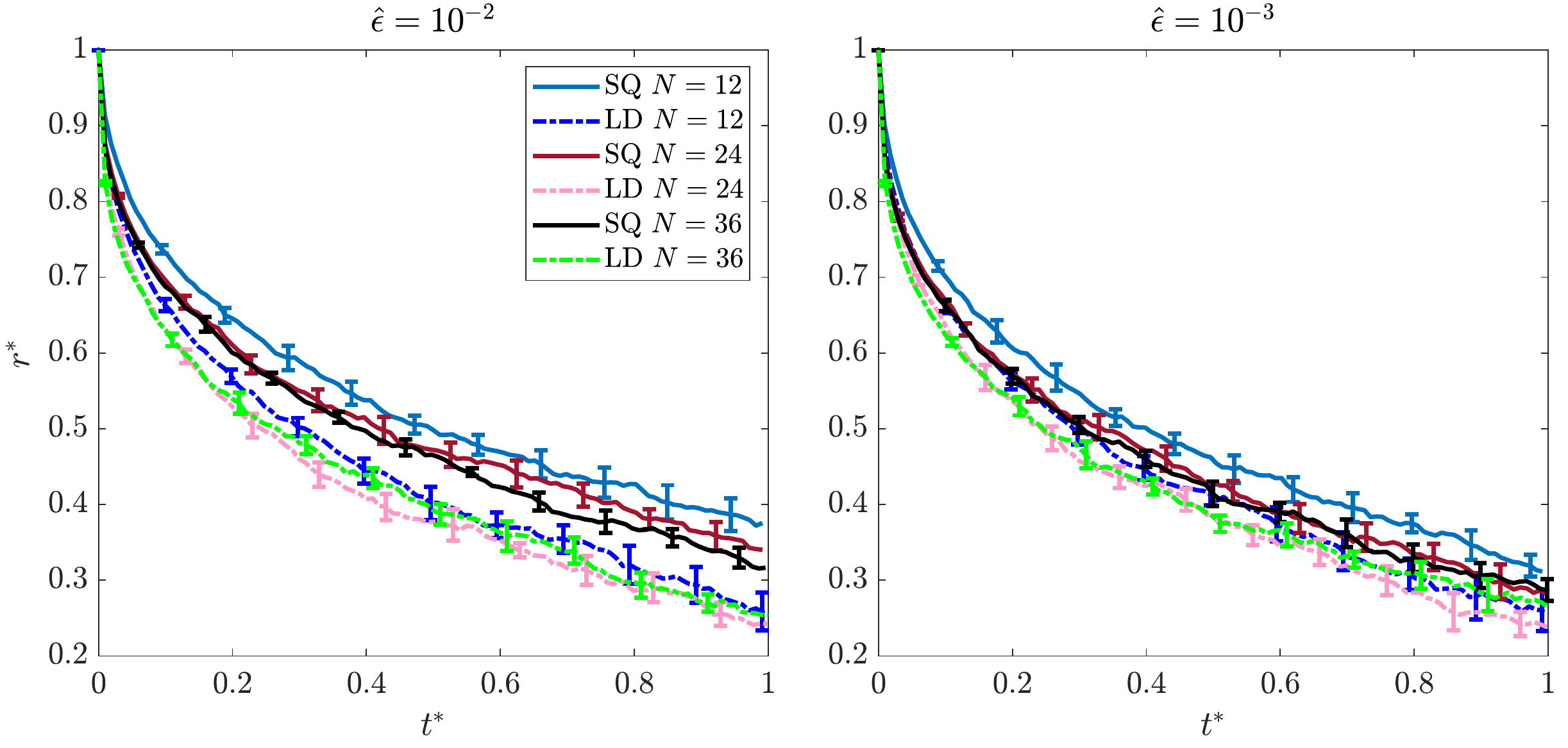}
\caption{\label{fig:LDSp} End-to-end distance trajectory comparing special quadrature (SQ) with the local drag theory\ \eqref{eq:UinnerA} with various numbers of collocation points. We use special quadrature with $N=36$, which for $\epsRS=10^{-2}$ gives the same trajectory as the blob-link code, as a reference result in both cases. The left plot shows $\epsRS=10^{-2}$, while the right shows $\epsRS=10^{-3}$.} 
\end{figure}

Figure\ \ref{fig:LDSp} shows how our local drag theory\ \eqref{eq:UinnerA} compares with special quadrature for $N=12$, $24$, and $36$ (using special quadrature with $N=36$ as a reference). The key takeaway is that local drag performs better as $\epsRS$ decreases, which is not surprising. When $\epsRS=10^{-2}$ and $N=12$, we see that local drag and special quadrature have roughly the same error at intermediate times, with local drag being more accurate at short times ($t^* < 0.1$) and special quadrature winning out on long times. For $N=24$ and $\epsRS=10^{-2}$, special quadrature is a clear winner, as the trajectory for local drag, which is unchanged when we increase to $N=36$, gives faster relaxation than the reference solution. Decreasing to $\epsRS=10^{-3}$, we see that local drag performs better: with $N=12$ it gives an error identical to special quadrature with $N=24$. This is likely just a coincidence, since with $N=24$ and $N=36$ the accuracy of local drag is clearly inferior to that of special quadrature. Still, it certainly provides a good approximation for $\epsRS=10^{-3}$. 

\begin{figure}
\centering
\includegraphics[width=\textwidth]{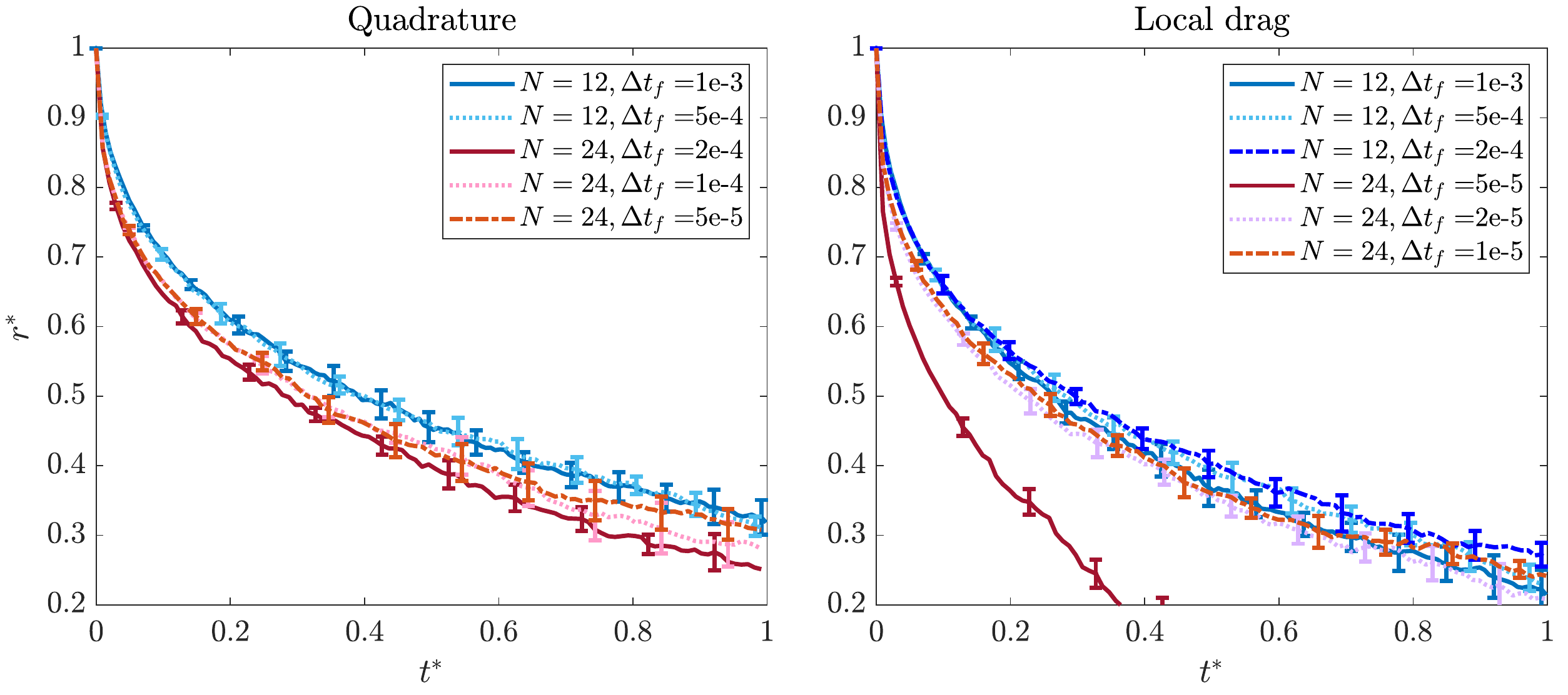}
\caption{\label{fig:LDSpTemp}Studying temporal convergence of special quadrature (left) and local drag (right). Blue colors show $N=12$, while red colors show $N=24$. This is for $\epsRS=10^{-3}$, where local drag is most effective. Time step sizes are reported in terms of $\D t_f =\D t/\tau_\text{fund}$, where $\tau_\text{fund}$ is defined in\ \eqref{eq:TauSm}.} 
\end{figure}

From a spatial perspective then, we see that local drag can be just as effective as special quadrature when $\epsRS$ is small. We still have not however, considered the temporal convergence of local drag, which we show in Fig.\ \ref{fig:LDSpTemp} for $\epsRS=10^{-3}$. In the left panel, we show how well behaved the special quadrature is: the data are basically converged when $\D t_f = 0.001$ for $N=12$ (and we could probably use an even larger $\D t$) and $\D t_f= 10^{-4}$ for $N=24$. The time step sizes for local drag are much smaller; the right panel shows that we need a time step size of $\D t_f = 5 \times10^{-4}$ for $N=12$ and $\D t_f = 2 \times 10^{-5}$ for $N=24$. Thus, while the special quadrature might be more expensive in space, it is ultimately the most efficient, since local drag requires smaller time step sizes to achieve convergence.


\section{Stochastic drift terms \label{sec:Drift}}
\setcounter{equation}{0}
\setcounter{figure}{0}    
\renewcommand{\thefigure}{C\arabic{figure}}
This appendix is devoted to discussing the stochastic drift term\ \eqref{eq:DriftGen} in the overdamped Ito Langevin equation for the tangent vector evolution. Our first goal is to show that this drift term is vital to correctly sample from the equilibrium distribution\ \eqref{eq:GBDist}. Because our midpoint temporal integrator (Section\ \ref{sec:OneSolveMP}) captures the drift term in expectation by moving the fibers to the midpoint, it is not so simple to turn off the drift term and see what happens. Because of this, in Appendix\ \ref{sec:NaiveRFD} we present a simple numerical scheme which captures the drift term via an expensive random finite difference (RFD) \cite{BrownianMultiBlobs}. In that scheme, we can easily turn off the drift term and show the disastrous effects that result when it is not included (see Appendix\ \ref{sec:PoorManResults}). After this, in Appendix\ \ref{sec:DriftOS} we turn to the more technical details of how the midpoint scheme captures the drift term\ \eqref{eq:DriftGen} in expectation. We note that the midpoint scheme as presented in the main text is optimized for dense linear algebra, as the term\ \eqref{eq:UMD} minimizes the number of evaluations of $\Mfor$ assuming a dense inversion is relatively cheap. Here we present an alternative method to generate the same drift term which applies when it is impossible to form $\Mfor$ as a dense matrix. We then compare our temporal integration scheme to the related scheme of Westwood et.\ al for rigid bodies \cite{westwood2021generalised} in Appendix\ \ref{sec:Brennan}. 

\subsection{RFD scheme \label{sec:NaiveRFD}}
We first present an RFD scheme which explicitly tells us what drift terms are being included, so that we can exclude them to demonstrate that\ \eqref{eq:ItoTau} \emph{without} the correct drift terms does not give the correct equilibrium statistical mechanics. This scheme begins by solving the system of equations
\begin{gather}
\label{eq:sysSimp}
\tdisc{\Msym}{n}\left(\M{L}\left(c\tdisc{\V{X}}{n+1,*} +(1-c)\tdisc{\V{X}}{n}\right) +\V{\Lambda}\right)+\tdisc{\V{U}}{n}_B = \tdisc{\M{K}}{n}\V \alpha \\ \nonumber
\left(\tdisc{\M{K}}{n}\right)^T  \V{\Lambda} = \V 0, \\ 
\label{eq:Xnp1st}
\text{where} \qquad \tdisc{\V{X}}{n+1,*} = \tdisc{\V{X}}{n} +\D t \tdisc{\M{K}}{n}\V{\alpha}
\end{gather}
is the linearized position at the next time step that we use in the implicit method determined by the coefficient $c$. Substituting\ \eqref{eq:Xnp1st} into\ \eqref{eq:sysSimp}, we begin by solving the saddle point system
\begin{gather}
\label{eq:SPNaive}
\tdisc{\begin{pmatrix} -\Msym & \left(\M{I}+c \D t \Msym \M{L} \right)\M{K}\\ \M{K}^T & \V 0 \end{pmatrix}}{n}
\tdisc{\begin{pmatrix} \V{\Lambda} \\ \V{\alpha} \end{pmatrix}}{n} =
\begin{pmatrix} -\Msym \M{L}\tdisc{\V{X}}{n}+\tdisc{\V{U}}{n}_B \\  \V{0}\end{pmatrix},
\end{gather}
where all time-dependent quantities are evaluated at time $n$, the Brownian velocity $\V{U}_B$ is defined in\ \eqref{eq:UB}, and $\V{\alpha}=\left(\V{\Omega},\Ump\right)$. Using the Schur complement to eliminate $\V{\Lambda}$, we obtain the solution
\begin{align}
\nonumber
\tdisc{\V{\alpha}}{n} &=\tdisc{\left(-\M{N}\M{K}^T \M{L}\V{X} + \sqrt{\frac{2 k_B T}{\D t}}\M{N}^{1/2}\gauss\right)}{n}+\mathcal{O}(\D t)\\ \label{eq:OmNaive}
& =\tdisc{\left( -\M{N}\Cbar^T \M{L}_{\Xs} \Xsbar +  \sqrt{\frac{2 k_B T}{\D t}}\M{N}^{1/2}\gauss\right)}{n}+\mathcal{O}(\D t)
\end{align}
where the terms of order $\D t$ come from using an implicit method. Comparing\ \eqref{eq:OmNaive} with the SDE\ \eqref{eq:ItoTau}, we see that we still need to capture the drift term\ \eqref{eq:DriftGen}. To do this, we add a random finite difference term \cite[Sec.~IIIC]{BrownianMultiBlobs} by setting
\begin{gather}
\label{eq:TauRFD}
\tdisc{\Xsbar}{n, \text{RFD}} = \begin{pmatrix}\text{rotate}\left(\tdisc{\Xs}{n}, \delta \tdisc{\V{\gauss}}{n, \text{RFD}}_\Bind{1:3N}\right) \\ \Xmp+\delta L \tdisc{\V{\gauss}}{n, \text{RFD}}_\Bind{3N+1:3N+3}\end{pmatrix} \\ 
\label{eq:OmegaRFD}
\tdisc{\V \alpha}{n,\text{RFD}} =\frac{k_B T}{\delta} \left(\M{N}\left(\tdisc{\Xsbar}{n}_\text{RFD} \right)-\M{N}\left(\tdisc{\Xsbar}{n}\right)\right)\begin{pmatrix}\tdisc{\V{\gauss}}{n, \text{RFD}}_\Bind{1:3N} \\ L^{-1}\tdisc{\V{\gauss}}{n, \text{RFD}}_\Bind{3N+1:3N+3} \end{pmatrix}
\end{gather}
where $\delta \ll 1$ is a dimensionless small parameter and $\tdisc{\V{\gauss}}{n, \text{RFD}}$ is generated independent of $\tdisc{\V{\gauss}}{n}$. The procedure\ \eqref{eq:TauRFD} corresponds to first rotating the tangent vectors around a random axis by a random small angle, which we follow in\ \eqref{eq:OmegaRFD} by solving two saddle point systems with the same random numbers on the right hand side. The Taylor expansion of the average $\V \alpha^\text{(RFD)}$ is (dropping the time superscript)
\begin{align*}
\Exp{\tdisc{\V \alpha}{\text{RFD}}}
& = \left(k_B T\right) \partial_j \left(\V{N}_{ik} \right)\Exp{\gauss^\text{(RFD)}_k \left(\Xsbar^\text{(RFD)}-\Xsbar\right)_j}+\mathcal{O}\left(\delta\right) \\ 
& = \left(k_B T\right) \partial_j \left(\V{N}_{ik} \right)\Exp{\gauss^\text{(RFD)}_k \gauss^\text{(RFD)}_m}  \Cbar_{j m}+\mathcal{O}\left(\delta\right)\\ 
&= \left(k_B T\right)\partial_j \left(\V{N}_{ik} \right)\Cbar_{j k}.
\end{align*}
which is exactly the drift term in\ \eqref{eq:DriftGen}. In the second equality, we use the fact that the update\ \eqref{eq:TauRFD} is $\Cbar \gauss$ to first order in $\delta$. Thus the RFD term $\tdisc{\V \alpha}{\text{RFD}}$ captures\ \eqref{eq:DriftGen} in expectation, and the fiber update 
\begin{gather}
\label{eq:RFDUpdate}
\tdisc{\Xs}{n+1}=\text{rotate}\left(\tdisc{\Xs}{n}, \D t\left( \tdisc{\V{\Omega}}{n}+\tdisc{\V \Omega}{n, \text{RFD}}\right), \right)\\ \nonumber \qquad \tdisc{\Xmp}{n+1}=\tdisc{\Xmp}{n}+\D t \left(\tdisc{\Ump}{n}+\tdisc{\Ump}{n, \text{RFD}}\right)
\end{gather}
is a weakly first-order accurate temporal integrator for\ \eqref{eq:ItoTau}.

\subsection{Drift terms in midpoint scheme \label{sec:DriftOS}}
Proceeding to the midpoint scheme in Section\ \ref{sec:OneSolveMP}, we now show that the tangent vector update\ \eqref{eq:OmMidP} generates the drift term\ \eqref{eq:DriftGen} in expectation. To simplify the notation, in this appendix we will assume that all time/position-dependent matrices and vectors are evaluated at time $n$, unless otherwise indicated. 
Substituting for $\V{U}_B$ and $\V{U}_\text{MD}$, we expand\ \eqref{eq:OmMidP} to read 
\begin{gather}
\label{eq:OneSolveAlphExp}
\tdisc{\V{\alpha}}{n+1/2,*}= -\tdisc{\left(\M{N}\M{K}^T\right)}{n+1/2,* }\M{L}\V{X}\\ \nonumber
+ \sqrt{ \frac{2k_B T}{\D t} }\tdisc{\left(\M{N}\M{K}^T \Msym^{-1}\right)}{n+1/2,*}\Msym^{1/2}\gauss\\ \nonumber
+\sqrt{ \frac{2k_B T}{\D t} }\tdisc{\left(\M{N}\M{K}^T \Msym^{-1}\right)}{n+1/2,*}\left(\tdisc{\Msym}{n+1/2,*}-\Msym \right) \left(\Msym^{-1/2}\right)^T \gauss
\end{gather}
We will Taylor expand the second and third line in succession. Beginning with the second line, we have
\begin{gather*}
\sqrt{ \frac{2k_B T}{\D t} }\left[\tdisc{\left(\M{N}\M{K}^T \Msym^{-1}\right)}{n+1/2,*}\Msym^{1/2} \gauss\right]_i=\sqrt{\frac{2 k_B T}{\D t}}\left(\M{N}\M{K}^T\Msym^{-1/2}\gauss \right)_i\\ 
\label{eq:DriftSec}
+\sqrt{\frac{2 k_B T}{\D t}}\partial_j \left(\M{N}\M{K}^T \Msym^{-1}\right)_{ih}\left(\tdisc{\Xsbar_j}{n+1/2,*}-\Xsbar_j\right) \Msym_{hp}^{1/2}\gauss_p+\mathcal{O}(\D t).
\end{gather*}
We recognize the first term from\ \eqref{eq:Nhalf} as $\M{N}^{1/2}\gauss$, which is the typical term associated with the Brownian noise. Substituting the expansion of $\tdisc{\Xs}{n+1/2,*}$ in\ \eqref{eq:TauTilde} into the second line, we get the average drift
\begin{align*}
\text{DriftUB}_i &= \sqrt{\frac{2 k_B T}{\D t}}\partial_j \left(\M{N}\M{K}^T \Msym^{-1}\right)_{ih}\sqrt{\frac{k_B T \D t }{2}}\Cbar_{ja} \M{K}^{-1}_{ab} \Msym_{bc}^{1/2} \Exp{\gauss_c \gauss_p} \Msym_{hp}^{1/2}\\
&= \left(k_B T\right) \partial_j \left(\M{N}\M{K}^T \Msym^{-1}\right)_{ih} \left(\Msym_{hb} \M{K}^{-T}_{ba}\Cbar^T_{aj}\right)
\end{align*}
where in the last equality we used $\Mfor = \Mfor^{1/2} \left(\Mfor^{1/2}\right)^T$. Rearranging this gives $\Mfor \left(\Mfor^{-1/2}\right)^T=\Mfor^{1/2}$, which alllows us to expand the last line in\ \eqref{eq:OneSolveAlphExp} as 
\begin{align}
\nonumber
\text{DriftUMD}_i &= \Exp{\sqrt{ \frac{2k_B T}{\D t} }\left(\M{N}\M{K}^T \Msym^{-1}\right)_{ih}\partial_j \left(\Msym_{hp} \right)\left(\tdisc{\Xsbar_j}{n+1/2,*}-\Xsbar_j\right)\Msym_{qp}^{-1/2} \gauss_{q}}+\mathcal{O}(\D t)\\
\nonumber
&=\left(k_B T\right)\left(\M{N}\M{K}^T \Msym^{-1}\right)_{ih}\partial_j \left(\Msym_{hp}\right)\Cbar_{ja} \M{K}^{-1}_{ab} \Msym_{bc}^{1/2} \Exp{\gauss_c \gauss_q} \Msym_{qp}^{-1/2} \\
\nonumber
&=\left(k_B T\right)\left(\M{N}\M{K}^T \Msym^{-1}\right)_{ih}\partial_j \left(\Msym_{hp}\right)\Cbar_{ja} \M{K}^{-1}_{ap} \\
\label{eq:DriftUMD}
& = \left(k_B T\right) \left(\M{N}\M{K}^T \Msym^{-1}\right)_{ih}\partial_j \left(\Msym_{hp}\Cbar_{ja} \M{K}^{-1}_{ap}\right) .
\end{align}
The last equality makes use of the identity
\begin{equation}
\label{eq:KinvKTDrift}
\M{N}_{ik}\M{K}^T_{kp} \partial_j \left( \Cbar_{jr} \M{K}^{-1}_{rp}\right)=0,
\end{equation}
which is proven later.

We can now add the two drift terms to obtain the total drift
\begin{align*}
\text{Drift}_i & = \text{DriftUB}_i +\text{DriftUMD}_i \\
& = k_B T \partial_j \left(\left(\M{N}\M{K}^T \Msym^{-1}\right) \left(\Msym \M{K}^{-T} \Cbar^T\right)\right)_{ij}\\
& = k_B T \partial_j \left(\M{N}\left(\Cbar\M{K}^{-1}\M{K}\right)^T\right)_{ij}
= k_B T \partial_j \left(\M{N}\Cbar^T\right)_{ij} = k_B T \partial_j \left(\M{N}_{ik}\right) \Cbar^T_{kj}.
\end{align*}
In the last line above, we used\ \eqref{eq:KinvK} in the second equality, thereby establishing that the one-solve scheme produces the drift\ \eqref{eq:DriftGen}, as desired.

\subsubsection{Proof of\ \eqref{eq:KinvKTDrift} \label{sec:KinvKTDrift}}
We now prove that\ \eqref{eq:KinvKTDrift} is indeed zero by substituting the definitions of $\M{K}$ and $\M{K}^{-1}$ to get
\begin{align*}
\M{N}_{ik}\M{K}^T_{kp} \partial_j \left( \Cbar_{jr} \M{K}^{-1}_{rp}\right)&=\M{N}_{ik}\Cbar_{ky}^T \X^T_{yp} \partial_j \left(\Cbar_{jr}\Cbar^T_{rb}\right) \X^{-1}_{bp}\\
&= \M{N}_{ik}\Cbar^T_{ky} \partial_j \left(\Cbar_{jr}\Cbar^T_{ry}\right).
\end{align*}
Now, the matrix $\left(\Cbar_{jr}\Cbar^T_{rb}\right)$ is a $3(N+1) \times 3(N+1)$ matrix composed of $N+1$ diagonal blocks. The first $N$ diagonal blocks are $\M{I}_3-\Xs_\Bind{p}\Xs_\Bind{p}$, while the last block is the $3 \times 3$ identity. It follows that the derivative 
\begin{equation}
\frac{\partial}{\partial \Xsbar_j} \left(\Cbar_{jr}\Cbar^T_{ry}\right) = -4 \Xsbar_y,
\end{equation}
for $b \leq 3N$, and zero otherwise (for the parts of $\Xsbar$ associated with constant motions). Multiplying by $\Cbar^T$ and using\ \eqref{eq:CXs} gives\ \eqref{eq:KinvKTDrift}.

\subsubsection{Alternate way of obtaining mobility drift \label{sec:MobDrAlt}}
An alternative way of obtaining the drift term\ \eqref{eq:UMD} is via an RFD, which avoids potentially expensive resistance problems. The alternative expression for $\V{U}_\text{MD}$ in this case is
\begin{gather}
\nonumber
\V \mu^\text{(RFD)} = \M{K}^{-1} \gauss^\text{(RFD)} \qquad 
\Xs^\text{(RFD)} = \text{rotate}\left(\Xs, \delta L \V \mu^\text{(RFD)} \right)\approx \Xs_p -\delta L \M{C}\M{K}^{-1}\gauss^\text{(RFD)} \\ 
\label{eq:MRFD}
\V{U}_\text{MD} =\frac{k_B T}{\delta L} \left(\Msym \left(\Xs^\text{(RFD)}\right)-\Msym\left(\Xs\right)\right) \gauss^\text{(RFD)},
\end{gather}
and it is not hard to show this generates the drift term\ \eqref{eq:DriftUMD} in expectation, for small $\delta$. Numerical tests show that using the expression\ \eqref{eq:UMD}, rather than\ \eqref{eq:MRFD}, gives a closer approximation to the equilibrium distribution for larger time step sizes, whereas for smaller time step sizes the two expressions perform similarly. In this paper, we use\ \eqref{eq:UMD}, but\ \eqref{eq:MRFD} will be useful for many-fiber systems in which we can apply $\Mfor$ fast, while inverting $\Mfor$ requires an iterative solver. 

\subsubsection{Comparison to scheme of Westwood et al.\ \cite{westwood2021generalised} \label{sec:Brennan}}
The midpoint method proposed here applies for any kinematic matrix $\M{K}\left(\V{X}\right)$ and mobility $\M{M}\left(\V{X}\right)$ and can therefore be used for suspensions of rigid bodies as well. In fact, our scheme is related to the generalized drift-correcting (gDC) scheme proposed by Westwood et al.\ \cite[Sec.~4.1]{westwood2021generalised} for rigid body suspensions.

For a single rigid body modeled as a rigid-multiblob with blob positions $\V{X}_i$ (generalization to many rigid bodies is straightforward) the authors of \cite{westwood2021generalised} define the kinematic matrix $\M{K}\left(\V{X}\right)$ through its action on a $6 \times 1$ rigid body velocity $\V{U} = \left[\V{u}, \V{\omega} \right]$, so that
\begin{equation}
\partial_{t} \V{X}_i = \V{u} + \V{\omega} \times \left(\V{X}_i - \V{q} \right) = \left[ \M{K}\left(\V{X}\right) \V{U} \right]_i,
\end{equation} 
where $\V{q}$ is an arbitrary tracking point for the rigid structure. Equation (E.5) in \cite{westwood2021generalised} gives the drift term for rigid bodies as  
\begin{equation}
\text{Drift}_i = k_B T \partial_j \M{N}_{ij}.
\end{equation} 
Similar to the midpoint scheme presented here, the gDC scheme captures the thermal drift term by splitting it into two parts
\begin{align}
\label{eq:gDCdrift}
\text{Drift}_i &= \left(k_B T\right) \partial_j \M{N}_{ij}  \\
\nonumber
             &= \left(k_B T\right) \partial_j \left( \M{N} \M K^{T} \M{M}^{-1/2}  \M{M}^{1/2} \M K^{-T} \right)_{ij} \\ 
\nonumber
             &= \left(k_B T\right) \partial_j \left( \M{N} \M K^{T} \M{M}^{-1/2}\right)_{ik} \left( \M{M}^{1/2} \M K^{-T} \right)_{kj} + k_B T \left( \M{N} \M K^{T} \M{M}^{-1/2}\right)_{ik} \partial_j \left( \M{M}^{1/2} \M K^{-T} \right)_{kj}, 
\end{align}
where $\M K^{-T} = \M K \left(\M K^{T} \M K \right)^{-1}$ is the left pseudoinverse of $\M K^{T}$ and can be computed very efficiently. By contrast, the midpoint scheme presented in this work splits the drift term according to 
\begin{align}
\nonumber
\text{Drift}_i &= \left(k_B T\right) \partial_j \M{N}_{ij}  \\
\nonumber
             &= \left(k_B T\right) \partial_j \left( \M{N} \M K^{T} \M{M}^{-1} \M{M}^{1/2} \M{M}^{1/2} \M K^{-T} \right)_{ij} \\ 
\label{eq:MIDdrift}
             &= \left(k_B T\right) \partial_j \left( \M{N} \M K^{T} \M{M}^{-1}\right)_{ik} \left( \M{M} \M K^{-T} \right)_{kj} + k_B T \left( \M{N} \M K^{T} \M{M}^{-1}\right)_{ik} \partial_j \left( \M{M} \M K^{-T} \right)_{kj}. 
\end{align}
While equations \eqref{eq:gDCdrift} and \eqref{eq:MIDdrift} appear very similar, the approach used in our midpoint scheme only requires one application of $\M{M}^{1/2}$, while the gDC scheme requires four. In some cases, saving on applications of $\M{M}^{1/2}$ can significantly reduce the computational cost of the whole scheme.   

To elaborate on why the gDC scheme requires extra applications of $\M{M}^{1/2}$, we will briefly summarize how the scheme captures each term in \eqref{eq:gDCdrift}. The second term is calculated efficiently by using an RFD for 
\begin{equation} \label{eq:gDCrfd}
\V{\nu}_k = \partial_j \left( \M{M}^{1/2} \M K^{-T} \right)_{kj}
\end{equation}
since it does not involve $\M N$, while the first term is captured using a midpoint time integration scheme of the form
\begin{align}
\tdisc{\V{X}}{n+1/2,*} - \tdisc{\V{X}}{n} &= \sqrt{ \frac{k_B T \D t}{2} } \M K^{-1} \M{M}^{1/2} \tdisc{\V{\eta}}{n} \\
\tdisc{\widetilde{\V{V}}}{n+1/2,*} &= \sqrt{ \frac{2 k_B T}{ \D t} } \tdisc{\left( \M{N} \M K^{T} \M{M}^{-1/2}\right)}{n+1/2,*} \tdisc{\V{\eta}}{n}, \label{eq:gDCmid}
\end{align}
so that 
\begin{align*}
\Exp{\tdisc{\widetilde{\V{V}}}{n+1/2,*}} &= \left(k_B T\right) \partial_j \left( \M{N} \M K^{T} \M{M}^{-1/2}\right)_{ik} \Exp{\tdisc{\V{\eta}}{n}_k \tdisc{\V{\eta}}{n}_r} \left( \M K^{-1} \M{M}^{1/2} \right)_{jr} \\
&= \left(k_B T\right) \partial_j \left( \M{N} \M K^{T} \M{M}^{-1/2}\right)_{ik} \left( \M{M}^{1/2} \M K^{-T} \right)_{kj}.
\end{align*}
The gDC scheme combines these parts to compute 
\begin{equation}
\text{Drift}_i = k_B T \left( \M{N} \M K^{T} \M{M}^{-1/2}\right)_{ik} \V{\nu}_k + \Exp{\widetilde{\V{V}}^{n+1/2,*}} = \left(k_B T\right) \partial_j \M{N}_{ij}.
\end{equation} 
Note that calculating $\tdisc{\widetilde{\V{V}}}{n+1/2,*}$ according to\ \eqref{eq:gDCmid} requires two applications of $\M{M}^{1/2}$ (one at each time level), and calculating $\V{\nu}_k$ using an RFD requires two additional applications of $\M{M}^{1/2}$. Hence the gDC scheme requires four total applications of $\M{M}^{1/2}$ while the midpoint scheme presented in this work only requires one.

\section{Small fluctuations of a curved fiber\label{sec:CurvX0}}
\setcounter{equation}{0}
\setcounter{figure}{0}    
\renewcommand{\thefigure}{D\arabic{figure}}
The goal of this appendix, which is meant to accompany the corresponding calculations for a free filament in Section\ \ref{sec:EqStat}, is to study the performance of our spatial and temporal discretization in the context of small fluctuations around an equilibrium state. The motivation for studying this problem is both physical and numerical: physically, it is useful to study because we are eventually interested in fibers that are constrained by cross linkers, which roughly hold them in place while they fluctuate. Numerically, starting with small fluctuations \cite{kantsler2012fluctuations}, and a curved, rather than straight, fiber allows us to confirm that our method correctly handles nonlinearities. We discuss the curved configuration we choose, and how we keep the dynamics from drifting away from it, in Section\ \ref{sec:smallsetup}.

Because the dynamics are linearized around a particular state, we can compute a theoretical covariance matrix and break the dynamics into a set of modes (from the linearized covariance). In Section\ \ref{sec:MCMC2nd}, we compute the equilibrium variance of each of these modes via Markov Chain Monte Carlo (MCMC) calculations using the equilibrium probability distribution\ \eqref{eq:GBDist}. By comparing our spectral discretization to a blob-link one, we verify that our equilibrium distribution\ \eqref{eq:GBDist} gives behavior similar to that of the more physical one with uniformly-spaced links. Following this, we take advantage of the small fluctuations to linearize the SDE\ \eqref{eq:ItoX} around the base state $\V{X}_0$, which we then diagonalize in Section\ \ref{sec:LinTScales} to compute relaxation timescales of the fundamental modes. This aids us in Section\ \ref{sec:LinTInt}, where we study the time step size required for our midpoint temporal integrator to accurately sample from the equilibrium probability distribution\ \eqref{eq:GBDist}.

\subsection{Set-up for small fluctuations \label{sec:smallsetup}}
\begin{figure}
\centering
\includegraphics[width=\textwidth]{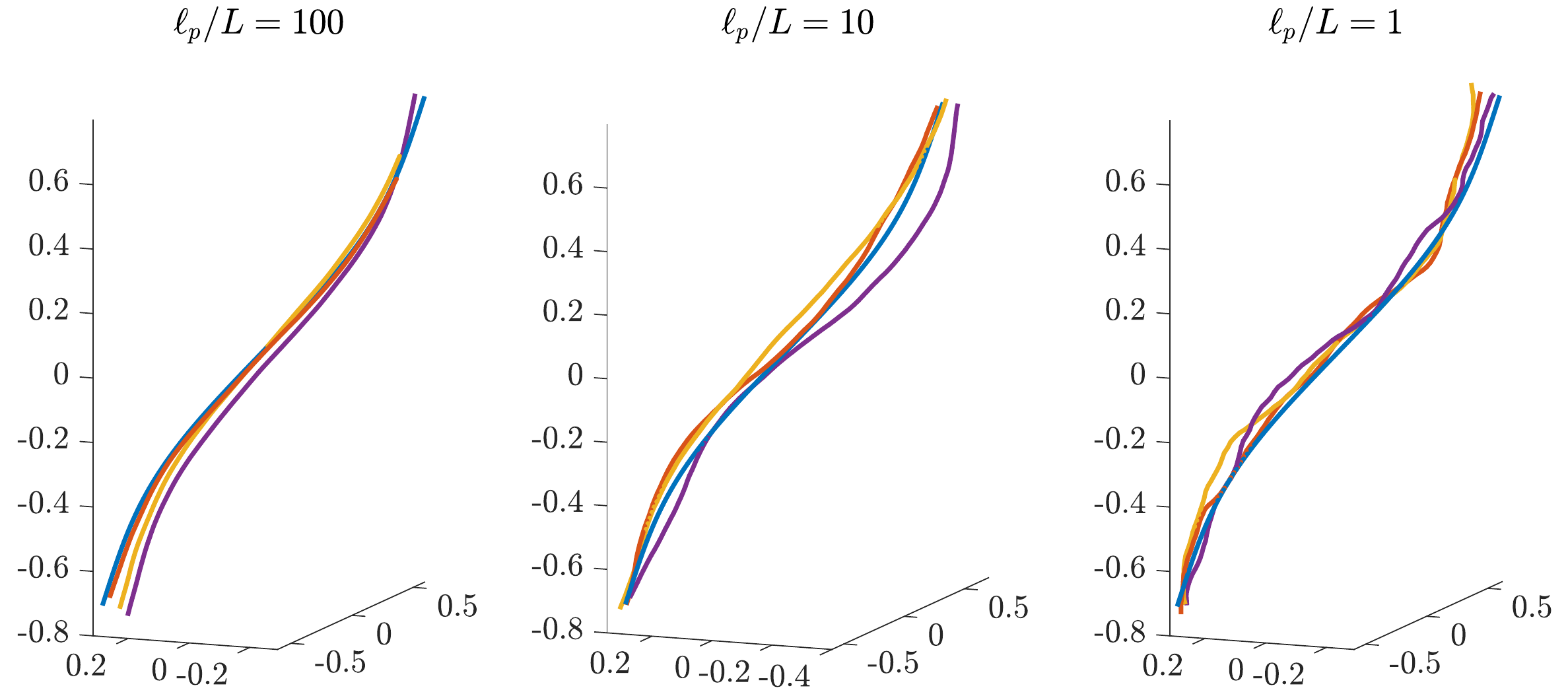}
\caption{\label{fig:DiscChains}Three samples of penalty-bound blob-link chains (with 100 links) for various $\ell_p/L$, sampled using MCMC.  We show the ``random'' filament shapes $\V{X}=\V{X}_0+\D \V{X}$, where $\D \V{X}$ is sampled from the equilibrium distribution\ \eqref{eq:EqdX}. The blue line shows $\V{X}_0$.}
\end{figure}

For the base state $\V{X}_0$, we introduce a curved configuration with $L=2$ $\mu$m and 
\begin{equation}
\label{eq:X0C}
\Xs_0(s) = \frac{1}{\sqrt{2}}\begin{pmatrix} \cos{\left(s^3(s-L)^3\right)}, & \sin{\left(s^3 (s-L)^3\right)}, & 1 \end{pmatrix}
\end{equation}
with $\V{X}_0$ defined as the integral of this on the $N+1$ point Chebyshev grid with $\Xmp=\V 0$. To keep the dynamics near this base state, we introduce an additional energy which penalizes the discrete squared $L^2$ norm of $\D \V{X}=\V{X}-\V{X}_0$, 
\begin{gather}
\label{eq:EP}
\mathcal{E}_P = \frac{P}{2} \D \V{X}^T \Wt \D\V{X}, \qquad P = (1.6 \times 10^4)\frac{k_B T}{L^3}. 
\end{gather}
This choice of $P$ ensures that the dimensionless discrete $L^2$ norm $\D \V X^T \Wt \D \V X/L^3$ remains constant when $k_B T$ and $L$ change. While this is the required scaling of $P$ to keep the relative magnitude of $\D \V{X}$ roughly constant (and small), throughout this section we will fix $k_B T=4.1 \times 10^{-3}$ pN$\cdot \mu$m and $L=2$ $\mu$m, thus setting $P=8.2$ pN/$\mu$m$^2$. To modify the filament's persistence length, we will modify the bending stiffness $\kappa$ in the modified bending energy that keeps the curved $\V{X}_0$ as an equilibrium configuration
\begin{gather}
\label{eq:EBendP}
\mathcal{E}_\text{bend} =\frac{1}{2} \kappa \D \V{X}^T  \left(\M{D}^2\right)^T \Wt \M{D}^2 \D\V{X}.
\end{gather}
In this case, the Ito Langevin equation\ \eqref{eq:ItoX} can be viewed as an equation for $\D \V{X}$ rather than $\V{X}$, and the matrix $\M{L}$ becomes 
\begin{gather}
\label{eq:LMatNp1}
\M{L}=P \Wt +  \kappa \left(\M{D}^2\right)^T \Wt \M{D}^2,
\end{gather}
which is the original definition\ \eqref{eq:Lmat} modified to account for the penalty force. Notice that the first term in the energy (the penalty term) is independent of $\kappa$, while the bending term is proportional to $\kappa$. Figure\ \ref{fig:DiscChains} shows $\V{X}_0$ along with some samples from the Gibbs-Bolztmann distribution\ \eqref{eq:GBDist} with this energy.

We will use the covariance in $\D \V{X}$ to determine the accuracy of our spatial and temporal discretizations. Because the energy is nonlinearly constrained, we cannot write down the exact covariance, but it is informative to project the modes onto the covariance we would obtain if we replace the nonlinear inextensibility constraint with the \emph{linearized} version $\D \V{X} = \M{K}\left[\V{X}_0\right] \V{\alpha}:=\M{K}_0 \V{\alpha}$. The resulting energy $\frac{1}{2}\D \V{X}^T \M{L}\D \V{X}=\frac{1}{2} \V{\alpha}^T \M{K}_0^T \M{L}\M{K}_0 \V{\alpha}$ dictates that $\V{\alpha}$ has covariance 
\begin{equation}
\M{C}_{\V{\alpha}} = \Exp{\V{\alpha} \V{\alpha}^T}=\left(k_B T\right)\left(\M{K}_0^T \M{L}\M{K}_0 \right)^\dagger.
\end{equation}
Now pre- and post-multiplying by $\M{K}_0$, we get the expected covariance for $\D \V{X}$ as
\begin{equation}
\label{eq:CovExp}
\M{C} =\Exp{\D \V{X} \D \V{X}^T}=\left(k_B T\right) \M{K}_0\left(\M{K}_0^T \M{L} \M{K}_0 \right)^\dagger\M{K}_0^T,
\end{equation}
which is valid in the limit of small fluctuations. 

To check how close the true covariance is to $\M{C}$, we will generate samples of $\D \V{X}$ through both MCMC sampling and Langevin dynamics, then use the procedure outlined in Section\ \ref{sec:CovExp} to project the resulting covariance onto the eigenmodes of $\M{C}$. This amounts to computing the Rayleigh quotient of the eigenmodes using the true covariance and normalizing with respect to the eigenvalues of $\M{C}$. We focus on the variance in each of the modes (diagonal entries of the covariance matrix projected onto the eigenmodes), as the off-diagonal entries for the smallest persistence length we consider ($\ell_p/L=1$) are zero within the statistical error of our MCMC calculations.

\subsubsection{Covariance calculations \label{sec:CovExp}}
To compare the covariance from MCMC to\ \eqref{eq:CovExp}, we first use a blob-link discretization to compute the matrix $\M{C}$ exactly as defined in\ \eqref{eq:CovExp}. In order for the eigenmodes of $\M{C}$ to be the same across all discretizations, we need to ensure that they are orthonormal in $L^2$ and not $\mathbb{R}^n$. This can be accomplished by computing the eigenvalues of $\Wt^{1/2} \M{C} \Wt^{1/2}$, which is a Hermitian matrix and therefore has eigendecomposition
\begin{equation}
\label{eq:CtildEig}
\Wt^{1/2} \M{C} \Wt^{1/2} = \widetilde{\M{V}}\M{\Lambda}\widetilde{\M{V}}^T \qquad \widetilde{\M{V}}^T \widetilde{\M{V}} = \M{I}
\end{equation}
The matrix $\Wt$ is used to compute the $L^2$ inner product, $\norm{\V{X}}_{L^2}^2 \approx \V{X}^T \Wt \V{X}$. In the case of the blob-link method, $\Wt = \D s \, \text{diag}\left(1/2, 1, \dots, 1, 1/2\right)$. Now, if we use\ \eqref{eq:CtildEig}, we observe that 
\begin{align*}
\nonumber
\M{C}&=\Wt^{-1/2} \widetilde{\M{V}}\M{\Lambda}\left(\Wt^{-1/2} \widetilde{\M{V}} \right)^T 
\end{align*}
Thus, if we define $\M{V}=\Wt^{-1/2} \widetilde{\M{V}}$, we have (using the orthonormality of $\widetilde{\M{V}}$ in $\mathbb{R}^n$) the desirable property that 
\begin{equation}
\label{eq:D2}
\M{C} = \M{V} \M{\Lambda}\M{V}^T \qquad \M{V}^T \Wt \M{V} = \M{I},
\end{equation}
which means the columns of $\M{V}$ are orthonormal in $L^2$. Then, to compute the covariance, we can project the modes of $\D \V{X}$ onto the columns of $\M{V}$ in $L^2$ as $\M{V}^T \Wt \D \V{X}$, since 
\begin{align}
\nonumber
\Exp{\D \V{X} \D \V{X}^T} &=\M{C}\\ 
\nonumber
\label{eq:DiagL2}
\Exp{\M{V}^T \Wt \D \V{X} \left(\M{V}^T \Wt \D \V{X}\right)^T }&= \M{V}^T \Wt\M{C} \Wt \M{V} =\M{\Lambda},
\end{align}
using\ \eqref{eq:D2} in the last equality.

\subsection{MCMC estimation of reference covariance\label{sec:MCMC2nd}}
To establish a reference result for the variance of each of the modes of $\M{C}$, we use MCMC to sample from the Gibbs-Boltzmann measure\ \eqref{eq:EqdX}, but with $\D \V X=\V{X}-\V{X}_0$ taking the place of $\V{X}$ in\ \eqref{eq:EqdX} and\ \eqref{eq:pacc}. 

\subsubsection{Blob-link MCMC \label{sec:MCMCSec}}
To establish a benchmark for the spectral code, we perform MCMC on a blob-link chain with $L=2$ and $N=100$ links. To do this, we initialize the chain with shape $\V{X}=\V{X}_0$. Then, at each Monte Carlo step, we propose a new chain $\widetilde{\V{X}}$ that is generated by rotating the tangent vectors by an oriented angle 
\begin{equation}
\label{eq:MCprop}
\widehat{\V{\Omega}}^\perp = 0.1 \sqrt{k_B T}\M{L}^{-1/2} \gauss,
\end{equation}
where $\gauss$ is a $3N+3$ vector of i.i.d standard normal random variables, and the constant of 0.1 is added to make the acceptance ratio roughly 40\%. Note that this proposal is informed by the bending energy (unconstrained covariance) matrix because there are many links in the blob-link chain, and so the multiplication of the Gaussian by $\M{L}^{-1/2}$ effectively projects the proposal onto the expected (unconstrained) covariance, which reduces the number of samples required to equilibrate the higher-order modes. When $P=0$ (for a free fiber), we compute the proposal\ \eqref{eq:MCprop} using the pseudo-inverse of $\M{L}$, since constant and linear modes have zero eigenvalues. Note that in this special case it is possible to construct direct (exact) independent samplers to obtain the chain configurations, though these were not used in this work. 

\subsubsection{Results}
Figure\ \ref{fig:DiscChains} shows some sample 100-link chains with various values of $\kappa$, which we quantify by the dimensionless persistence length $\ell_p/L=\kappa/(L k_B T)$. While the relative mean-square deviation of the chains from the blue base state is approximately the same (1\%) in all cases,  we see more and more ``wiggles'' in the chain as we decrease $\ell_p/L$.

\begin{figure}
\centering
\subfigure[Blob-link]{
\includegraphics[width=0.48\textwidth]{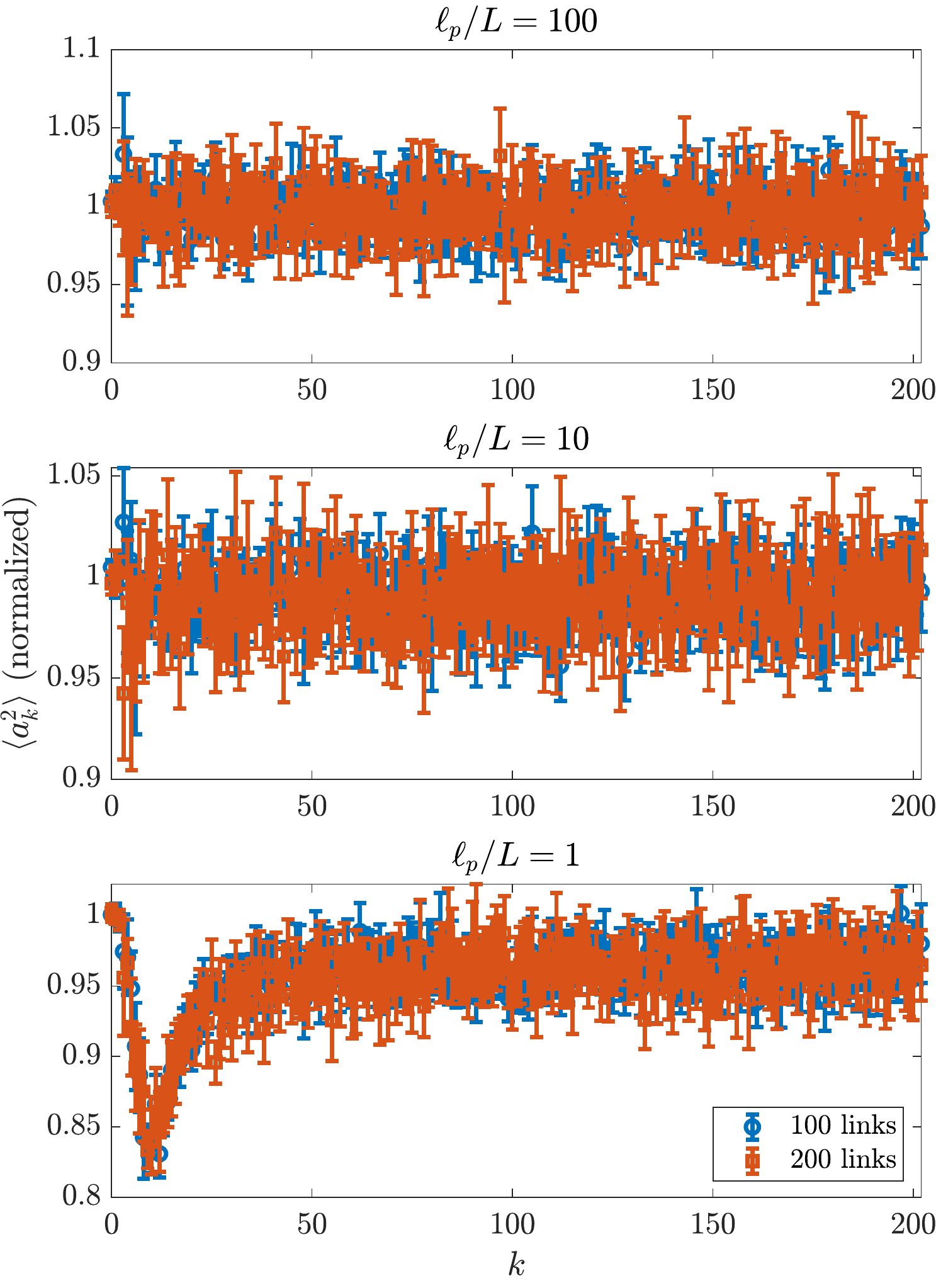}}
\subfigure[Spectral]{
\includegraphics[width=0.48\textwidth]{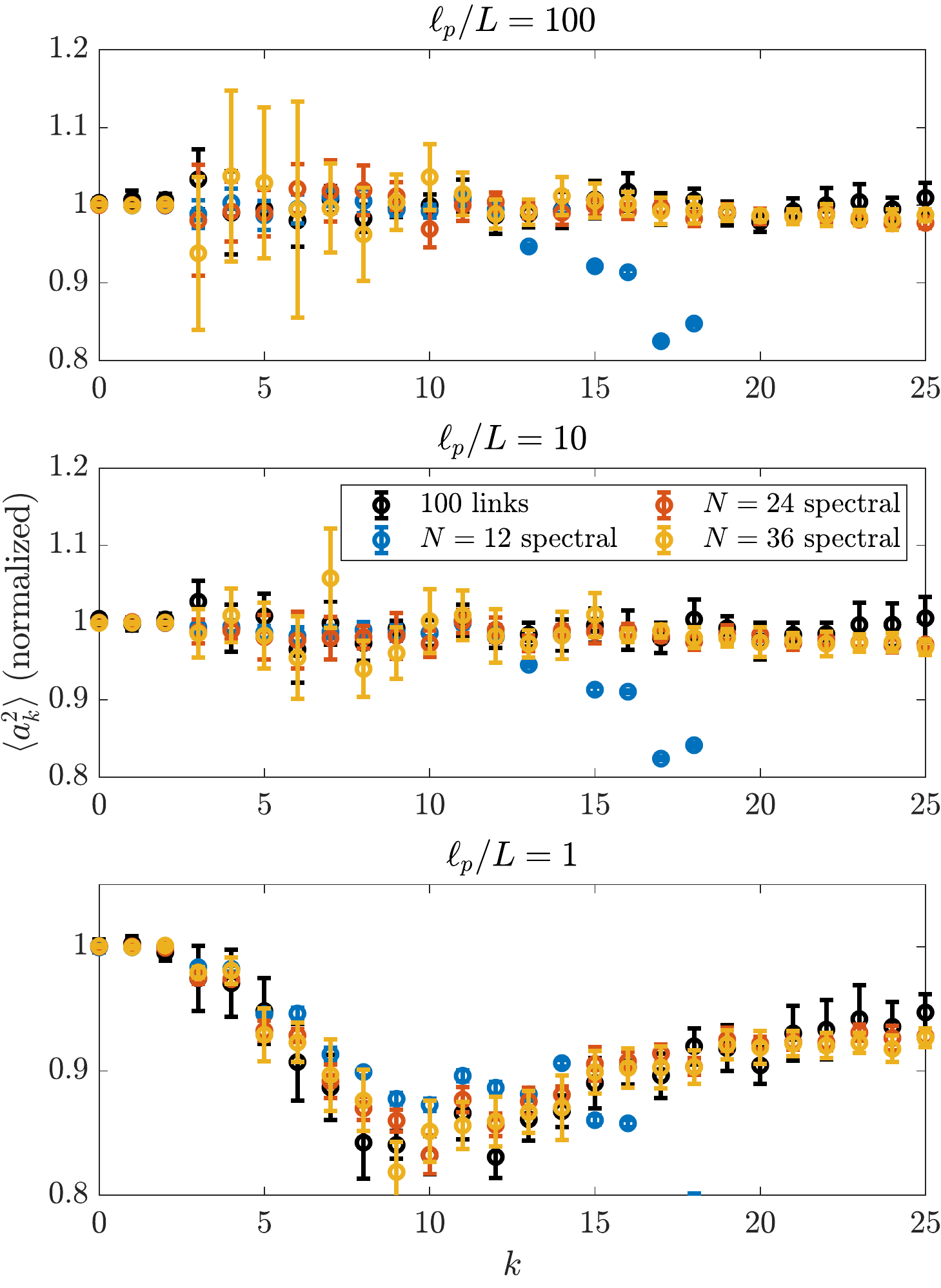}}
\caption{\label{fig:MCMCPen} Covariance of linearized modes for various $\ell_p/L$ using blob-link and spectral discretizations, estimated using MCMC sampling. We always project the variance onto the same set of $L^2$ orthonormal modes given by the eigenvectors of the linearized covariance\ \eqref{eq:CovExp}, as discussed in Section\ \ref{sec:CovExp}. (a) The covariance for a blob-link discretization with 100 links and 101 blobs (203 modes). Results are shown for 100 links (in blue), and 200 links (in red) to verify convergence. (b) Comparison of the variance of the first 25 modes using a spectral grid with those of the blob-link discretization. The black points show the results of MCMC calculations with 100 links, while the blue, red, and yellow symbols show results with $N=12, 24$, and 36 Chebyshev nodes, respectively. When $N=12$ and $k \geq 17$, the variance is reduced significantly to about 0.2 and is not shown.}
\end{figure}

To compute the variance of each mode and compare it to the expected result\ \eqref{eq:CovExp}, we use the MCMC sampler to generate $10^6$ samples of $\V{X}$ for $\ell_p/L=100,10$, and 1.  We then throw out the first 20\% of the samples and compute the covariance of the rest of the samples. We repeat this ten times to generate error bars, and show the variance in each eigenmode in Fig.\ \ref{fig:MCMCPen}(a), where we normalize by the corresponding eigenvalue and do the projection in $L^2$ as discussed in Section\ \ref{sec:CovExp}. In addition to running simulations with 100 links, we also run with 200 links to confirm that the results are converged in space. We see that, for $\ell_p/L \gtrsim 10$, the variance of each mode is exactly that predicted by the theory\ \eqref{eq:CovExp}. When we decrease $\kappa$ so that $\ell_p/L=1$, some of the intermediate modes (from about 5 to 20) are damped relative to the theoretical prediction. In this case the bending fluctuations are sufficiently large for nonlinear effects to matter.

\subsubsection{MCMC for spectral discretization \label{sec:MCMCSpec}}
Let us now perform MCMC on a chain discretized using Chebyshev collocation points, rather than uniformly-spaced nodes. To do this, we choose a number of Chebyshev collocation points $N$ and set $\V{X}=\V{X}_0$. Then, we propose a rotation of the tangent vectors at each of the Chebyshev nodes\footnote{Note that multiplication of $\gauss$ by $\M{L}^{-1/2}$ (c.f.\ \eqref{eq:MCprop}) is not necessary since in the spectral discretization the number of Chebyshev nodes is intended to be relatively small.}
\begin{equation}
\label{eq:OmPropSpec}
\widetilde{\V{\Omega}}_\Bind{p}=\frac{0.48}{N}\sqrt{\frac{L}{\ell_p}} \gauss_\Bind{p}
\end{equation}
where $\gauss_\Bind{p}$ is a vector of three i.i.d.~standard normal numbers for each node $p$, and the constant in front is chosen so that roughly 50\% of the samples are accepted if the tangent vectors are rotated by $\widetilde{\V{\Omega}}_p $ and the middle of the fiber is held in place.

The second part of the proposal is to update the fiber midpoint. To do this, we propose a new midpoint via 
\begin{equation}
\widetilde{\V{X}}_\text{MP} = \V{X}_\text{MP}+\left(7.5 \times 10^{-3}\right)L \gauss,
\end{equation}
where $\gauss$ is a vector of three i.i.d.~standard normal numbers. Now, notice that the scaling of $\Delta \V{X}_\text{MP}$ does not depend on $N$ or $\kappa$, since the penalty energy is independent of these two quantities. As such, the constant $7.5 \times 10^{-3}$ is chosen so that roughly 50\% of the samples are accepted if the tangent vectors are not updated and the fiber only translates. We combine the two proposals $\widetilde{\V{\Omega}}_p$ and $\widetilde{\V{X}}_\text{MP}$ into a single proposal 
\begin{equation}
\widetilde{\V{X}} = \X \begin{pmatrix} \text{rotate}\left(\Xs, \widetilde{\V{\Omega}}\right) \\ \widetilde{\V{X}}_\text{MP} \end{pmatrix},
\end{equation}
and compute the energy via $\frac{1}{2}\D \widetilde{\V{X}}^T \M{L}\D \widetilde{\V{X}}$; the total acceptance ratio is between 20 and 30\%. 

Repeating our MCMC procedure on the spectral discretization, we look at the variance of each of the eigenmodes of the linearized covariance\ \eqref{eq:CovExp} in Fig.\ \ref{fig:MCMCPen}(b). We show three different values of $N$, and observe that the spectral equilibrium distribution for $N=12$ has variance which diverges slightly from the blob-link variance for small mode numbers when $\ell_p/L=1$, and for larger mode numbers ($k \geq 17$) for all $\ell_p$. Increasing to $N=24$ and $N=36$, we see overlap with the variance from the blob-link discretization with 100 links for all values of $\ell_p$, which shows that the spectral method with a sufficiently large $N$ can properly reproduce the magnitude of the small equilibrium fluctuations of the blob-link chain (at least for this example). 

\subsection{Dynamics}
Our task now is to see if our temporal integrators can be used to sample from the Gibbs-Boltzmann distribution\ \eqref{eq:GBDist} via Langevin dynamics. Before doing that, we look at the fundamental timescales governing system dynamics. The largest of these is $\tau_\text{fund}$, which we use in the main body of the paper to non-dimensionalize the time step sizes.

\subsubsection{Linearized timescales \label{sec:LinTScales}}
One of the benefits of considering dynamics that are only slightly perturbed from the equilibrium configuration $\V{X}_0$ is that we can linearize the SDE\ \eqref{eq:ItoX} and diagonalize the result to compute the fundamental timescales in the system. Because the linearized mobility matrix is constant (evaluated at $\V{X}_0$), the drift term disappears and the linearized form of the SDE\ \eqref{eq:ItoX} around $\V{X}=\V{X}_0$ is
\begin{equation}
\label{eq:LangtoD}
\frac{d \D \V{X}}{d t} = -\widehat{\M{N}}_0\M{L}\D \V{X} + \sqrt{2 k_B T} \widehat{\M{N}}_0^{1/2}\Lop{W},
\end{equation}
where $\widehat{\M{N}}_0=\M{K}_0 \M{N}_0 \M{K}_0^T$, with $\M{N}$ defined in\ \eqref{eq:Ndef}. Since the null space of $\widehat{\M{N}}_0$ is the $N$ \emph{extensible} motions around $\V{X}=\V{X}_0$, there are $2N+3$ remaining directions, and the SDE\ \eqref{eq:LangtoD} diagonalizes when we consider the first $2N+3$ generalized eigenvectors of $\widehat{\M{N}}_0^\dagger$ and $\M{L}$, which satisfy
\begin{equation}
\label{eq:GenEig}
\widehat{\M{N}}_0^\dagger \M{V} = \M{L} \M{V}\M{\Lambda},
\end{equation}
where the $2N+3$ eigenvectors, which make up the $(3N+3) \times (2N+3)$ matrix $\M{V}$, are normalized such that
\begin{equation}
\M{V}^T \M{L} \M{V} = \M{I} \rightarrow \M{V}^T\widehat{\M{N}}_0^\dagger \M{V}=\M{\Lambda}.
\end{equation}
These equations imply that 
\begin{gather}
\widehat{\M{N}}_0^\dagger \M{V} =  \M{L}\M{V}\M{\Lambda}\M{V}^T \M{L} \M{V}  \rightarrow \widehat{\M{N}}_0^\dagger =\M{L} \M{V}\M{\Lambda}\M{V}^T \M{L} \rightarrow
\widehat{\M{N}}_0^{-1/2}= \M{L} \M{V}\M{\Lambda}^{1/2},
\end{gather}
which gives a square root of $\widehat{\M{N}}_0$ as $\widehat{\M{N}}_0^{1/2}=\widehat{\M{N}}_0\widehat{\M{N}}_0^{-1/2}$. Substituting this into the linearized Langevin equation\ \eqref{eq:LangtoD} and using the definition of the eigenvalues\ \eqref{eq:GenEig}, we obtain the SDE governing the evolution of the eigenmodes $\widehat{\V{X}}=\M{V}^{-1}\D \V{X}$, 
\begin{align}
\nonumber
\M{V} \frac{d \Xhat}{dt} &=- \widehat{\M{N}}_0\M{L}\M{V} \Xhat +\sqrt{2 k_B T}\widehat{\M{N}}_0 \M{L} \M{V}\M{\Lambda}^{1/2}  \Lop{W}\\ 
\label{eq:LangDiag}
 \frac{d \Xhat}{dt} &=- \M{\Lambda}^{-1} \Xhat +\sqrt{2 k_B T} \M{\Lambda}^{-1/2} \Lop{W}.
\end{align}
Thus, each of the $2N+3$ inextensible eigenmodes relaxes with characteristic timescale equal to its eigenvalue, and the modes satisfy the equipartition principle $\Exp{\Xhat \Xhat^T}=\left(k_B T\right) \M{I}$. The average elastic energy is then $\frac{1}{2}\Exp{\D \V{X}^T \M{L}\D \V{X}}=\frac{1}{2}\Exp{\Xhat^T \M{V}^T \M{L}\M{V} \Xhat}=\left(N+3/2\right)k_B T $.

\begin{figure}
\centering
\includegraphics[width=\textwidth]{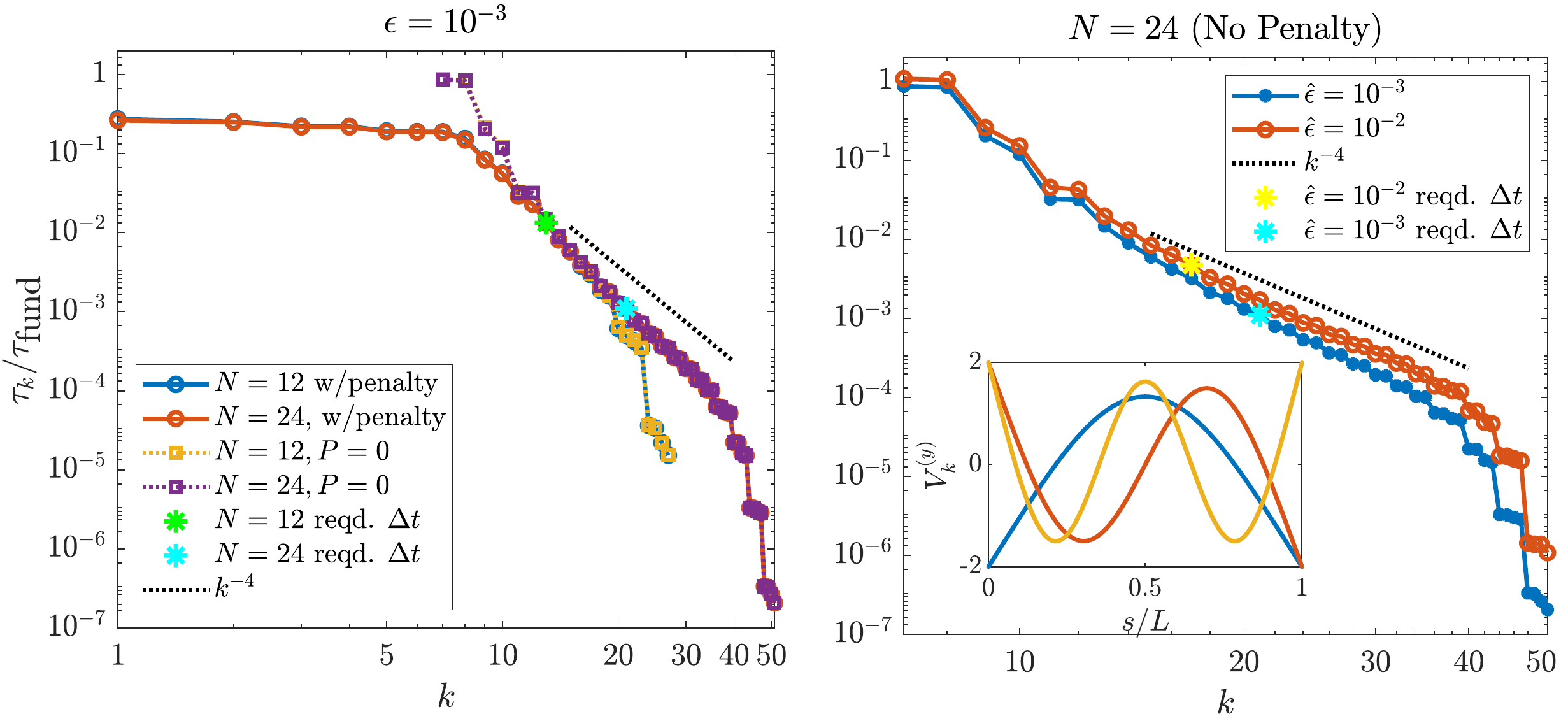}
\caption{\label{fig:Timescales} Timescales in the linearized filament problem. In each plot, we give the timescales (the $2N+3$ nonzero diagonal entries of $\M{\Lambda}$ in\ \eqref{eq:GenEig}) made dimensionless by $\tau_\text{fund}$ in\ \eqref{eq:TauSm} (the largest timescale for free fibers), then mark with a star the time step size required for the midpoint temporal integrator to have an error less than 5\% in the equilibrium variance when $\ell_p/L=10$ (see Section\ \ref{sec:LinTInt}). (Left) Dependence of timescales on $N$. We show the timescales for two different $N$ with (blue and red) and without (yellow and purple) the penalty force. (Right) Dependence on the mobility (aspect ratio $\epsRS$) with $P=0$. We show $\epsRS=10^{-3}$ in blue and $\epsRS=10^{-2}$ in red. Inset: the first three unique modes for a straight filament (compare to \cite[Fig.~2(a)]{kantsler2012fluctuations}). The eigenvalue of the blue mode is the timescale $\tau_\text{fund}$.}
\end{figure}

We first examine the dependence of the timescales on $N$ in the left panel of Fig.\ \ref{fig:Timescales}, beginning by plotting the timescales for $N=12$ (blue) and $N=24$ (red) with the penalty force included in $\M{L}$. Including the penalty force makes $\M{L}$ invertible, which gives finite timescales for the first six modes (three translation and three rotation) that have zero eigenvalues when $P=0$. After the first few modes, the timescales begin to decay with the expected $k^{-4}$ scaling \cite{poelert2012analytical} until the spatial discretization error causes shorter timescales than we would obtain in continuum. We see in particular that for a given $N$ we have about 6 incorrect timescales from spatial discretization error. Excepting these modes, we see that increasing $N$ only adds additional small timescales into the problem, since $N=12$ is sufficient to correctly give the timescales of the first 15 or so modes.  Since the penalty force contributes a fixed multiple of the identity to $\M{L}$, the largest eigenvalues are changed very little, and so the effect of the penalty force on the relaxation timescales is negligible for modes beyond the first 10 or so (see the left panel of Fig.\ \ref{fig:Timescales}). 

\subsection{Temporal accuracy \label{sec:LinTInt}}
\begin{figure}
\centering
\includegraphics[width=\textwidth]{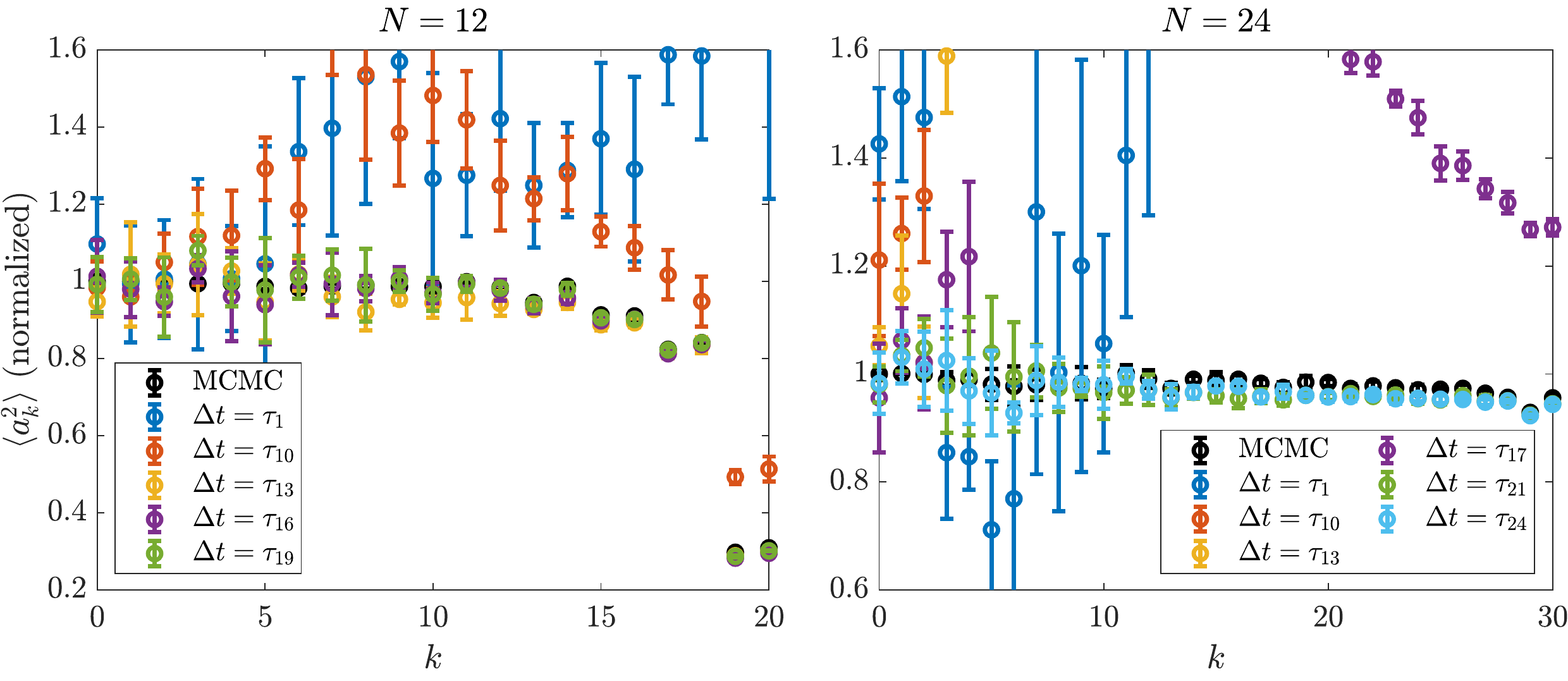}
\includegraphics[width=\textwidth]{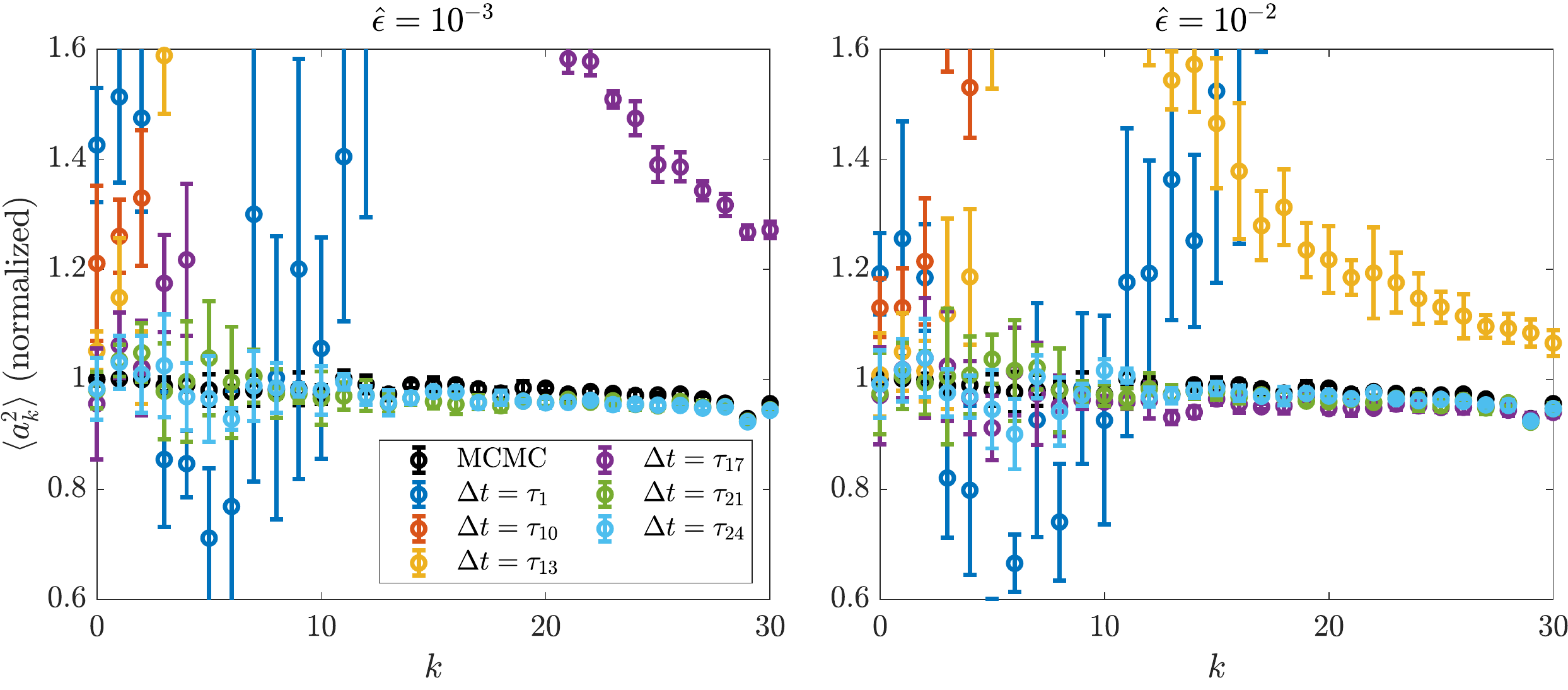}
\caption{\label{fig:TintLinConv} Convergence of the equilibrium covariance in our midpoint temporal integrator, with $\D t$ expressed in terms of relaxation timescales. In each case, we plot convergence to the MCMC results (see Fig.\ \ref{fig:MCMCPen}) as a function of the time step size, which is chosen from the modes in Fig.\ \ref{fig:Timescales}. (Top) Dependence of timescales on $N$ with fixed $\ell_p/L=10$ and $\epsRS=10^{-3}$. (Bottom) Dependence on the mobility (aspect ratio $\epsRS$) with fixed $N=24$ and $\ell_p/L=10$.}
\end{figure}

\begin{figure}
\centering
\includegraphics[width=\textwidth]{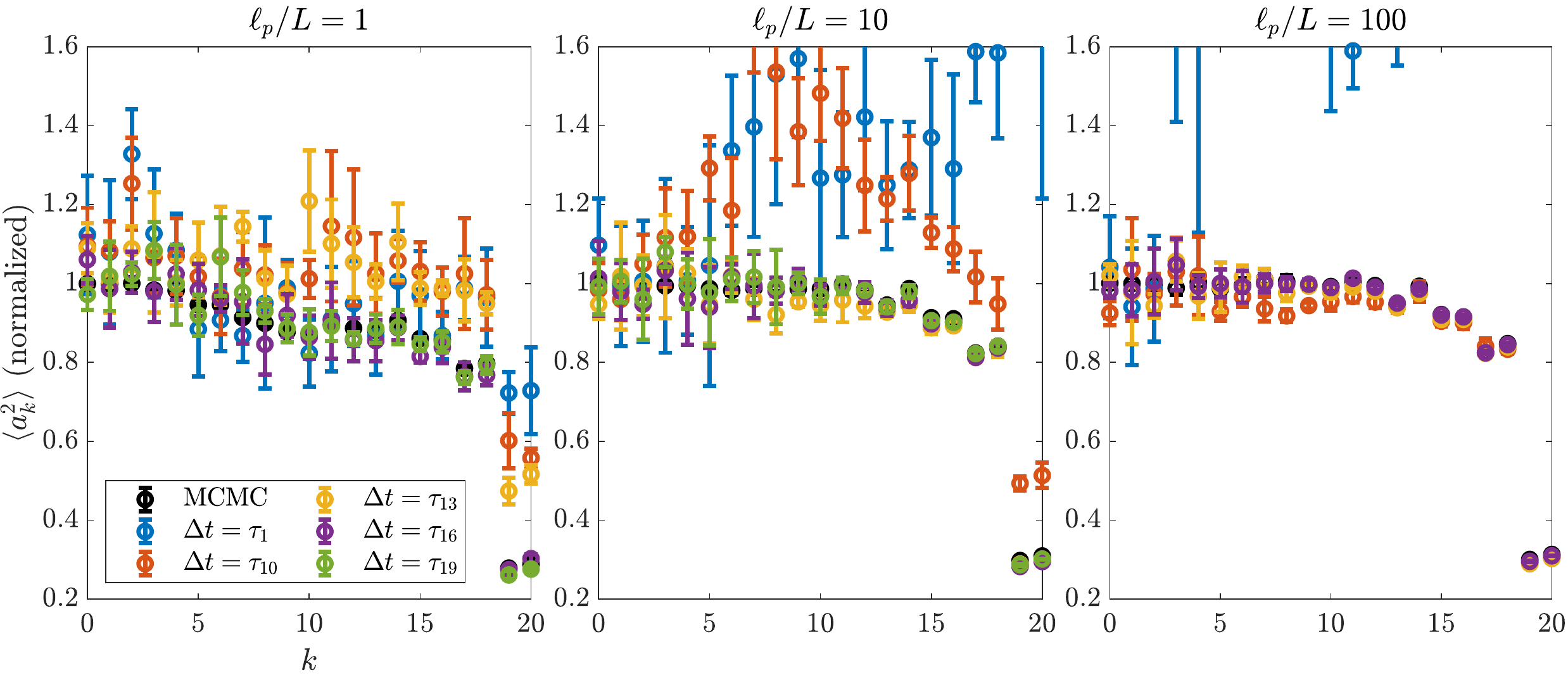}
\caption{\label{fig:TintLinConvLp} Convergence of the equilibrium covariance in our midpoint temporal integrator for changing $\kappa$, with $\D t$ expressed in terms of relaxation timescales. We fix $N=12$ and $\epsRS=10^{-3}$ and plot the results for changing $\kappa$ in terms of $\ell_p/L$. }
\end{figure}

We now use the timescale analysis of the previous section to understand the performance of our temporal integrator. For each set of parameters, we first compute the relevant timescales using the method of Section\ \ref{sec:LinTScales}. Then, we choose a subset of the modes and run the Langevin dynamics with $\D t=\tau_k$, for each $k$ in the set of modes chosen. To obtain statistics, we initialize the fiber in its equilibrium state and run until $10\tau_1$, removing the first $\tau_1$ as a burn-in period. We repeat this a total of 20 times to generate a mean, then perform the 20 trials a total of 5 times to generate error bars. The resulting variance of each mode is shown in Figs.\ \ref{fig:TintLinConv} and\ \ref{fig:TintLinConvLp}. Here we follow the pattern of Fig.\ \ref{fig:Timescales} to systematically vary the parameters $N$ (top of Fig.\ \ref{fig:TintLinConv}), $\epsRS$ (bottom), and $\kappa$ (Fig.\ \ref{fig:TintLinConvLp}), while keeping the others constant. 

\subsubsection{Effect of changing $N$}
Beginning with the top of Fig.\ \ref{fig:TintLinConv}, we see that the number of modes whose dynamics we need to resolve with $\D t$ depends on $N$. For $N=12$, where there are a total of 27 modes, we can use a time step corresponding to the $k=13$ mode and obtain an accurate variance for most modes, while for $N=24$, where there are 51 modes, we need a time step which matches the timescale for the $k=21$ mode. Thus the number of modes we need to resolve increases with $N$. Since the timescale of the modes scales as $k^{-4}$, doubling $N$ requires a time step refinement of at least 10. This is obviously sub-optimal, and tells us that the optimal way to simulate the dynamics of the first $k$ modes is to choose the minimum $N$ that resolves those $k$ modes.

\subsubsection{Effect of changing $\epsRS$}
We now move on to the variation in $\Delta t$ with changing $\epsRS$ for $N=24$. The bottom row of Fig.\ \ref{fig:TintLinConv} shows that, for a fixed relative time step size, $\epsRS=10^{-2}$ gives results closer to the equilibrium distribution. The required time step size for a 5\% error in the variance is about $\D t = \tau_{17}$ for $\epsRS=10^{-2}$ and $\D t = \tau_{21}$ for $\epsRS=10^{-3}$. Because the scaling of the timescales for intermediate modes with $\epsRS$ is somewhere between log scaling (factor of 1.5) and $\epsRS^{-1}$ scaling (factor of 10), the net result of this is that the absolute time step size we need decreases by a factor of six (see the stars in Fig.\ \ref{fig:Timescales}(b)). That said, changing aspect ratio cannot be accommodated by a simple rescaling of time since it scales each mode differently depending on its smoothness.

\subsubsection{Effect of changing $\kappa$}
Finally, we consider the variation with $\ell_p/L$ in Fig.\ \ref{fig:TintLinConvLp}, where we observe the expected behavior: as we increase $\ell_p/L$, we need to resolve the dynamics of fewer modes to get the correct variance. However, we find that the number of modes required does \emph{not} scale as $1/\kappa$; rather, we need $\approx 3$ less modes for every factor of 10 increase in $\kappa$ since the chain becomes smoother. Because increasing $\kappa$ causes the timescale of each mode to decrease, the net effect of this is that the actual time step size we need is roughly constant as $\kappa$ increases (for the range of parameters we consider here).

\subsection{Importance of drift terms\label{sec:PoorManResults}}
\begin{figure}
\centering
\includegraphics[width=\textwidth]{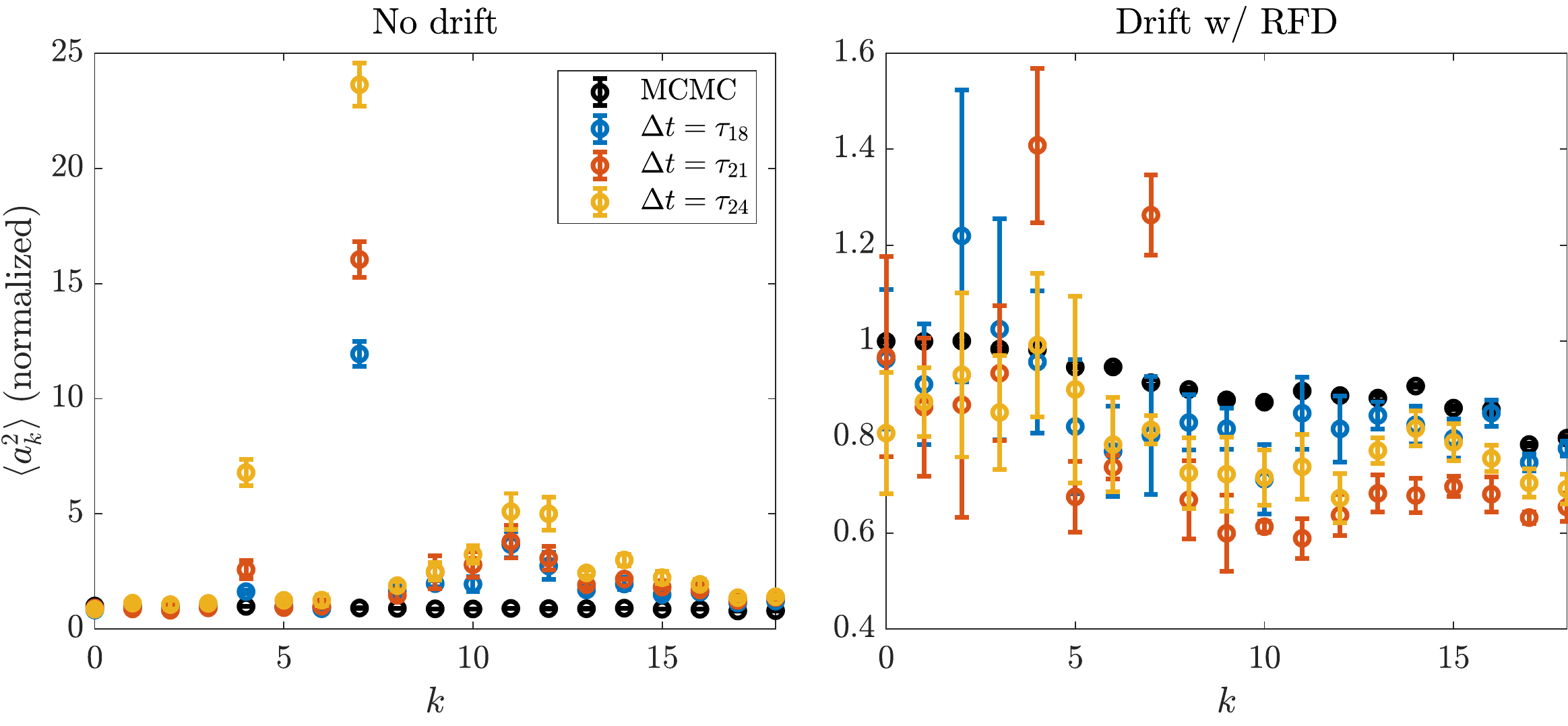}
\caption{\label{fig:PoorMan} Results for linearized fluctuations with constant mobility $\M{M}=\M{I}/(8\pi \mu)$. We fix $N=12$ and $\ell_p/L=1$ and show the variance of each of the $L^2$-orthonormal modes of the linearized covariance matrix\ \eqref{eq:CovExp} (as described in Appendix\ \ref{sec:CovExp}) using the simple RFD scheme in Section\ \ref{sec:NaiveRFD}. We show the results without (left) and with (right) the additional RFD drift term in\ \eqref{eq:OmegaRFD} with $\delta = 10^{-5}$. In all cases, the MCMC results in black are those obtained in Fig.\ \ref{fig:MCMCPen}(b). }
\end{figure}

To conclude, we look at what happens to the variance when we exclude the drift terms from the overdamped Langevin equation\ \eqref{eq:ItoX}. We do this with the RFD scheme of Appendix\ \ref{sec:NaiveRFD}, which allows us to explicitly exclude the drift terms to see what kind of covariance we obtain in $\V{X}$. To demonstrate that the drift terms are important even when the mobility is position-independent, we use the constant mobility $\M{M}=\M{I}/(8\pi \mu)$, which we apply to a fiber with $N=12$ tangent vectors and $\ell_p/L=1$. In Fig.\ \ref{fig:PoorMan}, we show the variance of each of the $L^2$ orthonormal modes with several different values of $\D t$, compared to the results from MCMC, which we obtained on the \emph{spectral} grid in Section\ \ref{sec:MCMCSpec}. In the left panel, we show the results when we simulate the dynamics \emph{without} adding the RFD term\ \eqref{eq:OmegaRFD}. As $\D t \rightarrow 0$, many of the modes (but especially 4 and 7) increase in variance without a visible bound. When we add the RFD term back in, the right panel shows that the dynamics converge to the MCMC dynamics as $\D t \rightarrow 0$. This shows that, even if the mobility is constant, proper handling of the stochastic drift terms that arise from the inextensibility constraints is necessary to obtain the correct dynamics. We note also that the RFD scheme appears to require \emph{all} of the modes (there are 27 for $N=12$) to be resolved to obtain accurate results, which is in contrast to the midpoint scheme which requires less than half.

\end{appendices}

\end{document}